\documentclass[10pt,reqno]{amsart}

\usepackage[latin1]{inputenc}
\usepackage{graphics,color}
\usepackage{amssymb,amsmath,amsthm,amscd}
\usepackage{latexsym,verbatim,graphicx,amsfonts}
\usepackage{hyperref}
\usepackage[mathscr]{euscript}
\usepackage{amsmath, amsthm, amssymb}
\usepackage{dsfont}
\usepackage{enumitem}
\usepackage{mathrsfs}
\usepackage{float}

\theoremstyle{plain}
\newtheorem{theorem}{Theorem}[section]

\newtheorem{corollary}[theorem]{Corollary}
\newtheorem{lemma}[theorem]{Lemma}
\newtheorem{proposition}[theorem]{Proposition}
\theoremstyle{definition}

\newtheorem{definition}[theorem]{Definition}
\newtheorem{example}[theorem]{Example}
\theoremstyle{remark}
\newtheorem{rem}[theorem]{Remark}

\newtheorem{question}[theorem]{Question}

\DeclareMathOperator{\spn}{span}

\newcommand{\norm}[1]
{\lVert#1\rVert}

\DeclareMathOperator*{\esssup}{ess\,sup}

\newcommand{\vertiii}[1]{{\left\vert\kern-0.25ex\left\vert\kern-0.25ex\left\vert #1 
    \right\vert\kern-0.25ex\right\vert\kern-0.25ex\right\vert}}

\font\sstext=ecss1000
\font\sssub=ecss1000 at 7pt
\font\sssubsub=ecss1000 at 5pt

\newfam\ssfam
\textfont\ssfam=\sstext
\scriptfont\ssfam=\sssub
\scriptscriptfont\ssfam=\sssubsub



\begin{document}

\title[Quantifying (non-)weak compactness]{Quantifying (non-)weak compactness of operators on $AL$- and $C(K)$-spaces}
\author[A. Acuaviva]{Antonio Acuaviva}
\address{School of Mathematical Sciences,
Fylde College,
Lancaster University,
LA1 4YF,
United Kingdom} \email{ahacua@gmail.com}
\author[A. B. Nasseri]{Amir Bahman Nasseri}
\address{Independent researcher} \email{amirbahman@hotmail.de}

\date{\today}

\begin{abstract}

We study the representation of non-weakly compact operators between $AL$-spaces. In this setting, we show that every operator admits a best approximant in the ideal of weakly compact operators. Using duality arguments, we extend this result to operators between $C(L)$-spaces where $L$ is extremally disconnected. We also characterize the weak essential norm for operators between $AL$-spaces in terms of factorizations of the identity on $\ell_1$. As a consequence, we deduce that the weak Calkin algebra $\mathscr{B}(E)/\mathscr{W}(E)$ admits a unique algebra norm for every $AL$-space $E$. By duality, similar results are obtained for $C(K)$-spaces. In particular, we prove that for operators $T: L_{\infty}[0,1] \to L_{\infty}[0,1]$ the weak essential norm, the residuum norm, and the De Blasi measure of weak compactness coincide, answering a question of Gonz\'alez, Saksman and Tylli.

\end{abstract}

\maketitle

\section{Introduction and Main Results}\label{Introduction}

A Banach lattice $E$ is said to be an Abstract Lebesgue space ($AL$-space) if it behaves like the classical $L_1$-spaces, in other words, if $\norm{x + y} = \norm{x} + \norm{y}$ for every $x, y \in E^+$, the positive cone of $E$. Important examples of $AL$-spaces naturally include $L_1$-spaces and the spaces $\mathcal{M}(K)$ of Borel regular measures on compact Hausdorff spaces $K$. In fact, from the remarkable representation theorems of Kakutani (see, for example \cite[Theorem 4.27]{aliprantis2006}), every $AL$-space can be represented (up to a lattice isometry) by an $L_1$-space, where the underlying measure space is \emph{strictly localizable } (see Definition \ref{def: strictly-localizable}). Thus, the study of $AL$-spaces reduces to the study of $L_1$-spaces with an underlying strictly localizable measure space. The same holds in the complex setting: using complex Banach lattices, every complex $AL$-space admits a representation as an $L_1$-space (see \cite[p. 135, Theorem 3]{lacey2008isometric}) with underlying strictly localizable measure space. In particular, for any measure space $(X, \Sigma, \mu)$, since $L_1(X, \Sigma, \mu)$ is an $AL$-space, there exists a strictly localizable measure space $(X', \Sigma', \mu')$ such that $L_1(X, \Sigma, \mu)$ and $L_1(X', \Sigma', \mu')$ are lattice isometrically isomorphic; for this special case, the construction of this lattice isometry can be explicitly found, for example, in the proof of \cite[Lemma 7.4]{kacena2013quantitative}. Consequently, requiring the measure space to be strictly localizable entails no loss of generality, and our results therefore remain valid for any measure space $(X,\Sigma,\mu)$. \\

The study of weak compactness in $AL$-spaces is classical in Banach space theory. Early works include those of Dunford and Pettis \cite{dunford1940linear} and Grothendieck's characterization of weak compactness in $\mathcal{M}(K)$ spaces \cite{grothendieck1953applications}. In recent years, there has been growing interest in the quantitative study of weak compactness (see, for example, \cite{angosto2009measures, kacena2013quantitative, kalenda2025quantitative}), leading to quantitative versions of classical results with consequences that cannot be directly deduced from the original statements. It is within this context that we situate our work.

Let us recall that an operator $T: E \to F$ is weakly compact if $T[B_E]$, the image of the unit ball of $E$, is a relatively weakly compact set in $F$. We define the \emph{weak essential norm} of an operator $T:E \to F$ as the distance of the operator to the ideal $\mathscr{W}(E, F)$ of weakly compact operators, that is
\begin{equation*}
    \norm{T}_w = \inf \{\norm{T - W}: W \in \mathscr{W}(E,F)\}.
\end{equation*}

We study the representation of weakly compact operators between $AL$-spaces and use this to give quantitative estimates of the weak essential norm. Our main results in this regard can be summarized in the following theorem.

\begin{theorem}\label{th: MainTheorem}
    Let $(X, \Sigma, \mu)$ and $(Y, \Gamma, \nu)$ be strictly localizable measure spaces, and $T: L_1(X, \Sigma, \mu) \to L_1(Y, \Gamma, \nu)$ be an operator. Then there exists a weakly compact operator $W:  L_1(X, \Sigma, \mu) \to L_1(Y, \Gamma, \nu)$ such that
    \begin{equation*}
        \norm{T}_w = \norm{T - W} = \inf_{\Theta\subseteq Y, \nu(\Theta)<\infty} \alpha_{\Theta}(T),
    \end{equation*}
    where $\alpha_\Theta (T) = \limsup_{\nu(B) \to 0, B \subseteq \Theta} \norm{\chi_{Y \backslash \Theta} T + \chi_B T}$.
\end{theorem}

Recall that if $K$ is a compact Hausdorff space, then the Riesz representation theorem gives $C(K)^* = \mathcal{M}(K)$. A simple duality argument now allows us to obtain the following. 

\begin{corollary}\label{cor: bestapproximationC(K)}
    Let $T: C(K) \to C(L)$ be an operator, where $K,L$ are compact Hausdorff spaces and $L$ is extremally disconnected. Then there exists a weakly compact operator $W: C(K) \to C(L)$ such that
    \begin{equation*}
         \norm{T^*}_w = \norm{T}_w = \norm{T - W}.
    \end{equation*}
\end{corollary}

Observe that the previous corollary applies, in particular, to $L_\infty[0,1]$.

\begin{rem}
    The results above fail if one considers the ideal of compact operators instead of the weakly compact operators. In this case, Feder \cite{feder1980certain} has shown that, in general, no best approximant exists.
\end{rem}

The previous theorem extends the main result of \cite{weis1984approximation}, which was originally proved under the assumption that \(X\) and \(Y\) are compact metric spaces, the \(\sigma\)-algebras are the Borel \(\sigma\)-algebras, and the measures are finite. We first generalize this result to arbitrary finite measure spaces, which in turn allows us to extend it further to all strictly localizable measure spaces. This broader extension has many interesting consequences for the study of non-weakly compact operators on $AL$- and $C(K)$-spaces, as we shall now discuss. \\

We start by giving a characterisation of the weak essential norm of operators between $AL$-spaces, in terms of factorizations of the identity operator on $\ell_1$. We will always use the convention that the supremum of positive quantities over the empty set is zero.

\begin{theorem}\label{th: AL-factorization}
    Let $E, F$ be $AL$-spaces and $T: E \to F$ be an operator. Then
    \begin{equation*}
        \norm{T}_w = \sup \{ \norm{U}^{-1}\norm{V}^{-1}: UTV = I_{\ell_1}\}.
    \end{equation*}
\end{theorem}

The previous theorem may be viewed as a quantitative version of Pe{\l}czy{\'n}ski's theorem on non-weakly compact operators with range in an $AL$-space \cite[Theorem 1]{pelczynski1965strictly2} in the special case where the domain is also an $AL$-space. This restriction is necessary, as no such quantitative statement involving the weak essential norm holds in general; see Example~\ref{ex: Not-main-L}. \\

By duality, a similar result can be obtained for $C(K)$-spaces, which again provides a quantitative version of  Pe{\l}czy{\'n}ski's classical theorem on weakly compact operators from $C(K)$-spaces \cite[Theorem 1]{pelczynski1965strictly}. 
Recall that the \emph{minimum modulus} (also called \emph{injectivity modulus}) of an operator $T: E \to F$ is defined by
\begin{equation*}
    m(T) = \inf \{ \norm{Tx}: x \in E, \hspace{5pt} \norm{x} = 1\}.
\end{equation*}

\begin{theorem}\label{th: C(K)-factorization}
    Let $T: C(L) \to C(K)$ be an operator, where $K,L$ are compact Hausdorff spaces. Then
    \begin{equation*}
        \norm{T^*}_w = \sup \{m(TV)|\hspace{5pt}  V: c_0 \to C(L), \hspace{5pt} \norm{V} \leq 1 \}.
    \end{equation*}
    In particular, if $K$ is sequentially compact, we have
    \begin{equation*}
        \norm{T^*}_w = \sup \{ \norm{U}^{-1}\norm{V}^{-1}: UTV = I_{c_0}\},
    \end{equation*}
    while if $L$ is extremally disconnected
    \begin{equation*}
        \norm{T^*}_w = \sup \{ \norm{U}^{-1}\norm{V}^{-1}: UTV = I_{\ell_\infty}\}.
    \end{equation*}
\end{theorem}

\begin{rem}
    Observe that if $K$ is assumed to be extremally disconnected, then $\norm{T^*}_w = \norm{T}_w$.
\end{rem}

Factorizations of this type have appeared in the literature for a significant time (see, for instance, \cite{meyer1995lower} and the discussion following \cite[Theorem 2.5]{Tylli1}, as well as \cite[Lemma 2.1]{Tylli1}). More recently, they have played a prominent role in the work of Johnson and Phillips \cite{johnsonandphillips}; see also the work of Arnott and Laustsen \cite{arnott2023uniqueness}. \\

Once a quantitative factorization has been established, one can extract significant information about the structure of norms in the weak Calkin algebra. Recall that an algebra norm $\norm{\cdot}$ on an algebra $\mathcal{A}$ is \emph{minimal} (respectively, \emph{maximal}) if for any other algebra norm $\vertiii{\cdot}$ on $\mathcal{A}$, there exists a constant $c > 0$ such that $c \norm{\cdot} \leq \vertiii{\cdot}$ (respectively, there exists $C > 0$ such that $\vertiii{\cdot} \leq C \norm{\cdot}$).

We say that an algebra $\mathcal{A}$ \emph{admits a unique} algebra norm if it possesses a norm that is both maximal and minimal. In this case, all algebra norms on $\mathcal{A}$ induce the same topology; that is, there is a unique topology arising from algebra norms. We emphasize that completeness of the algebra norms is not assumed. Using our factorization results, we obtain the following.

\begin{corollary}\label{cor: unique-algebra-normAL-space}
    Let $E$ be an $AL$-space. Then the weak Calkin algebra $\mathscr{B}(E) / \mathscr{W}(E)$ admits a unique algebra norm.
\end{corollary}

\begin{corollary}\label{cor: algebra-normC(K)-space-extremally}
    Let $L$ be an extremally disconnected compact Hausdorff space. Then the weak essential norm, $\norm{\cdot}_w$, is minimal on the weak Calkin algebra \break $\mathscr{B}(C(L)) / \mathscr{W}(C(L))$. 
\end{corollary}

\begin{corollary}\label{cor: algebra-normC(K)-space-sequentially}
    Let $K$ be a sequentially compact and compact Hausdorff space. Then the dual weak essential norm, $T \mapsto \norm{T^*}_w$, is minimal on the weak Calkin algebra $\mathscr{B}(C(K)) / \mathscr{W}(C(K))$. 
\end{corollary}

We note that Corollary \ref{cor: algebra-normC(K)-space-sequentially} gives a partial answer to a problem posed by Tylli, see \cite[Problems]{Tylli1}, who asked whether the weak essential norm is minimal on the weak Calkin algebra $\mathscr{B}(C([0,1])) / \mathscr{W}(C([0,1]))$. This naturally leads to the following question.

\begin{question}
   Is the weak essential norm dual for the weak Calkin algebra \break $\mathscr{B}(C[0,1])/\mathscr{W}(C[0,1])$, in other words, does there exists $c > 0$ such that $\norm{T^*}_w \geq c \norm{T}_w$ for all $T \in \mathscr{B}(C([0,1]))$?
\end{question}

Another property of interest in the study of algebra norms, this time of an isometric nature, is \emph{Bonsall's minimality property} (see \cite{Tylli1}). A norm $\norm{\cdot}$ on an algebra $\mathcal{A}$ is said to have Bonsall's minimality property if, for any other algebra norm $\vertiii{\cdot}$ satisfying $\vertiii{\cdot} \leq \norm{\cdot}$, we have equality $\vertiii{\cdot} = \norm{\cdot}$. We obtain the following.

\begin{corollary}\label{cor: BonsalMinimalityAL}
    Let $E$ be an $AL$-space. Then the weak essential norm has Bonsall's minimality property on the weak Calkin algebra $\mathscr{B}(E) / \mathscr{W}(E)$.
\end{corollary}

In passing, we note that this corollary is a consequence of a more general result about measures of weak noncompactness in $AL$-spaces, see Proposition \ref{prop: normalised-measures-weak-non} for details. 

\begin{corollary}\label{cor: BonsalMinimality-extremally}
    Let $L$ be an extremally disconnected compact Hausdorff space. Then the weak essential norm has Bonsall's minimality property on the weak Calkin algebra $\mathscr{B}(C(L)) / \mathscr{W}(C(L))$. 
\end{corollary}

\begin{corollary}\label{cor: BonsalMinimality-sequentially}
    Let $K$ be a sequentially compact and compact Hausdorff space. Then the dual weak essential norm, $T \mapsto \norm{T^*}_w$, has Bonsall's minimality property on the weak Calkin algebra $\mathscr{B}(C(K)) / \mathscr{W}(C(K))$. 
\end{corollary}

We observe that the same results about general measures of weak noncompactness can also be established in the case of $C(L)$-spaces, where $L$ is an extremally disconnected compact Hausdorff space (and dual ones in the case when $K$ is sequentially compact). \\

Of particular interest, the previous results give information about the \emph{De Blasi measure of weak noncompactness} and the \emph{residuum norm}. Let us recall that, for a subset $A$ of a Banach space $E$, we define the De Blasi measure of weak noncompactness of $A$ by
\begin{equation}\label{eq: deBlasi}
    \omega(A) = \inf \{ \varepsilon > 0: A \subseteq D + \varepsilon B_E, D \subseteq E \text{ weakly compact} \}.
\end{equation}
The function $\omega(\cdot)$ was introduced by De Blasi \cite{de1977property}, who showed that $\omega(A) = 0$ if and only if $A$ is relatively weakly compact. For an operator $T: E \to F$, we naturally define $\omega(T) = \omega(TB_E)$. 

For a Banach space $E$, denote by $E^{co}$ the quotient space $E^{**} / E$. The \emph{residuum operator} of the operator $T: E \to F$ is the operator $T^{co}: E^{co} \to F^{co}$ defined by
    \begin{equation*}
        T^{co}(x + E) := T^{**}(x) + F,
    \end{equation*}
where $x \in F^{**}$. In the literature, this operator has also been denoted by $R(T)$. It satisfies $\norm{T^{co}} \leq \norm{T}$, and $T^{co} = 0$ if and only if $T$ is weakly compact. Residuum operators have been extensively studied and play an important role in the theory of Tauberian operators; see, for instance, \cite[Chapter 3]{gonzalez2010tauberian}.

As an immediate consequence of Corollary \ref{cor: BonsalMinimality-extremally}, we obtain information about these quantities, answering questions of Tilly \cite[Remark 2.2, (ii)]{tylli1993spectral} and Gonz\'alez, Saksman and Tylli \cite[Remark after the proof of Theorem 2.2]{Gonzalez1995}.

\begin{theorem}\label{th: Residum-and-w}
    Let $T: L_\infty[0,1] \to  L_\infty[0,1]$ be an operator. Then
    \begin{equation*}
        \norm{T}_w = \norm{T^*}_w = \norm{T^{co}} = \omega(T).
    \end{equation*}
\end{theorem}

Before concluding the introduction, we would like to remark that, by the aforementioned theorems of Pe{\l}czy{\'n}ski, the ideals of weakly compact, strictly singular, and strictly cosingular operators coincide for $AL$-spaces. Consequently, any of our results concerning measures of weak noncompactness can be applied to measures of strict singularity or strict cosingularity. In the case of $C(K)$-spaces, the ideals of weakly compact and strictly singular operators coincide, and thus the same can be said for measures of strict singularity.

\bigskip
\section{Organization and notation}

We follow standard notational conventions. By an \emph{operator}, we mean a bounded linear map between normed spaces. For Banach spaces $E$ and $F$, we denote by $\mathscr{B}(E, F)$ the space of bounded operators from $E$ to $F$, and by $\mathscr{W}(E, F)$ the ideal of weakly compact operators. When $E = F$, we write simply $\mathscr{B}(E)$ and $\mathscr{W}(E)$. We write $\norm{T}_w$ for the weak essential norm of an operator $T \in \mathscr{B}(E,F)$, that is
\begin{equation*}
    \norm{T}_w = \norm{T + \mathscr{W}(E,F)} = \inf \{\norm{T - W}: W \in \mathscr{W}(E,F)\}.
\end{equation*}
We will denote by $B_E$ the unit ball of the Banach space $E$. We also denote by $I_E$ the identity operator on $E$; when the context is clear, we simply write $I$.

We denote by $(X, \Sigma, \mu)$ a measure space, where $X$ is a set, $\Sigma$ is a $\sigma$-algebra on $X$, and $\mu$ is a measure defined on $\Sigma$. As usual, $L_1(X, \Sigma, \mu)$ denotes the Banach space of equivalence classes of $\mu$-integrable functions.

For $A \in \Sigma$, we denote by $\mathds{1}_A$ the indicator function of $A$. We write $\chi_A$ for the operator
\begin{equation*}
    \chi_A: L_1(X, \Sigma, \mu) \to L_1(A, \tilde{\Sigma}, \mu), \hspace{10pt} \chi_A(f) = \mathds{1}_A f,
\end{equation*}
where $\tilde{\Sigma} = \{B \cap A : B \in \Sigma\}$.

Occasionally, we use $:=$ or $=:$ to stress that a term is being defined; for example, we write $L_1(X, \Sigma, \mu) =: E$ to emphasize that $E$ is defined by the expression on the left. More specialized notation will be introduced as and when needed. \\

We give an outline of the structure of this paper. In Section \ref{sec: preliminaries} we introduce some background material and associated terminology. Specifically, we introduce notation and known results about Riesz spaces, measure theory, integral and singular operators, and weak compactness in finite measure spaces. Section \ref{sec: mainTheorem} is devoted to the proof of Theorem \ref{th: MainTheorem} and also contains the proof of Corollary \ref{cor: bestapproximationC(K)}. For notational convenience, all results in this section are stated exclusively for the real case; however, the reader may verify that they extend easily to the complex case with the natural modifications.

Section \ref{sec: weak-calkinAL} discusses the consequences of Theorem \ref{th: main-with-Theta} for the weak Calkin algebra of $AL$-spaces. This section includes the proofs of Theorem \ref{th: AL-factorization} and Corollaries \ref{cor: unique-algebra-normAL-space} and \ref{cor: BonsalMinimalityAL}, as well as a brief aside on a general notion of measures of weak noncompactness for operators between $AL$-spaces, which may be of independent interest.

Finally, Section \ref{sec: weak-calkinC(K)} addresses the consequences of Theorem \ref{th: MainTheorem} for the weak Calkin algebra of $C(K)$-spaces. It contains the proofs of Theorem \ref{th: C(K)-factorization} and Corollaries \ref{cor: algebra-normC(K)-space-extremally}, \ref{cor: algebra-normC(K)-space-sequentially}, \ref{cor: BonsalMinimality-extremally}, and \ref{cor: BonsalMinimality-sequentially}.

\bigskip
\section{Preliminaries}\label{sec: preliminaries}

\subsection{Aspects about Riesz spaces.}\label{subsec: aspect-about-Riesz} 

We start with a very brief overview of normed Riesz spaces and Banach lattices; for a more detailed introduction, we refer the reader to \cite{abramovich2002}, \cite{aliprantis2006} and \cite{meyer-nieberg1991}.

An ordered vector space $E$ is called a \emph{Riesz space} if for any pair of elements $x, y \in E$, its supremum exists. We denote by $x\vee  y :=\sup\{x,y\}$ and $x\wedge y :=\inf\{x,y\}$. A Riesz space $E$ is called \emph{Dedekind complete} if every bounded set $A \subseteq E$ admits a supremum.

A Riesz space is said to have the \emph{countable sup property} if whenever an arbitrary subset $A$ has a supremum, then there exists an at most countable subset $C\subseteq A$ such that $\sup(C)=\sup(A)$. We note that we have the following, see \cite[comment after Lemma 2.6.1]{meyer-nieberg1991}.

\begin{proposition}\label{prop: countable-supremum-property}
    Let $(X, \Sigma, \mu)$ be a $\sigma$-finite measurable space. Then $L_1(X, \Sigma, \mu)$ has the countable sup property.
\end{proposition}

Given an element $x$ of a Riesz space $E$, we denote by $|x| := x \vee (-x)$ the \emph{modulus of} $x$. Two elements $x, y \in E$ are called \emph{orthogonal}, and denoted $x \perp y$, if $|x| \wedge |y| = 0$. This is equivalent to saying that $|x - y| = |x+y|$.
For $x,y \in E$ with $x \leq y$ then the \emph{order interval} $[x,y]$ is defined as $[x,y]:=\{x \in E: x\leq z\leq y\}$.
A subspace $A \subseteq E$ is called an \emph{order ideal} if for $y \in A$ and $|x| \leq |y|$  we have $x \in A$, or equivalently the order interval $[0, |y|] \subseteq A$ for $y \in A$. Let $\{ x_{\alpha} \}$ be a net in a Riesz space. A subset $A$ of a Riesz space is said to be \emph{order closed}, if for an order convergent net $\{x_\alpha\}\subseteq A$ with limit $x$ (in symbols $x_{\alpha}\overset{o}{\rightarrow}x$ ) we have $x \in A$. For the definition of order convergence, we refer the reader to \cite[p. 33]{aliprantis2006}.

An order closed ideal is called a \emph{band}. The disjoint complement of a subset $A\subseteq E$, denoted by $A^{d}$, is defined by
\begin{equation*}
    A^{d}=\{x\in E\ : \ |x|\wedge |y|=0 \text{ for each  }y\in A\}.
\end{equation*}
A band $A$ is called a \emph{projection band} if $A\oplus A^{d}=E$ holds. The corresponding projection to this decomposition is called \emph{band projection}. In a Dedekind complete Riesz space, every band is a projection band, see \cite[Theorem 1.42]{aliprantis2006}.

If $E$ is a Riesz space, then the set of all positive elements is denoted by $E^{+}:=\{x\in E\ |\ x\geq0\}$, and is called the \emph{positive cone}. If $T$ is a linear map between two ordered vector spaces, then $T$ is called \emph{positive} if $T(x) \geq 0$ for all $x\in E^{+}$. The vector space of all linear maps from $E$ to $F$, where $F$ is another ordered vector space, will be denoted by $\mathcal{L}(E, F)$, with the order $T\geq S$ whenever $T(x)\geq S(x)$ for all $x\in E^{+}$.\\

A linear map $T \in \mathcal{L}(E, F)$ is called \emph{regular} if it can be written as the difference of two positive linear maps. If $T$ is regular and $F$ is Dedekind complete then the modulus of $|T|$ exists \cite[Theorem 1.18]{aliprantis2006} and $T=T^{+}-T^{-}$, where $T^{+}$ (called the \emph{positive part}) and $T^{-}$ (called the \emph{negative part}) are positive linear maps. We denote by $\mathcal{L}_{r}(E,F)$ the set of all regular linear maps. If $F$ is Dedekind complete then $\mathcal{L}_{r}(E,F)$ is also Dedekind complete, see \cite[Theorem 1.18 and comments after]{aliprantis2006}.
\\

We say that a linear map $\Psi: E \to F$ is a \emph{lattice homomorphism} if $\Psi x \wedge \Psi y = 0$ whenever $x \wedge y = 0$. A bijective lattice homomorphism is called a \emph{lattice isomorphism}. Observe that lattice isomorphisms preserve projection bands.\\

A linear map $T: E\rightarrow F$ between two Riesz spaces is called \emph{order continuous} if $x_{\alpha}\overset{o}{\rightarrow}x$ implies $Tx_{\alpha}\overset{o}{\rightarrow}x$ for any order converging net $\{x_\alpha\}$. A linear map $T: E\rightarrow F$ between two Riesz spaces is called \emph{interval preserving} if $T$ is positive and  
$T([0,x])=[0, T(x)]$.

\begin{proposition}\label{prop: IntervalPreserv}
Let $E, F$ be Riesz-spaces and $\Psi: E \rightarrow F$ a lattice isomorphism. Then $\Psi$ is both order continuous and interval preserving.
\end{proposition}
\begin{proof}
    By \cite[Theorem 2.21 (2)]{aliprantis2006} we get that $\Psi$ is order continuous. By \cite[Theorem 2.15]{aliprantis2006} we have that $\Psi$ and $\Psi^{-1}$ are positive. Hence, it suffices to show that for an order interval we have $\Psi([0,x])=[0,\Psi(x)]$. The inclusion $\Psi([0,x])\subseteq[0,\Psi(x)]$ is clear, since $\Psi$ is positive while the reverse inclusion follows similarly since $\Psi^{-1}$ is positive.
\end{proof}

A Riesz space $E$ equipped with a norm $\norm{\cdot}$ with the property that $|x|\leq|y|$ implies $\norm{x}\leq\norm{y}$ is called a \emph{normed Riesz space}. If $E$ is also norm complete, we call $E$ a \emph{Banach lattice}. Moreover, when $E$ is an $AL$-space, then $E$ has an order continuous norm and therefore $E$ is Dedekind complete, see \cite[Definition 4.7, comments after Definition 4.20 and Corollary 4.10]{aliprantis2006}. \\ 

If $E, F$ are Banach lattices, then every regular linear map is also norm bounded and so $\mathcal{L}_r(E, F) \subseteq \mathscr{B}(E,F)$, see \cite[Theorem 4.3]{aliprantis2006}. 

Furthermore, if $E$ and $F$ are $AL$-spaces, then every operator is regular, see \cite[Theorem 3.9]{abramovich2002}. It follows that, for $AL$-spaces, we have that the space of operators coincides with the space of regular linear maps, in other words
\begin{equation}\label{eq: reg-equal-bounded}
    \mathscr{B}(E,F) = \mathcal{L}_r(E, F),
\end{equation}
whenever $E, F$ are $AL$-spaces. In particular, since $F$ is Dedekind complete, every operator $T$ can be written as the difference of its positive and negative parts.

Moreover, we have 
\begin{equation}\label{eq: adjoint-and-modulus-commute}
    |T|^* = |T^*|,
\end{equation}
see \cite[Theorem 2.28]{abramovich2002}, and
\begin{equation}\label{eq: modulus-and-norm-relation}
    \norm{T} = \norm{|T|} = \norm{|T|^*} = \norm{|T^*|},
\end{equation}
see \cite[Corollary 3.10]{abramovich2002}. Finally, we recall the Riesz-Kantorovich formula, see \cite[Theorem 1.16]{abramovich2002}; for any $f$ in the positive cone $f \in E^+$ we have  
\begin{equation}\label{eq: Riesz-Kantorovich}
    |T|f = \sup \{Tg: -f \leq g \leq f \}.
\end{equation}

\subsection{Measure theoretic aspects.}\label{subsec: aspects-about-Riesz}

We recall the following definition, see \cite[Definition 211E]{fremlin2001} for example.

\begin{definition}\label{def: strictly-localizable}
    A measure space $(X, \Sigma, \mu)$ is called \emph{strictly localizable} or \emph{decomposable} if the following holds:
    \begin{enumerate}
        \item $X = \bigcup_{i \in I} X_i$, where $X_i \in \Sigma$ are mutually disjoint measurable sets such that $\mu(X_i) < \infty$ for every $i \in I$.
        \item A set $A \subseteq X$ is measurable if and only if $A \cap X_i$ is measurable for every $i \in I$.
        \item For every $A \in \Sigma$ we have $\mu(A) = \sum_{i \in I} \mu(A \cap X_i)$.
    \end{enumerate}
    We call the family $((X_i, \Sigma_i, \mu))_{i \in I}$ a \emph{decomposition} of $(X, \Sigma, \mu)$ where $\Sigma_i = \{A \cap X_i: A \in \Sigma \}$.
\end{definition}

Observe that every $\sigma$-finite measure space is strictly localizable.

If $(X, \Sigma, \mu)$ is strictly localizable, then the decomposition $((X_i, \Sigma_i, \mu))_{i \in I}$ also gives a decomposition of the space $L_1(X, \Sigma, \mu)$. Indeed, we have a decomposition
\begin{equation*}
    L_1(X, \Sigma, \mu) = \left(\bigoplus_{i \in I} L_1(X_i, \Sigma_i, \mu) \right)_{\ell_1},
\end{equation*}
so that any $f \in L_1(X, \Sigma, \mu)$ can be expressed as $f = \sum_{i \in I} P_i f$, where $P_i$ is the canonical projection onto $L_1(X_i, \Sigma_i, \mu)$ and the sum is absolutely convergent.

\begin{definition}
    Let $(X, \Sigma, \mu)$ be a measure space. The \emph{measure algebra induced by $(X,\Sigma, \mu)$} is the algebra of sets in $\Sigma$ modulo null sets. We denote this algebra by $\Sigma_\mu$.
\end{definition}

Throughout the text, we will write $A \in \Sigma_\mu$ to refer to the equivalence class of $A \in \Sigma$ whenever appropriate. Algebraic operations in $\Sigma_\mu$ carry over naturally from $\Sigma$, and we continue to use the same symbols for the operations in the induced measure algebra as in the original $\sigma$-algebra. This convention allows us to manipulate elements of the measure algebra while conveniently referring to representatives in $\Sigma$. It is also consistent with the structure of $L_1$-spaces, since $\mathds{1}_A$ denotes the indicator function of the set $A$ up to null sets.

\begin{definition}
    Let $\Sigma_\mu$ and $\tilde{\Sigma}_{\tilde{\mu}}$ be measure algebras. We say that $\Sigma_\mu$ and $\tilde{\Sigma}_{\tilde{\mu}}$ are \emph{isomorphic} if there exists a bijection $\phi: \Sigma_\mu \to \tilde{\Sigma}_{\tilde{\mu}}$ such that
    \begin{equation*}
        \phi(A  \backslash B) = \phi(A) \backslash \phi(B), \hspace{5pt} \text{  } \hspace{5pt}  \phi \left( \bigcup_{n=1}^\infty A_n \right) = \bigcup_{n=1}^\infty \phi(A_n), \hspace{5pt} \text{ and } \hspace{5pt} \mu(A) = \tilde{\mu}(\phi(A))
    \end{equation*}
    for every $A, B \in \Sigma_{\mu}$, $(A)_{n=1}^\infty \subseteq \Sigma_{\mu}$.
\end{definition}

Observe that if $\Sigma_\mu$ and $\tilde{\Sigma}_{\tilde{\mu}}$ are isomorphic measure algebras induced by $(X,\Sigma, \mu)$ and $(\tilde{X}, \tilde{\Sigma}, \tilde{\mu})$ respectively, then the isomorphism $\phi: \Sigma_\mu \to \tilde{\Sigma}_{\tilde{\mu}}$ induces an isometric lattice isomorphism $\Phi: L_1(X, \Sigma, \mu) \to L_1(\tilde{X}, \tilde{\Sigma}, \tilde{\mu})$ such that $\Phi(\mathds{1}_A) = \mathds{1}_{\phi(A)}$ for every $A \in \Sigma$. 

Given a family of subsets $(A_i)_{i \in I}$ with $A_i \subseteq X$, we denote by $\sigma(\{A_i : i \in I\})$ the $\sigma$-algebra generated by this family, that is, the smallest $\sigma$-algebra on $X$ containing all the sets $A_i$.

\begin{definition}
    A $\sigma$-algebra $\Sigma$ is said to be \emph{countably generated} if there exists a countable generating family $\{A_n: n \in \mathbb{N} \} \subseteq \Sigma$, i.e. $\sigma(\{A_n: n \in \mathbb{N} \}) = \Sigma$.
\end{definition}

We will need the following characterization of measure spaces of finite measure and countably generated $\sigma$-algebras; see \cite[page 173, Theorem C]{P.R.H} for the non-atomic case. The atomic case follows by decomposing $\mu$ into its atomic and non-atomic parts. We recall that a measure space $(\tilde{X}, \tilde{\Sigma}, \tilde{\mu})$ is called \emph{standard} if $\tilde{X}$ is a compact metric space, $\tilde{\Sigma}$ is the $\sigma$-algebra of Borel sets and $\tilde{\mu}$ is a finite measure. 

\begin{theorem}\label{th: separablefinite-to-borel}
    Let $(X, \Sigma, \mu)$ be a measure space of finite measure and $\Sigma$ be countably generated. Then there exists a standard measure space $(\tilde{X}, \tilde{\Sigma}, \tilde{\mu})$ so that $\Sigma_\mu$ is isomorphic to $\tilde{\Sigma}_{\tilde{\mu}}$.
\end{theorem}

\begin{rem}\label{rem: mu-restricted-notation}
    To avoid unnecessary notation, given a measure space $(X, \Sigma, \mu)$ and a $\sigma$-algebra $\tilde{\Sigma} \subseteq \Sigma$ we will still denote the restricted measure $\mu\big|_{\tilde{\Sigma}}$ by $\mu$, where the meaning should be clear from the context.
\end{rem}

For future reference, we will need the following observation.

\begin{proposition}\label{prop: separable-in-separable-sublattice}
    Let $(X, \Sigma, \mu)$ be a measure space and $M \subseteq L_1(X, \Sigma, \mu)$ be a separable subset. 
    \begin{enumerate}[label = (\alph*), ref = (\alph*)]
        \item \label{it: finite-case} If $(X, \Sigma, \mu)$ has finite measure, then there exists a countably generated $\sigma$-subalgebra $\tilde{\Sigma} \subseteq \Sigma$ such that $L_1(X, \tilde{\Sigma}, \mu)$ is a sublattice of $L_1(X, \Sigma, \mu)$ containing $M$. 
        \item \label{it: general-case-sl} If $(X, \Sigma, \mu)$ is strictly localizable, then there exist a $\sigma$-finite measurable set $\Omega \subseteq X$ and a countably generated $\sigma$-algebra $\hat{\Sigma}$ on $\Omega$ such that $L_1(\Omega, \hat{\Sigma}, \mu)$ is a sublattice of $L_1(X, \Sigma, \mu)$ containing $M$.
    \end{enumerate}
    Observe that, in both cases, these sublattices are $1$-complemented in $L_1(X, \Sigma, \mu)$.
\end{proposition}
\begin{proof}
    Let $(f_n)_{n=1}^\infty$ be a dense sequence in $M$. For each $c \in \mathbb{Q}$ and $n \in \mathbb{N}$ define $A_{n,c} = f^{^{-1}}_n((c, \infty))$. Let $\tilde{\Sigma} = \sigma(\{A_{n,c}: n \in \mathbb{N}, c \in \mathbb{Q} \})$, which is clearly countably generated. It is clear that $L_1(X, \tilde{\Sigma}, \mu)$ is a separable sublattice of $L_1(X, \Sigma, \mu)$ which contains $M$. 
    
    In the strictly localizable case, let $\Omega = \bigcup_{n \in \mathbb{N}} \bigcup_{c \in \mathbb{Q}} A_{n,c}$ and define $\hat{\Sigma} = \{\Omega \cap B: B \in \tilde{\Sigma}\}$ so that $\Omega$ is $\sigma$-finite and $\hat{\Sigma}$ is a countably generated $\sigma$-algebra.  Again, it is clear that $L_1(\Omega, \hat{\Sigma}, \mu)$ is a separable sublattice of $L_1(X, \Sigma, \mu)$ which contains $M$. 
    
    The complementation result follows, for example, from \cite[Chapter 6, Theorem 6]{lacey2008isometric}.
\end{proof}

We also recall the following well-known result.

\begin{proposition}\label{prop: sigma-finite-to-finite}
     Let $(Y, \Gamma, \nu)$ be a strictly localizable measure space, and let $\Theta \subseteq Y$ be a $\sigma$-finite measurable set with $\nu(\Theta) \not = 0$. Then there exists a measure $\tilde{\nu}$ on $\Gamma$ such that $(Y, \Gamma, \tilde{\nu})$ is strictly localizable, $\tilde{\nu}(\Theta) = 1$, and there is an isometric lattice isomorphism
    $\Psi: L_1(Y, \Gamma, \nu) \to L_1(Y, \Gamma, \tilde{\nu})$. 
    
    Furthermore, $\Psi$ commutes with multiplication operators, i.e. for any \break $h\in L_\infty(Y,\Gamma,\nu)$ and any $f\in L_1(Y,\Gamma,\nu)$, we have $\Psi(hf)=h\Psi(f)$, where on the right hand side $h$ is considered as an element of $L_\infty(Y,\Gamma,\tilde{\nu})$.
\end{proposition}
\begin{proof}
    Let $(\Theta_n)_{n \in \mathbb{N}}$ be a partition of $\Theta$ satisfying $0 < \nu(\Theta_n) < \infty$ and define $g = \sum_{n = 1}^\infty 2^{-n} \nu(\Theta_n)^{-1} \mathds{1}_{\Theta_n}$. By construction $\int_\Theta g(y) d\nu(y) = 1$. It is now straightforward to verify that if we define $\nu_1$ by $\nu_1(B) = \nu(B \cap (Y \backslash \Theta))$ for $B \in \Gamma$ and $d\nu_2 = g d\nu$ then the measure $\tilde{\nu} = \nu_1 + \nu_2$ satisfies $\tilde{\nu}(\Theta) = 1$ and the map $\Psi: L_1(Y, \Gamma, \nu) \to L_1(Y, \Gamma, \tilde{\nu})$ given by $f \mapsto \mathds{1}_{Y \backslash \Theta}f + \mathds{1}_{\Theta} f/g$ is an isometric lattice isomorphism. Obviously $(Y, \Gamma, \tilde{\nu})$ is strictly localizable, since $(Y, \Gamma, \nu)$ was. That $\Psi$ commutes with the multiplication operator can be easily seen from the construction above.
\end{proof}

\subsection{Integral operators and singular operators on $\sigma$-finite $L_1$-spaces.}\label{subsec: integral-operators}

We recall how an operator $T: L_1(X, \Sigma, \mu) \to L_1(Y, \Gamma, \nu)$, where the underlying measure spaces are $\sigma$-finite, can be decomposed into its integral and singular parts. We will need the following definition, see \cite[Definition 5.1]{abramovich2002}.

\begin{definition}\label{def: integral}
    Let $(X, \Sigma, \mu)$ and $(Y, \Gamma, \nu)$ be $\sigma$-finite measure spaces. An operator $K:L_{1}(X,\Sigma,\mu)\rightarrow L_{1}(Y,\Gamma,\nu)$ is called an \emph{integral operator} if there exists a measurable function $k:X \times Y \rightarrow \mathbb{R}$ such that $|k(\cdot,y)f(\cdot)|\in L_{1}(X,\Sigma,\mu)$ ($\nu$-a.e.) for every $f\in L_{1}(X,\Sigma,\mu)$ and 
    \begin{equation*}
        (Kf)(y)=\int_{X}k(x,y)f(x)d\mu(x) \hspace{15 pt} (\nu\text{-a.e.}).
    \end{equation*}
    The function $k$ is called the \emph{kernel} of the integral operator $K$.
\end{definition}

\begin{rem}\label{rem: modulus-integral-operators}
    Observe that if $K : L_{1}(X,\Sigma,\mu) \to L_{1}(Y,\Gamma,\nu)$ is an integral operator with kernel 
    $k : X \times Y \to \mathbb{R}$, then its modulus $|K| : L_{1}(X,\Sigma,\mu) \to L_{1}(Y,\Gamma,\nu)$ 
    is also an integral operator, with kernel $|k|$. This follows from \cite[Theorems~5.10 and~5.11]{abramovich2002}, 
    together with the fact that all operators between $AL$-spaces are regular \eqref{eq: reg-equal-bounded}.
\end{rem}

The following definition will also be useful throughout the paper.

\begin{definition}
     Let $(X, \Sigma, \mu)$ and $(Y, \Gamma, \nu)$ be $\sigma$-finite measure spaces and $K: L_{1}(X,\Sigma,\mu)\rightarrow L_{1}(Y,\Gamma,\nu)$ be an integral operator with kernel $k$. We define the \emph{truncated integral approximations} $K_s$, where $s \in \mathbb{R^+}$, to be the integral operators with kernel given by
     \begin{equation*}
         k_s(x, y) = 
            \begin{cases}
            \displaystyle k(x,y) & \text{if } | k(x,y)| \geq s, \\
            0 & \text{otherwise}.
            \end{cases}
     \end{equation*}
     We also refer to the kernels $k_s$ as \emph{truncated kernels}.
\end{definition}

We will need the following result concerning integral operators between finite $L_1$-spaces.

\begin{proposition}\label{prop: KernelinL1}
    Let $(X, \Sigma, \mu)$ and $(Y, \Gamma, \nu)$ be measure spaces of finite measure and $K: L_1(X, \Sigma, \mu) \to L_1(Y, \Gamma, \nu)$ be an integral operator with kernel $k: X \times Y \to \mathbb{R}$. Then $k \in L_1(X \times Y, \Sigma \otimes \Gamma, \mu \otimes \nu)$.
\end{proposition}
\begin{proof}
   Since $(X, \Sigma, \mu)$ is finite, we have that $\mathds{1}_X \in L_1(X, \Sigma, \mu)$. It follows by using Remark \ref{rem: modulus-integral-operators} and Tonelli's Theorem that
    \begin{equation*}
       \int_{X \times Y}|k(x,y)| d(\mu(x) \otimes d\nu(y)) = \int_{Y} \int_X |k(x,y)| d\mu(x) d\nu(y) = \norm{|K| \mathds{1}_X }  < \infty,
    \end{equation*}
    that is $k \in L_1(X \times Y, \Sigma \otimes \Gamma, \mu \otimes \nu)$.
\end{proof}

In the following theorem, let $E = L_1(X, \Sigma, \mu)$ and $F = L_1(Y, \Gamma, \nu)$, where $(X, \Sigma, \mu)$ and $(Y, \Gamma, \nu)$ are $\sigma$-finite measure spaces. We denote by $\mathcal{L}_\kappa(E, F)$ the space of all regular integral operators. We have the following, where the first equality follows from \eqref{eq: reg-equal-bounded} and the second from \cite[Theorem~5.14]{abramovich2002}, together with the fact that $\mathcal{L}_{r}(E, F)$ is Dedekind complete and thus all bands are projection bands.

\begin{theorem}[Bukhvalov-Luxemburg-Schep-Zaanen] \label{th: Decomposition-integral-and-singular}
The vector space $\mathcal{L}_{\kappa}(E,F)$ is a projection band in $\mathcal{L}_{r}(E,F)$. In particular
\begin{equation*}
    \mathscr{B}(E,F)=\mathcal{L}_{r}(E,F)=\mathcal{L}_{\kappa}(E,F)\oplus\mathcal{L}_{\kappa}(E,F)^{d}.
\end{equation*}
\end{theorem}

An element $S \in \mathcal{L}_{\kappa}(E,F)^{d}$ is called \emph{singular} (not to be confused with strictly singular operators, which are precisely the weakly compact operators in the context of $L_1$-spaces). Observe that $S$ is a singular operator if and only if
\begin{equation}\label{eq: singular-if-only-if-orthogonal}
    |S| \wedge |K| = 0  \text{ for all integral operators } K,
\end{equation}
or equivalently
\begin{equation*}
    |S-K|=|S+K| \text{ for all integral operators } K.
\end{equation*}
We will see later that singular operators $S$ are far from being weakly compact.

Also, note that from the previous theorem, we have that any operator $T \in \mathscr{B}(E,F)$ can be decomposed into its singular and integral part, that is $T = S + K$.

\subsection{Weak compactness in $L_1$-spaces.}\label{subsec: weak-compactness-finite}

We recall the following classical measure of weak noncompactness for measure spaces of finite measure.

\begin{definition}
    Let $(X, \Sigma, \mu)$ be a measure space of finite measure, and let $A \subseteq L_1(X, \Sigma, \mu)$ be bounded. The \emph{modulus of equi-integrability} of $A$ is defined by
    \begin{equation*}
        \Delta(A) = \limsup_{\mu(B) \to 0} \sup_{f \in A} \norm{\chi_B f}.
    \end{equation*}
\end{definition}

The modulus of equi-integrability is a measure of weak noncompactness of bounded subsets of $L_1(X, \Sigma, \mu)$. In particular, a bounded set $A \subseteq L_1(X, \Sigma, \mu)$ is relatively weakly compact if and only if $\Delta(A) = 0$, see \cite[Theorem 5.2.9]{albiac2006topics}. If $T: E \to L_1(X, \Sigma, \mu)$ is an operator, we define $\Delta(T) = \Delta(TB_E)$ so that $T$ is weakly compact if and only if $\Delta(T) = 0$.

\begin{rem}
    We would like to emphasize that the modulus of equi-integrability, as previously defined, only characterises weak compactness for measure spaces of finite measure. Indeed, consider for example $X = \mathbb{R}$, $\Sigma$ its Borel $\sigma$-algebra and $\lambda$ the Lebesgue measure. Then the set $A = \{\mathds{1}_{[n, n+1]}: n \in \mathbb{N}\} \subseteq L_1(\mathbb{R}, \Sigma, \lambda)$ is clearly not weakly compact but $\Delta(A) = 0$. Therefore, a different quantity will be needed in the general case. 
\end{rem}

\bigskip
\section{The proof of Theorem \ref{th: MainTheorem}}\label{sec: mainTheorem}

We begin by briefly outlining the steps in the hope that this may be helpful to the reader and provide a sense of direction. The argument proceeds in four stages, according to the structure of the underlying measure spaces: we first treat standard measure spaces, then move to measure spaces of finite measure and countably generated $\sigma$-algebra, then to arbitrary measure spaces of finite measure, and finally to strictly localizable measure spaces. \\

In Section~\ref{subsect: Standardmeasurespace}, we treat the case of standard measure spaces. This reduces to introducing a result of Weis, Theorem~\ref{th: Weis-standard-measure-case}, which essentially establishes Theorem \ref{th: MainTheorem} in this setting. \\

In Section~\ref{subsect: countablygenerated}, we address measure spaces of finite measure and countably generated $\sigma$-algebra. The idea is to use Theorem ~\ref{th: separablefinite-to-borel} to reduce to the case of standard measure spaces, allowing us to apply Weis's theorem. To ensure that this reduction works, we must show that the isomorphism preserves singular operators and the structure of kernels of integral operators. This is shown in Lemma~\ref{lmm: TransKTruncatedK}. Using this, we establish Theorem~\ref{th: countably-generated-finite-case}, extending Weis's result to finite measure spaces with countably generated $\sigma$-algebra. \\

In Section~\ref{subsection: finitemeasurespace}, we treat arbitrary spaces of finite measure by reducing to the previous case via two technical lemmas, Lemmas~\ref{lmm: technical1} and~\ref{lmm: technical2}, which rely on the notion of rich families; see Definition~\ref{def: richFamily}. To be precise, Lemma~\ref{lmm: technical1} shows that, given a separable subset $M \subseteq L_1(X, \Sigma, \mu)$ and a sequence of operators $(T_s)_{s \in \mathbb{N}}$, we can find a $L_1$-sublattice containing $M$ and induced by a measure space with a countably generated $\sigma$-algebra, while preserving the moduli $|T_s|$ for all $s \in \mathbb{N}$. Building on this, Lemma~\ref{lmm: technical2} ensures that for a singular operator $S$ and a separable set $M$, there exists an $L_1$-sublattice, induced by a measure space with a countably generated $\sigma$-algebra, which contains $M$ and in which the restriction of $S$ remains singular. This allows us to handle the singular part of the operator in Proposition~\ref{prop: singular-operators-finite-case}, while Lemma~\ref{lmm: kernel-operator-finite-case} addresses the integral part. Since these are the essential ingredients we will need from the finite measure setting, we postpone the proof of the main theorem for finite measures to the next section, where it will follow as a corollary of the localized version; see Corollary~\ref{cor: general-finitec-case}. \\

Finally, in Section~\ref{subsection: strictly-localizable}, we address the strictly localizable case. The idea is to localize the behaviour of $T$ to an arbitrary measurable set $\Theta \subseteq Y$ of finite measure, reduce to the finite measure case, and then show that this localized version suffices. To this end, we introduce a localized version of $\Delta(T)$, given by $\inf_{\Theta \subseteq Y, \nu(\Theta) < \infty} \alpha_{\Theta}(T)$, and focus on the study of $\alpha_{\Theta}(T)$. We prove Proposition~\ref{prop: alpha-representation}, which provides a localized representation that serves as the counterpart of the expression for $\Delta(T)$ in the finite measure setting. Using this, we establish a localized form of the main theorem, Theorem~\ref{th: main-with-Theta}. The rest of the section is now devoted to using this localized version to get the general one. For technical reasons, we need to explore the behaviour in $\sigma$-finite localized versions (instead of finite ones); note, however, that Proposition \ref{prop: sigma-finite-to-finite} allows us to go back to finite localizations, at least in terms of operators. In particular, we prove Proposition~\ref{prop: first-part-MainTheorem} which shows the existence of a best approximant, and Proposition~\ref{prop: second-part-MainTheorem} which establishes that
\begin{equation*}
    \norm{T}_w = \inf_{\substack{\Theta \subseteq Y,  \nu(\Theta) < \infty}} \alpha_{\Theta}(T).
\end{equation*}
These two propositions yield Theorem~\ref{th: MainTheorem}.

\subsection{Standard measure spaces.}\label{subsect: Standardmeasurespace} Let us focus on the particular case in which the underlying spaces $(X,\Sigma,\mu)$ and $(Y,\Gamma,\nu)$ are standard measure spaces. In this setting, for an operator $T: L_{1}(X,\Sigma,\mu)\rightarrow L_{1}(Y,\Gamma,\nu)$, by \cite[Theorem 2.1]{weis1984representation},  there exist measures  $(\mu_{y})_{y\in Y}$ on $(X, \Sigma)$ such that
\begin{equation*}
    (Tf)(y)=\int_{X}f(x)d\mu_{y}(x).
\end{equation*}
We can use Lebesgue decomposition to obtain $\mu_y = \mu_y^i + \mu_y^s$ where $\mu_y^i$ is $\mu$-absolutely continuous and $\mu_y^s$ is $\mu$-singular. By \cite[Proposition 2.6]{weis1984representation}, this gives a decomposition $T = S + K$ defined by
\begin{equation*}
     (Sf)(y) = \int_{X}f(x)d\mu^s_{y}(x) \hspace{5 pt} \text{ and } \hspace{5pt} (Kf)(y) = \int_{X}f(x)d\mu^i_{y}(x) = \int_{X}f(x) k(x,y) d\mu(x),
\end{equation*}
where $k( \cdot ,y) = \frac{d\mu_y^i}{d\mu_y}$ is the Radon-Nikodym derivative. 

There is a similar representation for $T^*: L_{\infty}(Y,\Gamma,\nu) \to L_{\infty}(X,\Sigma,\mu)$ by measures $(\nu_{y})_{y\in Y}$ on $(Y, \Gamma)$ so that
\begin{equation*}
    (T^*g)(y)=\int_{X}g(x)d\nu_{y}(x).
\end{equation*}
Again, we can decompose $\nu_y = \nu_y^i + \nu_y^s$ where $\nu_y^i$ is $\nu$-absolutely continuous and $\nu_y^s$ is $\nu$-singular. Using this decomposition, the adjoint of $S$ and $K$ are given by
\begin{align*}
    (S^* g)(x)=
    \int_{Y}g(y)d\nu^s_{x}(y) \hspace{5pt} \text{ and } \hspace{5pt} (K^* g)(x)=
    \int_{Y}g(y)d\nu^i_{x}(y) = \int_{Y}g(y)h(x,y)d\nu(y),
\end{align*}
where $h(x, \cdot) = \frac{d\nu_x^i}{d\nu_x}$ is the Radon-Nikodym derivative. The previous representations are unique and satisfy $|\mu_{y}^{s}|\wedge|\mu_{y}^{i}|=0$ $\nu$-a.e. and $|\nu_{x}^{s}|\wedge|\nu_{x}^{i}|=0$ $\mu$-a.e., see \cite[Remark 2.3]{weis1984representation}. By Tonelli's theorem, one easily checks that actually $h(x,y) = k(x,y)$.

Observe that this decomposition is a particular case of the decomposition into integral and singular part, given by Theorem \ref{th: Decomposition-integral-and-singular}. In particular, in a standard measure space, an operator is singular if and only if it has a singular representation in the sense of measures discussed before.

In the setting of standard measure spaces, Weis \cite{weis1984approximation} has shown how these representations can be used to give a formula for $\Delta(T)$.

\begin{theorem}[Weis]\label{th: Weis-standard-measure-case}
    Let $(X, \Sigma, \mu)$ and $(Y, \Sigma, \nu)$ be standard measure spaces, and $T: L_1(X, \Sigma, \mu) \to L_1(Y, \Sigma, \nu)$ be an operator. Let $(\nu_y^i)_{x \in X}$ and $(\nu_y^s)_{x \in X}$ be the absolutely continuous and singular parts of the measure representation of $T^*$. Then there exists a weakly compact operator $H$ such that
    \begin{equation*}
        \norm{T}_w = \norm{T - H} = \Delta(T) =  \inf_{n \in \mathbb{N}}\esssup_{x\in X}\{\rVert \nu_{x}^s \lVert +\int_{Y}|k_{n}(x,y)|d\nu(y)\},
    \end{equation*}
    where $k_n$ are the truncated kernels
    \begin{equation*}
        k_n(x, y) = 
            \begin{cases}
            \displaystyle \frac{d\nu_x^i}{d\nu}(y) & \text{if } \left| \frac{d\nu_x^i}{d\nu}(y) \right| \geq n, \\
            0 & \text{otherwise}.
            \end{cases}
    \end{equation*}
\end{theorem}

\begin{rem}
    We note that $k(x, y)$, as previously defined, is in fact the kernel of the integral operator $K$ from the decomposition $T = S + K$ given by Theorem~\ref{th: Decomposition-integral-and-singular}. Using this decomposition, the previous expression for $\Delta(T)$ can also be expressed as
    \begin{equation*}
        \Delta(T) =  \inf_{n \in \mathbb{N}}\esssup_{x\in X}\{|S|^*\mathds{1}_Y(x) +\int_{Y}|k_{n}(x,y)|d\nu(y)\}.
    \end{equation*}
\end{rem}

\subsection{Measure spaces of finite measure with countably generated $\sigma$-algebra.}\label{subsect: countablygenerated}

Our goal is to extend the result from standard measure spaces to measure spaces of finite measure and countably generated $\sigma$-algebra. This will be achieved by combining Weis's theorem with Theorem~\ref{th: separablefinite-to-borel}.

We now argue that singular and integral operators are preserved under lattice isomorphisms, and how the truncated approximations of integral operators transform under such isomorphisms.

\begin{lemma}\label{lmm: TransKTruncatedK}
    Let $(X, \Sigma, \mu)$ and $(Y, \Gamma, \nu)$ be measure spaces of finite measure with $\Sigma$ and $\Gamma$ countably generated. Let $S: L_1(X, \Sigma, \mu) \to L_1(Y, \Gamma, \nu)$ be a singular operator, and $K: L_1(X, \Sigma, \mu) \to L_1(Y, \Gamma, \nu)$ be an integral operator.

    \begin{enumerate}[label=(\roman*) , ref=(\roman*)]
        \item \label{it: a1} Suppose $(\tilde{X}, \tilde{\Sigma}, \tilde{\mu})$ is a measure space such that $\phi: \tilde{\Sigma}_{\tilde{\mu}} \to \Sigma_\mu $ is an isomorphism and let $\Phi: L_1(\tilde{X}, \tilde{\Sigma}, \tilde{\mu}) \to L_1(X, \Sigma, \mu)$ be the induced lattice isomorphism. Then the operator $\tilde{S} = S \Phi$ is a singular operator and $\tilde{K} = K \Phi$ is an integral operator with $\tilde{K}_n = K_n \Phi$ as its truncated integral approximations.
        \item \label{it: a2} Suppose that $(\tilde{Y}, \tilde{\Gamma}, \tilde{\nu})$ is a measure space such that $\psi: \Gamma_\nu \to \tilde{\Gamma}_{\tilde{\nu}}$ is an isomorphism and let  $\Psi: L_1(Y,\Gamma,\nu) \to L_1(\tilde{Y}, \tilde{\Gamma}, \tilde{\nu})$ be the induced lattice isomorphism. Then the operator $\hat{S} = \Psi S$ is a singular operator and $\hat{K} = \Psi K$ is an integral operator with $\hat{K}_n = \Psi K_n$ as its truncated integral approximations.
    \end{enumerate}
\end{lemma}
\begin{proof}
    \ref{it: a1} We start by proving $\tilde{K} = K \Phi$ is an integral operator and show that the truncated approximations are transformed in the natural way. Without loss of generality, decomposing $K = K^+ - K^-$ and considering $K^+$ and $K^-$ separately, we may assume that $K$ is positive. Let $k: X \times Y \to \mathbb{R}$ be the kernel of $K$. By Proposition \ref{prop: KernelinL1} we have $k \in L_1(X \times Y, \Sigma \otimes \Gamma, \mu \otimes \nu)$. Therefore, we can choose $(k^{(m)})_{m \in \mathbb{N}}$ a sequence of simple functions $k^{(m)}:  X \times Y \to \mathbb{R}$ converging to $k$ in norm,
    \begin{equation*}
        k^{(m)} =  \sum_{l=1}^{N_m} a_{m,l}\cdot\mathds{1}_{A_{m,l}\times B_{m,l}}.
    \end{equation*}

    Define the simple functions $\tilde{k}^{(m)}: \tilde{X} \times Y \to \mathbb{R}$ by
    \begin{equation*}
        \tilde{k}^{(m)} =  \sum_{l=1}^{N_m} a_{m,l}\cdot\mathds{1}_{\phi^{-1}(A_{m,l})\times B_{m,l}}.
    \end{equation*}
    Since the sequence $(k^{(m)})_{m \in \mathbb{N}}$ is Cauchy, so is $(\tilde{k}^{(m)})_{m \in \mathbb{N}} \subseteq L_1(\tilde{X} \times Y, \tilde{\Sigma} \otimes \Gamma, \tilde{\mu} \otimes \nu)$, thus it is convergent. Denote by $\tilde{k}$ the limit of this sequence and define the integral operator $\tilde{H} \tilde{f}(y) = \int_{\tilde{X}} \tilde{k}(\tilde{x}, y) \tilde{f}(\tilde{x}) d\tilde{\mu}(\tilde{x})$, we will show that $\tilde{H} = \tilde{K}$, so that $\tilde{K}$ is in fact an integral operator.
    
    Passing to a subsequence, we can assume that $k^{(m)}$ and $\tilde{k}^{(m)}$ converge pointwise a.e. Define $k_y = k(\cdot, y)$, $\tilde{k}_y = \tilde{k}(\cdot, y)$, $k^{(m)}_y = k^{(m)}(\cdot, y)$, $\tilde{k}^{(m)}_y = \tilde{k}^{(m)}(\cdot, y)$. It follows that for $\nu$-almost all $y$ we have that $k^{(m)}_y \to k_{y}$ pointwise $\mu$-a.e and $\tilde{k}^{(m)}_y \to \tilde{k}_{y}$ pointwise $\tilde{\mu}$-a.e.

   It is easy to see that for any $c \in \mathbb{R}$ we have
   \begin{equation*}
       \{\tilde{k}_{y}^{(m)}>c\}=\phi^{-1}(\{k_{y}^{(m)} > c\}),
   \end{equation*}
   where the equality is taken in the equivalence classes modulo null sets. Therefore, for $c \in \mathbb{R}$ we get
   \begin{align*}
        \{\tilde{k}_{y}>c\}&=\bigcup_{n_{0}=1}^{\infty}\bigcap_{m=n_0}^{\infty}\{\tilde{k}^{(m)}_{y}>c\} =\bigcup_{n_{0}=1}^{\infty}\bigcap_{m=n_0}^{\infty}\phi^{-1}(\{k^{(m)}_{y}>c\})\\
        &=\phi^{-1}(\bigcup_{n_{0}=1}^{\infty}\bigcap_{m=n_0}^{\infty}\{k^{(m)}_{y}>c\}) =\phi^{-1}(\{k_{y}>c\}).
    \end{align*}
    Thus, for $f = \Phi \tilde{f}$ we have
    \begin{align*}
        (\tilde{H}_{n}\tilde{f})(y)&=\int_{\{\tilde{k}_y>n\}}\tilde{k}(\tilde{x},y)\tilde{f}(\tilde{x})d\tilde{\mu}(\tilde{x}) =\int_{\phi^{-1}(\{k_y>n\})}\tilde{k}(\tilde{x},y)\tilde{f}(\tilde{x})d\tilde{\mu}(\tilde{x})\\
        &=\int_{\{k_y>n\}}k(x,y)f(x)d\mu(x)=(K_n f)(y),
    \end{align*}
    that is $\tilde{H}_n \tilde{f} = K_n \Phi \tilde{f}$. In particular, for $n = 0$ we get that $\tilde{H} = K\Phi = \tilde{K}$, as desired. Moreover, this proves that the truncated kernels transform in the natural way.
    The preceding shows, in particular, that right multiplication by $\Phi$ preserves integral operators. \\
    
    To show that $\Phi$ also preserves singular operators, it suffices to verify that right multiplication by $\Phi$ defines a lattice isomorphism on the space of operators. Once this is established, it follows that right multiplication by $\Phi$ preserves projection bands. Moreover, as we showed before, multiplication by $\Phi$ is a bijection from \break $\mathcal{L}_\kappa(L_1(X,\Sigma,\mu),L_1(Y,\Gamma,\nu))$ onto $\mathcal{L}_\kappa(L_1(\tilde{X},\tilde{\Sigma},\tilde{\mu}),L_1(Y,\Gamma,\nu))$, and the result now follows immediately from Theorem \ref{th: Decomposition-integral-and-singular}.

    Thus, we are only left to show that right multiplication by $\Phi$,
        \begin{equation*}
            \mathscr{B}(L_1(X, \Sigma, \mu), L_1(Y, \Gamma, \nu)) \to \mathscr{B}(L_1(\tilde{X}, \tilde{\Sigma}, \tilde{\mu}), L_1(Y, \Gamma, \nu)), \hspace{5pt} T \mapsto T\Phi
        \end{equation*}
    is a lattice isomorphism. Observe that since $\Phi$ is a lattice isomorphism, it is interval preserving by Proposition \ref{prop: IntervalPreserv}. By \cite[Theorem 2.16 (2)]{aliprantis2006}, it follows that right multiplication by $\Phi$ is a lattice homomorphism, while it is clearly bijective and thus a lattice isomorphism. \\

   \ref{it: a2} 
    We proceed by first defining the simple functions
    $\hat{k}^{(m)} : X \times \tilde{Y} \to \mathbb{R}$ by
    \begin{equation*}
        \hat{k}^{(m)} = \sum_{l=1}^{N_m} a_{m,l} \,\mathds{1}_{A_{m,l} \times \psi(B_{m,l})}.
    \end{equation*}
    As above, the sequence $(\hat{k}^{(m)})_{m \in \mathbb{N}} \subseteq L_1(X \times \tilde{Y}, \Sigma \otimes \tilde{\Gamma}, \mu \otimes \tilde{\nu})$ is Cauchy, thus it is convergent; denote by $\hat{k}$ the limit of this sequence. Let $\hat{H}$ be the integral operator with kernel $\hat{k}$, and we will show that $\hat{H} = \hat{K}$, so that $\hat{K}$ is an integral operator.
    
    The argument then proceeds analogously to that for \ref{it: a1}, but now applied to the adjoint operator. We can define
    \begin{equation*}
        k_x = k(x, \cdot), \quad \hat{k}_x = \hat{k}(x, \cdot), \quad k_x^{(m)} = k^{(m)}(x, \cdot), \quad \hat{k}_x^{(m)} = \hat{k}^{(m)}(x, \cdot).
    \end{equation*}

     Arguing as above, it is straightforward to see that $\{\hat{k}_{x} > c\} = \psi(\{k_{x} > c\})$. Consider the adjoint $(\hat{H}_n)^*: L_\infty (\tilde{Y}, \tilde{\Gamma}, \tilde{\nu}) \to L_\infty(X, \Sigma, \mu)$. Similarly as before, $\Psi$ is a positive and order continuous operator, and thus it is interval preserving by Proposition \ref{prop: IntervalPreserv}. By \cite[Theorem 2.19]{aliprantis2006}, $\Psi^*$ is a lattice homomorphism, while it is clearly an isometry and thus it is a lattice isometric isomorphism. For any $\tilde{g} \in L_\infty(\tilde{Y}, \tilde{\Gamma}, \tilde{\nu})$, define $g = \Psi^* \tilde{g}$, we have
        \begin{align*}
        ((\hat{H}_{n})^{*}\tilde{g})(x)&=\int_{\{\hat{k}_x>n\}}\hat{k}(x,\tilde{y})\tilde{g}(\tilde{y})d\tilde{\nu}(\tilde{y}) =\int_{\psi(\{k_x>n\})}\hat{k}(x,\tilde{y})\tilde{g}(\tilde{y})d\tilde{\nu}(\tilde{y})\\
        &=\int_{\{k_y>n\}}k(x,y)g(y)d\nu(y)=((K_{n})^{*}g)(y) = (((K_n)^*\Psi^{*})\tilde{g})(x),
    \end{align*}
    and therefore $(\hat{H}_{n})^{*}=(K_n)^*\Psi^{*}=(\Psi K_n)^{*}$, so that $\hat{H}_{n}=\Psi K_n$. Taking $n = 0$ gives that $\hat{H} = \Psi K = \hat{K}$, so that $\hat{K}$ is an integral operator and the truncated approximations transform as claimed. The preceding shows that left multiplication by $\Psi$ preserves integral operators.
    
    It remains to show that $\Psi$ also preserves singular operators. As before, since left multiplication by $\Psi$ is a bijection from $\mathcal{L}_\kappa(L_1(X, \Sigma, \mu), L_1(Y, \Gamma, \nu))$ onto \break $\mathcal{L}_\kappa(L_1(X, \Sigma, \mu), L_1(\tilde{Y}, \tilde{\Gamma}, \tilde{\nu}))$, it suffices to verify that left multiplication by $\Psi$ defines a lattice isomorphism on the space of operators.

    Thus, we need to show that left multiplication by $\Psi$,
    \begin{equation*}
        \mathscr{B}(L_1(X, \Sigma, \mu), L_1(Y, \Gamma, \nu)) \to \mathscr{B}(L_1(X, \Sigma, \mu), L_1(\tilde{Y}, \tilde{\Gamma}, \tilde{\nu})), \hspace{5pt} T \mapsto \Psi T
    \end{equation*}
    is a lattice isomorphism. $\Psi$ is order continuous and thus by \cite[Theorem 2.17]{aliprantis2006} left multiplication by $\Psi$ is a lattice homomorphism, while it is clearly bijective and thus a lattice isomorphism. This completes the proof.
\end{proof}

We are ready for our main theorem in the setting of measure spaces of finite measure with countably generated $\sigma$-algebras.

\begin{theorem}(Countably generated $\sigma$-algebra case) \label{th: countably-generated-finite-case}
    Let $(X, \Sigma, \mu)$ and $(Y, \Gamma, \nu)$ be measure spaces of finite measure with $\Sigma$ and $\Gamma$ countably generated. Let $T: L_1(X, \Sigma, \mu) \to L_1(Y, \Gamma, \nu)$ be an operator and decompose $T = S + K$, where $S$ is a singular operator and $K$ is an integral operator. Then
    \begin{equation*}
        \Delta(T)=\inf_{n \in \mathbb{N}}\esssup_{x\in X}\{|S|^{*}\mathds{1}_{Y}(x)+\int_{Y}|k_{n}(x,y)|d\nu(y)\},
    \end{equation*}
    where $k_n(x,y)$ is the truncated kernel of $K$. In particular, for a singular operator $S$, $\Delta(S) = \norm{S^*} = \norm{S}$.
\end{theorem}
\begin{proof}
    Applying Theorem \ref{th: separablefinite-to-borel}, we can find standard measure spaces $(\tilde{X}, \tilde{\Sigma}, \tilde{\mu})$ and $(\tilde{Y}, \tilde{\Gamma}, \tilde{\nu})$, as well as isometric lattice isomorphisms
    \begin{equation*}
        \Phi: L_1(\tilde{X}, \tilde{\Sigma}, \tilde{\mu}) \to L_1(X, \Sigma, \mu)
        \quad \text{and} \quad
        \Psi: L_1(Y, \Gamma, \nu) \to L_1(\tilde{Y}, \tilde{\Gamma}, \tilde{\nu}),
    \end{equation*}
    mapping characteristic functions to characteristic functions. Decompose $T = S + K$ into its singular and integral parts and consider the operator
    \begin{equation*}
        \tilde{T} = \Psi T \Phi = \Psi S \Phi + \Psi K \Phi = \tilde{S} + \tilde{K}
    \end{equation*}
    where $\tilde{S} = \Psi S \Phi$ is a singular operator and $\tilde{K} = \Psi K \Phi$ is an integral operator by Lemma \ref{lmm: TransKTruncatedK}

     Let $(\tilde{\nu}_{\tilde{x}})_{\tilde{x}\in \tilde{X}}$ be the representation of the operator $\tilde{T}^{*}$, with $(\tilde{\nu}_{\tilde{x}}^{s})_{\tilde{x}\in \tilde{X}}$ its singular part and $\tilde{k}(\cdot,\cdot)$ the kernel induced by the absolutely continuous part. The operator $\tilde{S}^{*}$ is also singular (see \cite[Theorem 5.3]{weis1984representation}) and its measure representation is equal to $(\tilde{\nu}_{\tilde{x}}^{s})_{\tilde{x}\in \tilde{X}}$ since the decomposition into singular and absolutely continuous part of the measure representation is unique and the space of integral operators is a projection band (see Theorem \ref{th: Decomposition-integral-and-singular} and \cite[Remark 2.3]{weis1984representation}).

     By Lemma \ref{lmm: TransKTruncatedK}, the truncated integral approximations $\tilde{K}_n$ of $\tilde{K}$ can be expressed as $\tilde{K}_n = \Psi K_n \Phi$ where $K_n$ are the truncated integral approximations of $K$.

     Define $V_n = S + K_n$ and $\tilde{V}_n = \tilde{S} + \tilde{K}_n = \Psi V_n \Phi$. Since $\norm{|V_n|^*} = \norm{|\tilde{V}_n|^*}$ we get
     \begin{align*}
        \esssup_{x\in X}\{|S|^{*}\mathds{1}_{Y}(x)+ \int_{Y}|k_{n}(x,y)|d\nu(y)\}&=\norm{|S|^{*}+|K_n|^{*}}\\
        &=\norm{|V_{n}|^{*}} = \norm{|\tilde{V_n}|^{*}} = \norm{|\tilde{S}|^{*}+|\tilde{K}_{n}|^{*}}\\
        &=\esssup_{\tilde{x}\in \tilde{X}}\{\norm{\tilde{\nu}_{\tilde{x}}^{s}} + \int_{Y}|\tilde{k}_{n}(\tilde{x},\tilde{y})|d\tilde{\nu}(\tilde{y})\}.
    \end{align*}
    From the definition of $\Delta$, it is easy to see that $\Delta(T) = \Delta(\tilde{T})$. Applying Theorem \ref{th: Weis-standard-measure-case} to $\tilde{T}$ gives
    \begin{align*}
        \Delta(T) &= \Delta(\tilde{T}) \stackrel{\text{Thm. \ref{th: Weis-standard-measure-case}}}{=}\inf_{n \in \mathbb{N}}\esssup_{\tilde{x}\in \tilde{X}}  \{\norm{\tilde{\nu}_{\tilde{x}}^{s}} + \int_{Y}|\tilde{k}_{n}(\tilde{x},\tilde{y})|d\tilde{\nu}(\tilde{y})\} \\
        &=\inf_{n \in \mathbb{N}}\esssup_{x\in X}\{|S|^{*}\mathds{1}_{Y}(x)+ \int_{Y}|k_{n}(x,y)|d\nu(y)\}.
    \end{align*}    
    This completes the first part of the proof.  If $S$ is a singular operator, the previous result immediately gives $\Delta(S) = \norm{S^*} = \norm{S}$.
\end{proof}

\subsection{Measure spaces of finite measure.}\label{subsection: finitemeasurespace} 

We now move to the setting of arbitrary finite measure spaces. We start by recalling the notion of rich families, introduced by Borwein and Moors \cite{borwein2000Separable}, which will allow us to formulate the forthcoming lemmas in a unified way.

\begin{definition}\label{def: richFamily}
Let $X$ be a Banach space. We say that a family $\mathcal{R}$ of closed separable subspaces of $X$ is \emph{rich} if it satisfies the following conditions:
\begin{enumerate}[label = (\alph*), ref = (\alph*)]
        \item \label{it: def-rich-a} For every separable subset $M\subseteq X$ there is an $E\in\mathcal{R}$, such that $M\subseteq E$.
        \item \label{it: def-rich-b} For every increasing sequence $(E_n)_{n\in\mathbb{N}} \subseteq \mathcal{R}$, the subspace $E:=\overline{\bigcup_{n\in\mathbb{N}} E_n}\in \mathcal{R}$.
    \end{enumerate}
\end{definition}

Recall that, as mentioned in Remark \ref{rem: mu-restricted-notation}, given a measure space $(X,\Sigma,\mu)$ and $\sigma$-algebra $\tilde{\Sigma} \subseteq \Sigma$, we always denote the restricted measure $\mu|_{\tilde{\Sigma}}$ simply by $\mu$.

\begin{example} Let $(X,\Sigma,\mu)$ be a finite measure space. Then  Proposition \ref{prop: separable-in-separable-sublattice} shows that
\begin {align*} \{L_1(X,\tilde{\Sigma},\mu): \tilde{\Sigma}\subseteq \Sigma \text { countably generated } \sigma\text{-algebra} \}
\end{align*}
is a rich family of closed separable subspaces of $L_{1}(X,\Sigma,\mu)$.
\end{example}

One advantage of rich families is that the countable intersection of rich families is again a rich family, see \cite[Proposition 1.1]{borwein2000Separable}.
\begin{proposition}\label{prop: IntersectionRFamily} Let $X$ be a Banach space and let $\mathcal{R}_{n}$ be a rich family of closed separable subspaces of $X$ for each $n\in\mathbb{N}$. Then
\begin{equation*}
\mathcal{R}:=\bigcap_{n\in\mathbb{N}}\mathcal{R}_n
\end{equation*}
is a rich family of closed separable subspaces of $X$.
\end{proposition}

We focus on the singular part of the operator. Before doing so, however, we require the following technical result. From now on, if $E$ is a closed subspace of $L_1(X, \Sigma, \mu)$, we will denote by $J_E: E \to L_1(X, \Sigma, \mu)$ the inclusion operator.

\begin{lemma}\label{lmm: technical1}
    Let $(X, \Sigma, \mu)$ and $(Y, \Gamma, \nu)$ be measure spaces of finite measure. Then the following assertions hold:
     \begin{enumerate}[label=(\roman*) , ref=(\roman*)]
        \item \label{it: rich-families-a1}
        For any operator $T: L_1(X, \Sigma, \mu) \to L_1(Y, \Gamma, \nu)$ the family
         \begin{align*}
            \mathcal{R}_{T}:= \{E : \hspace{3pt} & E = L_{1}(X,\tilde{\Sigma},\mu) \text{ for }\tilde{\Sigma}\subseteq \Sigma \text{ countably generated } \sigma\text{-algebra}, \\
            & \text{and }  |TJ_{E}|=|T|J_{E}\} 
        \end{align*}
         is a rich family of closed separable subspaces of $L_1(X, \Sigma, \mu)$.
        \item \label{it: rich-families-a2} For a sequence of operators $(T_s)_{s\in\mathbb{N}}$
        \begin{equation*}
\mathcal{R}:=\bigcap_{s\in\mathbb{N}}\mathcal{R}_{T_s}
        \end{equation*}
        is a rich family of closed separable subspaces of $L_1(X, \Sigma, \mu)$, where $\mathcal{R}_{T_s}$ is defined as in \ref{it: rich-families-a1}.
    \end{enumerate}
\end{lemma}
\begin{proof}
    \ref{it: rich-families-a1} We will first show that condition \ref{it: def-rich-a} of Definition \ref{def: richFamily} holds. For that, let $M\subseteq L_1(X,\Sigma, \mu)$ be a separable set. By Proposition \ref{prop: separable-in-separable-sublattice}, we can find a countably generated $\sigma$-algebra $\Sigma_1\subseteq \Sigma$, such that the sublattice $E_1 := L_1(X, \Sigma_1, \mu) \subseteq L_1(X, \Sigma, \mu)$ contains $M$. Let $(f_n)_{n \in \mathbb{N}}$ be a norm dense subset in $E_1^+$. Using the Riesz-Kantorovich formula \eqref{eq: Riesz-Kantorovich}, we have
    \begin{equation*}
        |T|f_n = \sup \{T g: -f_n \leq g \leq f_n, g \in L_1(X, \Sigma, \mu) \}.
    \end{equation*}
    By Proposition \ref{prop: countable-supremum-property},  $L_1(X,\Sigma,\mu)$ has the countable sup property and therefore, for each $n \in \mathbb{N}$, we can find a sequence $(g_{n}^{(r)})_{r \in \mathbb{N}} \subseteq L_1(X, \Sigma, \mu)$ such that
    \begin{equation*}
        |T|f_n = \sup_{r \in \mathbb{N}} \{T g_{n}^{(r)}: -f_n \leq g_{n}^{(r)} \leq f_n \}.
    \end{equation*}
   Again, by Proposition \ref{prop: separable-in-separable-sublattice}, we can find a countably generated $\sigma$-algebra $\Sigma_2$ with $\Sigma_1\subseteq\Sigma_2$, so that $E_2 = L_1(X, \Sigma_2, \mu)$ satisfies 
    \begin{equation*}
        \{g_{n}^{(r)}: n,r \in \mathbb{N} \} \cup \{f_n: n \in \mathbb{N} \}\subseteq E_2,
    \end{equation*}
    and thus
    \begin{equation*}
        |T|f_n = \sup_{r \in \mathbb{N}} \{Tg_{n}^{(r)}: -f_n \leq g_{n}^{(r)} \leq f_n \} \leq \sup\{T h: -f_n \leq h \leq f_n, \hspace{5pt} h \in L_1(X, \Sigma_2, \mu) \}.
    \end{equation*}
    It follows that for every $n \in \mathbb{N}$ we have
    \begin{equation*}
        |T| f_n = |T J_{E_2}|f_n,
    \end{equation*}
    and since the sequence $(f_n)_{n \in \mathbb{N}}$ is dense in $E_1^+$ we get
    \begin{equation*}
        |T| f = |T J_{E_2}|f
    \end{equation*}
    for all $f \in E_1$.

    Arguing inductively, we can construct an increasing sequence $(\Sigma_m)_{m \in \mathbb{N}}$ of countably generated $\sigma$-algebras such that the sublattices $(E_m)_{m \in \mathbb{N}} = (L_1(X, \Sigma_m, \mu))_{m \in \mathbb{N}}$ satisfy
    \begin{equation*}
        M \subseteq E_1 \subseteq E_2 \subseteq E_3 \subseteq \ldots \subseteq E_{m-1} \subseteq E_m \subseteq \ldots
    \end{equation*}
    and
    \begin{equation*}
        |T| f = |T J_{E_m}|f
    \end{equation*}
    for all $f \in E_{m-1}$.

    Observe that if we let $\tilde{\Sigma}$ be the smallest $\sigma$-algebra on $X$ containing $\bigcup_{m \in \mathbb{N}} \Sigma_m$ then $\tilde{\Sigma}$ is countably generated. Let $E = L_1(X, \tilde{\Sigma}, \mu)$ and note that $\bigcup_{m \in \mathbb{N}} E_m$ is dense in $E$. For each $m \in \mathbb{N}$ and $f \in E_m^+$, it follows that
    \begin{equation*}
     |TJ_E| f \geq |TJ_{E_{m+1}}| f = |T| f.
    \end{equation*} 
    Since this holds for any $f \in E^+_{m}$, it follows that $|TJ_E|f = |T|f$ for all $f \in E_m$. As the union of the $E_m$ is dense in $E$, we get that $|TJ_E|f = |T|f$ for all $f \in E$; in other words, $|TJ_E| = |T|J_E$. By construction $M \subseteq E = L_1(X, \tilde{\Sigma}, \mu)$, so that $\mathcal{R}_T$ satisfies condition~\ref{it: def-rich-a} of Definition~\ref{def: richFamily}.
    
    We now prove condition~\ref{it: def-rich-b} of Definition~\ref{def: richFamily}. Let $(G_n)_{n \in \mathbb{N}} := (L_1(X, \tilde{\Sigma}_n, \mu))_{n \in \mathbb{N}} \subseteq \mathcal{R}_T$ be an increasing sequence, and let
    \begin{equation*}
        G = \overline{\bigcup_{n\in\mathbb{N}} G_n}.
    \end{equation*}
    Since $(G_n)_{n \in \mathbb{N}}$ is increasing, it is easy to see that $G = L_1(X, \hat{\Sigma}, \mu)$, where $\hat{\Sigma}$ is the smallest $\sigma$-algebra containing $\bigcup_{n \in \mathbb{N}} \tilde{\Sigma}_n$; in particular, $\hat{\Sigma}$ is countably generated. Therefore, to verify that $G \in \mathcal{R}_T$, we only need to check that
    \begin{equation*}
        |TJ_G| = |T| J_G.
    \end{equation*}
    This argument is identical to the one previously done for $E$, this time using that $\bigcup_{n \in \mathbb{N}} G_n$ is dense in $G$, together with the fact that
    \begin{equation*}
        |T|f = |TJ_{G_n}|f
    \end{equation*}
    for all $f \in G_{n}$ and $G_{n} \in \mathcal{R}_T$. Details are omitted.
    
\ref{it: rich-families-a2} Follows immediately from Proposition \ref{prop: IntersectionRFamily}.
\end{proof}

Using the previous lemma, we can now show that it is possible to pass to a separable sublattice $E$ such that the restriction of the operator $S$ remains singular.

\begin{lemma}\label{lmm: technical2}
    Let $(X, \Sigma, \mu)$ and $(Y, \Gamma, \nu)$ be measure spaces of finite measure and let $S: L_1(X, \Sigma, \mu) \to L_1(Y, \Gamma, \nu)$ be a singular operator. Then 
    \begin{align*}
    \mathcal{R}^S = \{E\times F: & \hspace{3pt} E = L_{1}(X,\tilde{\Sigma},\mu) \text{ and } F = L_{1}(Y,\tilde{\Gamma},\nu) \text{ for some } \tilde{\Sigma}\subseteq  \Sigma, \tilde{\Gamma}\subseteq \Gamma\\
    &\text{ countably generated } \sigma\text{-algebras, } S(E)\subseteq F, \\
    &\text{ and } S:E\rightarrow F \text{ is singular }\}  
    \end{align*}
    is a rich family of closed separable subspaces of $L_1(X,\Sigma,\mu)\times L_1(Y,\Gamma,\nu)$.
\end{lemma}
\begin{proof} 

We start by proving part~\ref{it: def-rich-a} of Definition~\ref{def: richFamily}. Let $A$ be a separable subset of $L_1(X,\Sigma,\mu)\times L_1(Y,\Gamma,\nu)$. Without loss of generality, we may assume that $A = M \times N$ for some separable sets $M \subseteq L_1(X,\Sigma,\mu)$ and $N \subseteq L_1(Y,\Gamma,\nu)$. Indeed, let $(f_n, g_n)_{n\in\mathbb{N}}$ be a dense sequence in $A$, and let $M$ and $N$ be the sublattices generated by the $f_n$ and $g_n$, respectively. Then $A \subseteq M \times N$. Therefore, to verify part~\ref{it: def-rich-a}, it suffices to show that for every pair of separable subsets $M \subseteq L_1(X,\Sigma,\mu)$ and $N \subseteq L_1(Y,\Gamma,\nu)$, there exists $E \times F \in \mathcal{R}^S$ such that $M \times N \subseteq E \times F$. \\

By Proposition \ref{prop: separable-in-separable-sublattice}, we can find a countably generated $\sigma$-algebra $\Sigma_1\subseteq\Sigma$ such that the sublattice $E_1 = L_1(X, \Sigma_1, \mu) \subseteq L_1(X, \Sigma, \mu)$ contains $M$. Similarly, we can find a countably generated $\sigma$-algebra  $\Gamma_1 \subseteq \Gamma$ such that the sublattice $F_1 = L_1(Y, \Gamma_1, \nu) \subseteq L_1(Y, \Gamma, \nu)$ satisfies $S[E_1]+ N \subseteq F_1$.

   Observe that if $K: E_1 \to F_1$ is an integral operator with kernel $k: X \times Y \to \mathbb{R}$, then Proposition \ref{prop: KernelinL1} guarantees that $k \in G_1 = L_1(X \times Y, \Sigma_1 \otimes \Gamma_1, \mu \otimes \nu)$. Find a dense sequence $(k^{(m)})_{m \in \mathbb{N}}$ of simple functions in $G_1$. We note that, since simple functions are essentially bounded, $k^{(m)}$ is the kernel of an integral operator $K^{(m)}: L_1(X, \Sigma, \mu) \to L_1(Y, \Gamma, \nu)$. 

   Define a sequence of operators $V_s: L_1(X, \Sigma, \mu) \to L_1(Y, \Gamma, \nu)$ by
   \begin{equation*}
       V_{2m - 1} = S - K^{(m)} \text{ for } s = 2m-1 \text{ and } V_{2m} = S + K^{(m)} \text{ for } s = 2m.
   \end{equation*}
   By Lemma \ref{lmm: technical1} \ref{it: rich-families-a2} we can find a countably generated $\sigma$-algebra $\tilde{\Sigma}_1 \subseteq \Sigma$ such that the sublattice $\tilde{E}_1 = L_1(X, \tilde{\Sigma}_1, \mu)$ contains $E_1$ and 
   \begin{equation*}
       |V_s J_{\tilde{E}_1}| = |V_s|J_{\tilde{E}_1}
   \end{equation*}
   for all $s \in \mathbb{N}$. In other words
   \begin{equation}\label{eq: a1}
       |SJ_{\tilde{E}_{1}} \pm K^{(m)}J_{\tilde{E}_{1}}|=|S \pm K^{(m)}|J_{\tilde{E}_{1}} 
   \end{equation}
   for all $m \in \mathbb{N}$. 
   
    Since $(k^{(m)})_{m \in \mathbb{N}}$ is dense in $G_1$, for any integral operator $K$ with kernel $k \in G_1$ we can find a subsequence of $(k^{(m)})_{m \in \mathbb{N}}$, which for simplicity we still denote by $k^{(m)}$, such that $\norm{k^{(m)} - k} \to 0$. Fix $f \in \tilde{E}_1^+$ an essentially bounded function. We have
    \begin{align*}
        \big|\hspace{3pt}|SJ_{\tilde{E}_1} + K^{(m)} J_{\tilde{E}_1}|f - |SJ_{\tilde{E}_1} + K J_{\tilde{E}_1}|f\hspace{3pt} \big| &\leq  |SJ_{\tilde{E}_1} + K^{(m)} J_{\tilde{E}_1} - ( SJ_{\tilde{E}_1} + K J_{\tilde{E}_1})|f\hspace{3pt} \\
        &= \hspace{3pt}|K^{(m)} J_{\tilde{E}_1} - K J_{\tilde{E}_1}|f,
    \end{align*}
    and thus
    \begin{align}\label{eq: a3}
        \big\lVert \hspace{3pt}|SJ_{\tilde{E}_1} + K^{(m)} J_{\tilde{E}_1}|f - &|SJ_{\tilde{E}_1} + K J_{\tilde{E}_1}|f\hspace{3pt}  \big\rVert \leq \hspace{3pt} \big\lVert |K^{(m)}J_{\tilde{E}_1}- K J_{\tilde{E}_1}|f \hspace{3pt} \big\rVert
         \notag \\
         &\leq \int_Y \int_{X} |k^{(m)}(x,y) - k(x,y)| \hspace{2pt} f(x) d\mu(x) d\nu(y) \to 0,
    \end{align}
    since $\norm{k^{(m)} - k} \to 0$ and $f$ is essentially bounded. 
    
    Combining \eqref{eq: a1} and \eqref{eq: a3}, we obtain
   \begin{align}
       |S + K| f \geq |SJ_{\tilde{E}_{1}}  + KJ_{\tilde{E}_{1}}|f &\stackrel{\eqref{eq: a3}}{=} \lim_{m \to \infty} |SJ_{\tilde{E}_{1}} + K^{(m)}J_{\tilde{E}_{1}}|f \notag \\
       &\stackrel{\eqref{eq: a1}}{=} \lim_{m \to \infty}|S + K^{(m)}|f = |S + K|f,
   \end{align}
   where the last limit exists using the same argument as before.
   Since this is true for any essentially bounded positive function, we obtain
   \begin{equation*}
       |SJ_{\tilde{E}_{1}} + KJ_{\tilde{E}_{1}}| = |S + K|J_{\tilde{E}_{1}};
   \end{equation*}
   an identical argument shows $|SJ_{\tilde{E}_{1}} - KJ_{\tilde{E}_{1}}| = |S - K|J_{\tilde{E}_{1}}$.
   
   Therefore, $\eqref{eq: a1}$ holds with $K^{(m)}$ replaced by any integral operator $K$ with kernel in $G_1$. We can now find a countably generated $\sigma$-algebra $\tilde{\Gamma}_1 \subseteq \Gamma$ such that the sublattice $\tilde{F_1} = L_1(Y, \tilde{\Gamma}_1, \nu)$ satisfies $S[\tilde{E}_1] + F_1 \subseteq \tilde{F}_1$. 

   Arguing in the same manner and using Proposition \ref{prop: separable-in-separable-sublattice}, we can find a countably generated $\sigma$-algebra $\Sigma_2 \subseteq \Sigma$ such that $\tilde{E}_1 \subseteq L_1(X, \Sigma_2, \mu) = E_2$ and
   \begin{equation*}
       |SJ_{E_2} \pm K J_{E_2}| = |S \pm K|J_{E_2}
   \end{equation*}
    for all integral operators with kernel in $\tilde{G}_1 = L_1(X \times Y, \tilde{\Sigma}_1 \otimes \tilde{\Gamma}_1, \mu \otimes \nu)$. As before, we can choose a countably generated $\sigma$-algebra $\Gamma_2 \subseteq \Gamma$ such that the sublattice $F_2 = L_1(Y, \Gamma_2, \nu)$ satisfies $S[E_2] + \tilde{F}_1 \subseteq F_2$.

    Inductively, we can construct a sequence of increasing countably generated $\sigma$-algebras $(\Sigma_m)_{m \in \mathbb{N}}$, $(\tilde{\Sigma}_m)_{m \in \mathbb{N}}$, $( \Gamma_m)_{m \in \mathbb{N}}$ and $(\tilde{\Gamma}_m)_{m \in \mathbb{N}}$, each contained in $\Sigma$ and $\Gamma$ respectively, such that
    \begin{enumerate}[label=(\alph*)]
        \item $M \subseteq E_1 \subseteq \tilde{E}_1 \subseteq E_2 \subseteq \tilde{E}_2 \subseteq \ldots$,
        \item $N\subseteq F_1 \subseteq \tilde{F}_1 \subseteq F_2 \subseteq \tilde{F}_2 \subseteq\ldots$,
        \item $G_1 \subseteq \tilde{G}_1 \subseteq G_2 \subseteq \tilde{G}_2 \subseteq \ldots$,
        \item $S[E_{2m}] + \tilde{F}_{2m-1} \subseteq F_{2m}$ for all $m \geq 1$,
        \item \label{it: technical1} $|SJ_{E_{2m}} \pm KJ_{E_{2m}}| = |S \pm K|J_{E_{2m}}$ for all integral operators $K$ with kernel $k \in \tilde{G}_{2m-1}$ and all $m \geq 1$,
    \end{enumerate}
    where
    \begin{equation*}
        \begin{aligned}
            E_m &= L_1(X, \Sigma_m, \mu), \quad & \tilde{E}_m &= L_1(X, \tilde{\Sigma}_m, \mu), \\
            F_m &= L_1(Y, \Gamma_m, \nu), \quad & \tilde{F}_m &= L_1(Y, \tilde{\Gamma}_m, \nu), \\
            G_m &= L_1(X \times Y, \Sigma_m \otimes \Gamma_m, \mu \otimes \nu), \quad &
            \tilde{G}_m &= L_1(X \times Y, \tilde{\Sigma}_m \otimes \tilde{\Gamma}_m, \mu \otimes \nu).
        \end{aligned}
    \end{equation*}

    Let $\tilde{\Sigma} \subseteq \Sigma$ be the $\sigma$-algebra generated by the $\Sigma_m$ and $\tilde{\Gamma} \subseteq \Gamma$ the $\sigma$-algebra generated by the $\Gamma_m$, and observe that both $\tilde{\Sigma}$ and $\tilde{\Gamma}$ are countably generated $\sigma$-algebras. Define $E = L_1(X, \tilde{\Sigma}, \mu)$ and $F = L_1(Y, \tilde{\Gamma}, \nu)$.

    By construction, $M \subseteq E$ and $N \subseteq F$ so that $M \times N \subseteq E \times F$. Therefore, to finish, we need to show that $E \times F \in \mathcal{R}^S$.
    
    Note that that since $\bigcup_{m \in \mathbb{N}} E_m$ is dense in $E$ and $S[E_{2m}] \subseteq F_{2m}$, it is easy to check that $S[E] \subseteq F$. Therefore, to show that $E \times F \in \mathcal{R}^S$, we only need to verify that $S\big|_E: E \to F$ is still a singular operator.

    Define $G = L_1(X \times Y, \tilde{\Sigma} \otimes \tilde{\Gamma}, \mu \otimes \nu)$ and let $H: E \to F$ be an integral operator with kernel $h \in G$. Observe that $\tilde{G}_{m} \subseteq \tilde{G}_{m+1}$ and $\bigcup_{m \in \mathbb{N}} \tilde{G}_m$ is norm dense in $G$, and thus we can find a sequence of simple functions $h^{(m)} \in \tilde{G}_{2m-1}$ such that $\norm{h^{(m)} - h} \to 0$. Denote by $H^{(m)}: E \to F$ the integral operators with kernel $h^{(m)}$.

    Fix $n \in \mathbb{N}$ and an essentially bounded positive function $f \in E^+_n$. Arguing as in the proof of \eqref{eq: a3} we have
    \begin{align*}
        \big\lVert \hspace{3pt}|SJ_E - H^{(m)}J_E|f - |SJ_E - HJ_E|f \hspace{3pt}  \big\rVert \to 0, \hspace{5pt} \text{ and } \norm{|S - H^{(m)}|f - |S - H|f} \to 0
    \end{align*}
    since $\norm{h^{(m)} - h} \to 0$ and $f$ is essentially bounded.
    
    Thus, for all $m$ with $2m \geq n$ we get
    \begin{align*}
        |S - H|f &\geq |SJ_E - HJ_E|f = \lim_{m \to \infty} |SJ_{E} - H^{(m)} J_E|f \\
        &\geq   \lim_{m \to \infty} |SJ_{E_{2m}} - H^{(m)} J_{E_{2m}}|f \overset{\ref{it: technical1}}{=} \lim_{m \to \infty}|S - H^{(m)}|f = |S - H|f,
    \end{align*}
    that is $|S-H|f = |SJ_E - HJ_E|f = |SJ_E - H|f$. A similar argument also gives $|SJ_E + H|f = |S + H|f$. 

    Since $S$ is singular and $H$ is an integral operator, we get
    \begin{equation*}
        |SJ_E - H|f = |S - H|f = |S + H|f = |SJ_E + H|f,
    \end{equation*}
    for all essentially bounded functions $f \in E_n$. Since $n \in \mathbb{N}$ was arbitrary and $\bigcup_{n \in \mathbb{N}} E_n$ is dense in $E$, the above equation holds for all $f \in E$. In other words, we have $|SJ_E| \wedge |H| = 0$ for all integral operators $H$ with kernels in $G$, which, by \eqref{eq: singular-if-only-if-orthogonal}, is equivalent to $SJ_E$ being a singular operator. This shows that $SJ_E$ is singular, so $E \times F \in \mathcal{R}^S$, and hence condition \ref{it: def-rich-a} of Definition \ref{def: richFamily} is satisfied.

    To prove part \ref{it: def-rich-b} of Definition \ref{def: richFamily}, consider an increasing sequence
    \begin{equation*}
        (E_n \times F_n)_{n \in \mathbb{N}} := \big(L_1(X, \Sigma_n, \nu) \times L_1(Y, \Gamma_n, \nu)\big)_{n \in \mathbb{N}} \subseteq \mathcal{R}^S,
    \end{equation*}
    where the reader should be careful not to confuse $\Sigma_n$, $\Gamma_n$, $E_n$, and $F_n$ here with those used in the construction above. We aim to show that
    \begin{equation*}
        \overline{\bigcup_{n \in \mathbb{N}} E_n \times F_n } = \overline{ \bigcup_{n \in \mathbb{N}} E_n } \times \overline{ \bigcup_{n \in \mathbb{N}} F_n } \in \mathcal{R}^S.
    \end{equation*}
    
    Define
    \begin{equation*}
        E = \overline{ \bigcup_{n \in \mathbb{N}} E_n } \quad \text{and} \quad F = \overline{ \bigcup_{n \in \mathbb{N}} F_n },
    \end{equation*}
    so that the previous condition can simply be expressed as $E \times F \in \mathcal{R}^S$.
    
    Observe that since $(E_n \times F_n)_{n \in \mathbb{N}}$ is an increasing family, the same holds for $(E_n)_{n \in \mathbb{N}}$ and $(F_n)_{n \in \mathbb{N}}$. In particular, it is easy to see that this implies
    \begin{equation*}
        E = L_1(X, \tilde{\Sigma}, \mu) \quad \text{and} \quad F = L_1(Y, \tilde{\Gamma}, \nu),
    \end{equation*}
    where $\tilde{\Sigma} \subseteq \Sigma$ and $\tilde{\Gamma} \subseteq \Gamma$ are the smallest $\sigma$-algebras containing $\bigcup_{n \in \mathbb{N}} \Sigma_n$ and $\bigcup_{n \in \mathbb{N}} \Gamma_n$, respectively.
    
    Since $\bigcup_{n \in \mathbb{N}} E_n$ is dense in $E$ and $S(E_n) \subseteq F_n$, it is clear that $S(E) \subseteq F$.
    To see that $E \times F \in \mathcal{R}^S$, we only need to check that $ S\big|_{E}: E \to F \text{ is singular}$.
    This follows similarly as before, using the density of $\bigcup_{n \in \mathbb{N}} E_n$ in $E$ together with the fact that $E_n \times F_n \in \mathcal{R}^S$. Since the arguments are similar, and condition \ref{it: def-rich-b} will not be needed in future results, we omit the details.
\end{proof}

The previous lemma allows us to reduce the analysis of singular operators on spaces of finite measure to the setting of measure spaces with a countably generated $\sigma$-algebra, as we show now.

\begin{proposition}\label{prop: singular-operators-finite-case}
    Let $(X, \Sigma, \mu)$ and $(Y, \Gamma, \nu)$ be measure spaces of finite measure and $S: L_1(X, \Sigma, \mu) \to L_1(Y, \Gamma, \nu)$ be a singular operator. Then we have
    \begin{equation*}
        \Delta(S) = \norm{S} = \norm{S}_w.
    \end{equation*}
\end{proposition}
\begin{proof}
    It is enough to show $\Delta(S) \geq \norm{S}$, since the other inequality is always true. Let $(f_n)_{n \in \mathbb{N}}$ be a sequence of norm one functions such that $\norm{Sf_n} \to \norm{S}$ and define $M = \spn \{ f_n: n \in \mathbb{N}\}$. By Lemma \ref{lmm: technical2} the family $\mathcal{R}^S$ is a rich family of closed separable subspaces of $L_1(X, \Sigma, \mu) \times L_1(Y, \Gamma, \nu)$ and thus by Definition \ref{def: richFamily} \ref{it: def-rich-a} applied to the separable set $M\times\{0\}$, we can find countably generated $\sigma$-algebras $\tilde{\Sigma}$ and $\tilde{\Gamma}$ such that $M \subseteq L_1(X, \tilde{\Sigma}, \mu) =: E$, $S[E] \subseteq L_1(Y, \tilde{\Gamma}, \nu) =: F$ and $S\big|_E: E \to F$ is singular. We note that $\norm{SJ_E} = \norm{S}$ since $M \subseteq E$.
    Theorem \ref{th: countably-generated-finite-case} applied to $SJ_E: E \to F$ gives $\Delta(S J_E) = \norm{SJ_E}$ and therefore
    \begin{equation*}
        \Delta(S) \geq \Delta(SJ_E) = \norm{SJ_E} = \norm{S},
    \end{equation*}
    finishing the proof.
\end{proof}

We conclude by describing how to handle the integral part of the operator when working with general measure spaces of finite measure.

\begin{lemma}\label{lmm: kernel-operator-finite-case}
    Let $(X, \Sigma, \mu)$ and $(Y, \Gamma, \nu)$ be measure spaces of finite measure and $K: L_1(X, \Sigma, \mu) \to L_1(Y, \Gamma, \nu)$ be an integral operator with kernel $k$. Let $\varepsilon>0$, then for each $n\in \mathbb{N}$ we can find $G_{n}\subseteq Y$ with $\nu(G_{n})\leq\frac{2\rVert K \lVert}{n}$ and $M_{n}\in\Sigma$ of positive measure such that
    \begin{equation*}
      \int_{G_n}|k(x,y)|d\nu(y)\geq\int_{Y}|k_{n}(x,y)|d\nu(y)-\varepsilon \text{ for all } x\in M_{n},
    \end{equation*}
    where $k_{n}$ are the truncated kernels of $k$.
\end{lemma}
\begin{proof}
    Without loss of generality, we can assume that the operator $K$ is positive. We divide the proof into three cases.
    \begin{enumerate}   [label=\emph{Case~\arabic*:}, leftmargin=*]
        \item \textit{(Standard measure spaces)} Suppose first that $(X, \Sigma, \mu)$ and $(Y, \Gamma, \nu)$ are standard measure spaces. In this case, the result was proven in the first proposition of \cite{weis1984approximation}.
        \item \textit{(Measure spaces of finite measure with countably generated $\sigma$-algebra)} Suppose now that $(X, \Sigma, \mu)$ and $(Y, \Gamma, \nu)$ are measure spaces of finite measure with $\Sigma$ and $\Gamma$ countably generated. By Theorem \ref{th: separablefinite-to-borel} we can find standard measure spaces $(\tilde{X}, \tilde{\Sigma}, \tilde{\mu})$ and $(\tilde{Y}, \tilde{\Gamma}, \tilde{\nu})$ and measure algebra isomorphisms $\phi: \tilde{\Sigma}_{\tilde{\mu}} \to \Sigma_\mu$ and $\psi: \Gamma_\nu \to \tilde{\Gamma}_{\tilde{\nu}}$; denote by $\Phi: L_1(\tilde{X}, \tilde{\Sigma}, \tilde{\mu}) \to L_1(X, \Sigma, \mu)$ and $\Psi: L_1(Y, \Gamma, \nu) \to L_1(\tilde{Y}, \tilde{\Gamma}, \tilde{\nu})$ the induced isometric lattice isomorphisms.

        Define $\tilde{K} = \Psi K \Phi$, which is an integral operator with truncated integral approximations $\tilde{K}_n = \Psi K_n \Phi$ by Lemma \ref{lmm: TransKTruncatedK}. Fix any $\varepsilon > 0$.

        By the previous case, for each $n \in \mathbb{N}$, we can find $\tilde{G}_n \subseteq \tilde{Y}$ with $\tilde{\nu}(\tilde{G}_n) < \frac{2\norm{\tilde{K}}}{n}$ and $\tilde{M}_n \in \tilde{\Sigma}$ of positive measure such that
        \begin{equation*}
            \tilde{K}^* \mathds{1}_{\tilde{G}_n}(\tilde{x}) = \int_{\tilde{G}_n}\tilde{k}(\tilde{x},\tilde{y})d\tilde{\nu}(\tilde{y})\geq \int_{\tilde{Y}}\tilde{k}_{n}(\tilde{x},\tilde{y})d\tilde{\nu}(\tilde{y})-\varepsilon = \tilde{K}^* \mathds{1}_{\tilde{Y}}(\tilde{x})- \varepsilon,
        \end{equation*}
        for all $\tilde{x} \in \tilde{M}_n$. In other words,
        \begin{equation*}
            \chi_{\tilde{M}_n} \tilde{K}^* \mathds{1}_{\tilde{G}_n} \geq \chi_{\tilde{M}_n}\tilde{K}_n^* \mathds{1}_{\tilde{Y}} - \varepsilon\mathds{1}_{\tilde{M}_n},
        \end{equation*}
        equivalently
        \begin{equation*}
            \chi_{\tilde{M}_n} \Phi^* K^* \Psi^* \mathds{1}_{\tilde{G}_n} \geq \chi_{\tilde{M}_n} \Phi^* K^*_n \Psi^* \mathds{1}_{\tilde{Y}} - \varepsilon\mathds{1}_{\tilde{M}_n}.
        \end{equation*}
        Using the definition of $\Psi$, we get
        \begin{equation*}
            \chi_{\tilde{M}_n} \Phi^* K^*  \mathds{1}_{\psi^{-1}(\tilde{G}_n)} \geq \chi_{\tilde{M}_n} \Phi^* K^*_n  \mathds{1}_{Y} - \varepsilon\mathds{1}_{\tilde{M}_n},
        \end{equation*}
        and since $\Phi$, and thus $\Phi^*$, is a lattice isomorphism, this is equivalent to
        \begin{equation*}
            (\Phi^*)^{-1}\chi_{\tilde{M}_n} \Phi^* K^*  \mathds{1}_{\psi^{-1}(\tilde{G}_n)} \geq (\Phi^*)^{-1}\chi_{\tilde{M}_n} \Phi^* K^*_n  \mathds{1}_{Y} - \varepsilon(\Phi^*)^{-1}\mathds{1}_{\tilde{M}_n}.
        \end{equation*}
        Since $\Phi$ satisfies $\chi_{\tilde{M}_n}\Phi^* = \Phi^* \chi_{\phi(\tilde{M}_n)}$ and $(\Phi^*)^{-1}\mathds{1}_{\tilde{M}_n} = \mathds{1}_{\phi(\tilde{M}_n)}$, we get
        \begin{equation*}
            \chi_{\phi(\tilde{M}_n)} K^*  \mathds{1}_{\psi^{-1}(\tilde{G}_n)} \geq \chi_{\phi(\tilde{M}_n)}  K^*_n \mathds{1}_{Y} - \varepsilon\mathds{1}_{\phi(\tilde{M}_n)}.
        \end{equation*}
        Setting $M_n = \phi(\tilde{M}_n)$ and $G_n = \psi^{-1}(\tilde{G}_n)$ this gives
        \begin{equation*}
            \chi_{M_n} K^*  \mathds{1}_{G_n} \geq \chi_{M_n} K^*_n \mathds{1}_{Y} - \varepsilon \mathds{1}_{M_n}.
        \end{equation*}
        In other words, we have
        \begin{equation*}
        \int_{G_n} k(x,y) d\nu(y) = K^* \mathds{1}_{G_n}(x) \geq K_n^* \mathds{1}_Y - \varepsilon = \int_{Y} k_n(x,y) d\nu(y) - \varepsilon
        \end{equation*}
        for all $x \in M_n$.
        
        Since $\tilde{M}_n$ has positive measure, then so does $M_n$, while
        \begin{equation*}
            \nu(G_n) = \nu(\psi^{-1}(\tilde{G}_n)) = \tilde{\nu}(\tilde{G}_n) < \frac{2 \norm{\tilde{K}}}{n} = \frac{2 \norm{K}}{n},
        \end{equation*}
        which finishes this case.
        
        \item \textit{(Finite measure)} Finally, consider the case when $(X, \Sigma, \mu)$ and $(Y, \Gamma, \nu)$ are measure spaces of finite measure. By Proposition \ref{prop: KernelinL1}, we have that $k \in L_1(X \times Y, \Sigma \otimes \Gamma, \mu \otimes \nu)$. Let $k^{(m)}$ be a sequence of simple functions defined on rectangles,
        \begin{equation*}
            k^{(m)} = \sum_{l=1}^{N_m} a_{m, l} \mathds{1}_{A_{m,l} \times B_{m,l}},
        \end{equation*}
        such that $\norm{k^{(m)} - k} \to 0$. 
       Let $\tilde{\Sigma}$ be the $\sigma$-algebra generated by $\{A_{m,l}: m, l \in \mathbb{N} \}$ and $\tilde{\Gamma}$ be the $\sigma$-algebra generated by $\{B_{m,l}: m, l \in \mathbb{N} \}$. Observe that $\tilde{\Sigma}$ and $\tilde{\Gamma}$ are countably generated and $k \in L_1(X \times Y, \tilde{\Sigma} \otimes \tilde{\Gamma}, \mu \otimes \nu)$. Now, $K$ can be seen as an integral operator $K: L_1(X, \tilde{\Sigma}, \mu) \to L_1(Y, \tilde{\Gamma}, \nu)$. By the previous case, for any $\varepsilon > 0$ and any $n \in \mathbb{N}$ we can find measurable sets $G_n  \subseteq Y$ with $\nu(G_n) < \frac{2\norm{K}}{n}$ and $M_n  \subseteq X$ of positive measure such that
       \begin{equation*}
           \int_{G_n} k(x,y) d\nu(y) \geq \int_{Y} k_n(x,y) d\nu(y) - \varepsilon
       \end{equation*}
       for all $x \in M_n$. This finishes the proof.  \qedhere
    \end{enumerate}
\end{proof}

\subsection{Strictly localizable measure spaces.} \label{subsection: strictly-localizable}

We are finally ready to move to the general case, where the underlying measure spaces are strictly localizable. Recall that if $(X, \Sigma, \mu)$ is strictly localizable then we have a decomposition $((X_i, \Sigma_i, \mu))_{i \in I}$ so that each $(X_i, \Sigma_i, \mu)$ has finite measure and
\begin{equation*}
    L_1(X, \Sigma, \mu) = \left(\bigoplus_{i \in I} L_1(X_i, \Sigma_i, \mu)\right)_{\ell_1}.
\end{equation*}

\begin{rem}
    The previous decomposition as an $\ell_1$-sum is actually valid for any measure space $(X,\Sigma,\mu)$, if one allows the sets $(X_i)_{i\in I}$ to form an \emph{almost decomposition}, in the sense that
\begin{equation*}
        \mu(X_i \cap X_j) = 0 \quad \text{for all } i \neq j.
\end{equation*}
    Details of such a construction can be found, for example, in \cite[Lemma~7.4]{kacena2013quantitative}. Therefore, even the concrete formulas we prove for measures remain valid for arbitrary measure spaces $(X,\Sigma,\mu)$, with only minor modifications of the proofs. For simplicity of presentation, we state and prove the results in the case of strictly localizable measure spaces. We thank the referee for pointing out that the results can naturally be extended to arbitrary measure spaces.
\end{rem}

From the previous decomposition, to study an operator $T$ from $L_1(X, \Sigma, \mu)$, it is enough to study the operators $T\big|_{L_1(X_i, \Sigma_i, \mu)}$ for each $i \in I$. This motivates us to state the following result, which is essentially well known.

\begin{lemma}\label{lmm: finite-to-stricly-localizable}
    Let $(X, \Sigma, \mu)$ be a measure space of finite measure and $(Y, \Gamma, \nu)$ be a strictly localizable measure space with decomposition $((Y_j, \Sigma_j, \nu))_{j \in J}$. 
    Then the range of any operator $T: L_1(X, \Sigma, \mu) \to L_1(Y, \Gamma, \nu)$ is contained in the direct sum of $L_1(Y_j, \Gamma_j, \nu)$ for at most countably many $j \in J$.
\end{lemma}
\begin{proof}
    We have the decomposition
    \begin{equation*}
        L_1(Y, \Gamma, \nu) = \Bigl(\bigoplus_{j \in J} L_1(Y_j, \Gamma_j, \nu)\Bigr)_{\ell_1}.
    \end{equation*}
    Since $T$ is weakly continuous, the image of any weakly compact set under $T$ is again weakly compact. By \cite[Lemma 7.2 (ii)]{kacena2013quantitative}, this implies that the image of a weakly compact set is contained in a countable sum of the subspaces $L_1(Y_j, \Gamma_j, \nu)$. Finally, since $\mu$ is finite, $L_1(X, \Sigma, \mu)$ is weakly compactly generated, and the result follows.
\end{proof}

Let $(X, \Sigma, \mu)$ and $(Y, \Gamma, \nu)$ be strictly localizable measure spaces and let \break $((X_i, \Sigma_i, \mu))_{i \in I}$ be a decomposition of $(X, \Sigma, \mu)$. The previous lemma shows that the range of the restricted operator  $T|_{L_{1}(X_{i},\Sigma_{i},\mu)}$ is contained in some $\sigma$-finite $L_1$-sublattice of $L_{1}(Y,\Gamma,\nu)$. In particular, Theorem \ref{th: Decomposition-integral-and-singular} gives a decomposition $T|_{L_{1}(X_{i},\Sigma_{i},\mu)} = S^{(i)} + K^{(i)}$ for some singular operator $S^{(i)}$ and integral operator $K^{(i)}$ so that $|S^{(i)}|\wedge |K^{(i)}|=0$. In other words, for each $f \in L_1(X_i, \Sigma_i, \mu)$ we have
\begin{equation*}
    T|_{L_{1}(X_{i},\Sigma_{i},\mu)}f = S^{(i)}f + K^{(i)}f = S^{(i)}f + \int_{X_{i}}k^{(i)}(x,y)f(x)d\mu(x),
\end{equation*}
where $k^{(i)}: X_i \times Y \to \mathbb{R}$ is the kernel of the integral operator $K^{(i)}$. Define the operators $S, K: L_1(X,\Sigma,\mu)\rightarrow L_{1}(Y,\Gamma,\nu)$ by $Sf:=\sum_{i\in  I}S^{(i)}P_{i}f$ and $Kf:=\sum_{i\in  I}K^{(i)}P_{i}f$, where $P_i: L_1(X, \Sigma, \mu) \to L_1(X_i, \Sigma_i, \mu)$ is the canonical projection. In passing, we note that the absolute convergence of $\sum_i P_i f$, together with the facts that absolute convergence is preserved by $T$ and that $|S^{(i)} P_i| \leq |T|$ and $|K^{(i)} P_i| \leq |T|$ for all $i \in I$, implies the absolute convergence of the preceding series defining $S$ and $K$. It is clear that $T=S+K$ and $(Kf)(y)=\int_{X}k(x,y)f(x)d\mu(x)$, where $k(x,y)=\sum_{i\in I}k^{(i)}(x,y)$.\\

To get a formula for the norm of $T$, we will first determine the value  $|K|^*\mathds{1}_{Y}(x)$. For any $x\in X$ we can find $i \in I$ such that $x\in X_i$. By Lemma \ref{lmm: finite-to-stricly-localizable} the image of $K^{(i)}$ is contained in a sublattice $L_1(\Theta, \gamma, \nu)\subseteq L_1(Y,\Gamma,\nu)$ for some $\sigma$-measurable set $\Theta$ and $\gamma = \{A \cap \Theta: A \in \Sigma \}$. It follows that we can use Tonelli's Theorem to determine the adjoint operator of $K^{(i)}$, in particular
\begin{align*}
    |K|^*\mathds{1}_{Y}(x)=
    P_i^{*}|K|^*\mathds{1}_{Y}(x)&
    =|KP_i|^*\mathds{1}_{Y}(x)
    =|K^{(i)}|^*\mathds{1}_{\Theta}(x)\\
    &=\int_{\Theta}|k^{(i)}(x,y)|d\nu(y)
    =\int_{Y}|k(x,y)|d\nu(y),
\end{align*}
where the first equality holds since $x \in X_i$ and the second follows from \eqref{eq: adjoint-and-modulus-commute} together with the fact that $P_i$ is a band projection. The fourth equality uses Tonelli's Theorem, while the last equality holds because $|k(x,y)| = 0$ whenever $y \notin \Theta$ and $|k(x,y)| = |k^{(i)}(x,y)|$ for all $y \in \Theta$, by the definition of $k$ and the fact that $x \in X_i$.

\begin{rem}
    The reader is advised to exercise caution when analyzing this situation, as many standard results do not automatically carry over if one attempts to directly extend Definition~\ref{def: integral} to strictly localizable measure spaces. For instance, in the previous calculation, we used Tonelli's theorem to conclude that the kernels of $K^{(i)}$ and $(K^{(i)})^*$ coincide. This relies on the spaces being $\sigma$-finite, and may fail in a more general setting.
\end{rem}

Since all $S^{(i)}$ are singular, we have $|S^{(i)}| \wedge |K^{(j)}| = 0$ for all $i,j \in I$, and hence $|S| \wedge |K| = 0$. Therefore, using the previous representation of $|K^*| \mathds{1}_{Y}(x)$, we obtain
\begin{align}\label{eq: OperatorOfOneFunction}
|T|^{*}\mathds{1}_{Y}(x) = |S|^{*}\mathds{1}_{Y}(x)+|K|^*\mathds{1}_{Y}(x) 
    =|S|^{*}\mathds{1}_{Y}(x)+\int_{Y}|k(x,y)|d\nu(y),
\end{align}
and \eqref{eq: modulus-and-norm-relation} gives
\begin{align}\label{eq: norm-of-operator}
\rVert T \lVert=\rVert T^{*} \lVert = \rVert|T|^{*} \lVert&=\esssup_{x\in X}\{|T|^{*}\mathds{1}_{Y}(x)\}  \notag\\
    &=\esssup_{x\in X}\{|S|^{*}\mathds{1}_{Y}(x)+\int_{Y}|k(x,y)|d\nu(y)\}.
\end{align}

We will refer to this previous construction when we speak of the decomposition $T = S + K$ into its singular and integral parts, in the context of strictly localizable measure spaces. Similar considerations apply when referring to integral or singular operators, kernels, truncated kernels, and related notions. \\

For any measurable subset $\Theta \subseteq Y$ with finite measure, consider the measure spaces $(\Theta, \gamma, \nu)$ and $(Y \backslash \Theta, \tilde{\gamma}, \nu)$, where $\gamma = \{A \cap \Theta: A \in \Gamma \}$ and $\tilde{\gamma} = \{A \cap (Y \backslash \Theta): A \in \Gamma \}$. Observe that we have a decomposition $L_1(Y, \Gamma, \nu) = L_1(\Theta, \gamma, \nu) \oplus L_1(Y \backslash \Theta, \tilde{\gamma}, \nu)$ so that $L_1(\Theta, \gamma, \nu)$ and $L_1(Y \backslash \Theta, \tilde{\gamma}, \nu)$ are projection bands. We define
\begin{equation*}\label{eq: definition-alpha-theta}
    \alpha_{\Theta}(T):=\limsup_{\nu(B)\rightarrow 0, B\subseteq \Theta} \rVert \chi_{Y \backslash \Theta} T+\chi_{B}T \lVert.
\end{equation*}

\begin{rem}
    The condition that $\nu(\Theta) < \infty$ is crucial to guarantee that $\alpha_\Theta$ is a good measure of weak compactness. Further, observe that if $Y$ has finite measure, then taking $\Theta = Y$ gives $\alpha_\Theta(T) = \Delta(T)$, that is, we recover the modulus of equi-integrability.
\end{rem}

We also wish to define a version of the weak essential norm restricted to a subset $\Theta \subseteq Y$. To this end, we consider weakly compact operators whose range lies in $L_1(\Theta, \gamma, \nu)$. Since this will be needed later, we allow $\Theta$ to be $\sigma$-finite in the following definition. We define
\begin{equation*}
    \norm{T}_{w, \Theta} = \inf \{ \norm{T - W} \hspace{5pt} | W: L_1(X, \Sigma, \mu) \to L_1(\Theta, \gamma, \nu) \text{ weakly compact} \}.
\end{equation*}

Observe that for any weakly compact operator $W$ with range in $L_1(\Theta, \gamma, \nu)$ we have $\alpha_\Theta(T + W) = \alpha_\Theta(T)$ and that $\alpha_\Theta(W) = 0$ if and only if $W$ is a weakly compact operator with range in $L_1(\Theta, \gamma, \nu)$.

\begin{proposition}\label{prop: alpha-to-essential-theta}
    Let $(X, \Sigma, \mu)$ and $(Y, \Gamma, \nu)$ be strictly localizable measure spaces and $T: L_1(X, \Sigma, \mu) \to L_1(Y, \Gamma, \nu)$ be an operator. Then for any measurable subset $\Theta \subseteq Y$ of finite measure, we have
    \begin{equation*}
        \alpha_{\Theta}(T) \leq \norm{T}_{w, \Theta}.
    \end{equation*}
\end{proposition}
\begin{proof}
 It is clear that for any $T: L_1(X, \Sigma, \mu) \to L_1(Y, \Gamma, \nu)$ we have 
     \begin{equation*}
        \alpha_\Theta(T) \leq \norm{T}.
    \end{equation*}
     Hence for any weakly compact operator $W: L_1(\Sigma, X, \mu) \to L_1(\Theta, \gamma, \nu)$  it follows that
    \begin{equation*}
        \alpha_\Theta(T) = \alpha_\Theta(T + W) \leq \norm{T + W}.
    \end{equation*}
    Since this is true for any weakly compact operator $W: L_1(\Sigma, X, \mu) \to L_1(\Theta, \gamma, \nu)$, taking the infimum yields $\alpha_\Theta(T) \leq \norm{T}_{w, \Theta}$.
\end{proof}

We now establish a representation for $\alpha_{\Theta}(T)$. The reader may wish to compare it with the earlier representations of $\Delta(T)$ given in Theorems~\ref{th: Weis-standard-measure-case} and~\ref{th: countably-generated-finite-case}.

\begin{proposition}\label{prop: alpha-representation}
     Let $(X, \Sigma, \mu)$ and $(Y, \Gamma, \nu)$ be strictly localizable measure spaces and $T: L_1(X, \Sigma, \mu) \to L_1(Y, \Gamma, \nu)$ be an operator. Express $T = S + K$ where $S$ is a singular operator and $K$ is an integral operator with kernel $k$. Then for any measurable subset $\Theta \subseteq Y$ of finite measure, we have
     \begin{equation*}
         \alpha_{\Theta}(T)=\inf_{n\in\mathbb{N}}\esssup_{x\in X}\{|S^*|\mathds{1}_{Y}(x)+\int_{Y\backslash \Theta}|k(x,y)|d\nu(y)+\int_{\Theta}|k_{n}(x,y)|d\nu(y)\},
     \end{equation*}
     where $k_n$ are the truncated kernels of the operator $K$.
\end{proposition}
\begin{proof}
    For convenience, let
    \begin{equation*}
         \eta = \inf_{n\in\mathbb{N}}\esssup_{x\in X}\{|S^*|\mathds{1}_{Y}(x)+\int_{Y\backslash \Theta}|k(x,y)|d\nu(y)+\int_{\Theta}|k_{n}(x,y)|d\nu(y)\}.
    \end{equation*}
    We first show that $\alpha_\Theta(T) \leq \eta$. For this, choose a sequence $B_m \subseteq \Theta$ of measurable sets with $\nu(B_m) \to 0$ and $\norm{\chi_{Y \backslash \Theta}T + \chi_{B_m}T} \to \alpha_{\Theta}(T)$. Since $\chi_{Y \backslash \Theta}$ and $\chi_{B_m}$ are band projections, we have
    \begin{equation*}
        |\chi_{Y \backslash \Theta}T + \chi_{B_m} T| = (\chi_{Y \backslash \Theta} + \chi_{B_m})|T|.
    \end{equation*}
    Therefore, by decomposing $T$ into its positive and negative parts and arguing over these, we can assume without loss of generality that $T$ is a positive operator. In particular, $K$ and $S$ are also positive operators.

    We have
    \begin{align*}
        \norm{\chi_{Y \backslash \Theta} T + \chi_{B_m}T} &= \norm{T^*\chi^*_{Y \backslash \Theta} + T^* \chi^*_{B_m}} \\
        &= \norm{S^*\chi^*_{Y \backslash \Theta} + S^* \chi^*_{B_m} + K^*\chi^*_{Y \backslash \Theta} + K^* \chi^*_{B_m}}\\
        &= \esssup_{x \in X} \{S^*\mathds{1}_{Y \backslash \Theta}(x) + S^* \mathds{1}_{B_m}(x) + K^*\mathds{1}_{Y \backslash \Theta}(x) + K^* \mathds{1}_{B_m}(x)\}.
    \end{align*}
    Since $S^*$ is a positive operator, we have \begin{equation*}
        S^* \mathds{1}_{Y \backslash \Theta} + S^* \mathds{1}_{B_m} \leq S^* \mathds{1}_Y
    \end{equation*}
    and thus
    \begin{equation*}
        \norm{\chi_{Y \backslash \Theta} T + \chi_{B_m}T} \leq \esssup_{x \in X} \{S^*\mathds{1}_{Y}(x) + K^*\mathds{1}_{Y \backslash \Theta}(x) + K^* \mathds{1}_{B_m}(x)\}.
    \end{equation*}
    For each fixed $n \in \mathbb{N}$ we have
    \begin{align*}
        K^* \mathds{1}_{B_m}(x) &= \int_{B_m} k(x,y) d \nu(y) \leq \int_{B_m} k_n(x,y) d \nu(y) + n \cdot \nu(B_m) \\
        &\leq \int_{\Theta} k_n(x,y) d \nu(y) + n \cdot \nu(B_m),
    \end{align*}
    and thus
    \begin{equation*}
        \norm{\chi_{Y \backslash \Theta} T + \chi_{B_m}T} \leq  \esssup_{x \in X} \{S^*\mathds{1}_{Y}(x) + \int_{Y \backslash \Theta} k(x,y) d\nu(y) + \int_{\Theta} k_n(x,y) d \nu(y) + n \cdot \nu(B_m)\}.
    \end{equation*}
    Letting $m \to \infty$ gives
    \begin{equation*}
        \alpha_\Theta(T) \leq \esssup_{x \in X} \{S^*\mathds{1}_{Y}(x) + \int_{Y \backslash \Theta} k(x,y) d\nu(y) + \int_{\Theta} k_n(x,y) d \nu(y) \},
    \end{equation*}
    for each $n \in \mathbb{N}$. Taking the infimum over $n \in \mathbb{N}$ gives $\alpha_\Theta(T) \leq \eta$.

    To see the converse inequality, define the operators
    \begin{equation*}
        V_n = S^* \chi^*_{Y \backslash \Theta} + K^* \chi^*_{Y \backslash \Theta} + S^* \chi^*_\Theta + K^*_n \chi^*_\Theta, 
    \end{equation*}
    where $K_n$ are the truncated integral approximations of $K$. Observe that since $K$ is a positive operator, so are each of its truncated integral approximations $K_n$, and thus $V_n$ is positive for each $n \in \mathbb{N}$. Note that we have 
    \begin{equation*}
        \norm{V_n} = \esssup_{x \in X} V_n\mathds{1}_Y(x),
    \end{equation*}
    and
    \begin{equation*}
        \lim_{n \to \infty} \norm{V_n} = \eta.
    \end{equation*}
    Since $(X, \Sigma, \mu)$ is strictly localizable, we can find measurable sets $M_n \subseteq X$ of finite positive measure such that
    \begin{equation*}
        V_n \mathds{1}_{Y}(x) \geq \esssup_{x \in X} V_n\mathds{1}_Y(x) - \frac{1}{n} = \norm{V_n} - \frac{1}{n},
    \end{equation*}
    for all $x \in M_n$. For each $n \in \mathbb{N}$, let $\Sigma_n = \{A \cap M_n: A \in \Sigma\}$ and recall that $J_{L_1(M_n, \Sigma_n, \mu)}$ is the inclusion operator of the projection band $L_1(M_n, \Sigma_n, \mu)$ into $L(X, \Sigma, \mu)$. Consider the operators
    \begin{equation*}
        R_n: L_1(M_n, \Sigma_n, \mu) \to L_1(\Theta, \gamma, \nu)
    \end{equation*}
    given by $R_n = \chi_{\Theta}SJ_{L_1(M_n, \Sigma_n, \mu)}$. We note that since $L_1(M_n, \Sigma_n, \mu)$ and $L_1(\Theta, \gamma, \nu)$ are projection bands, then $R_n$ is still a singular operator. By Proposition \ref{prop: singular-operators-finite-case} we have $\Delta(R_n) = \norm{R_n} = \norm{R_n^*}$. For each $n \in \mathbb{N}$, choose a measurable set $B_n \subseteq \Theta$ such that $\nu(B_n) < 1/n$ and
    \begin{equation*}
        \esssup_{x \in X} R_n^* \mathds{1}_{B_n}(x) = \norm{R_n^* \chi^*_{B_n}} = \norm{\chi_{B_n}R_n} \geq \Delta(R_n) - \frac{1}{n} = \norm{R_n^*} - \frac{1}{n}.
    \end{equation*}
    We can find $\tilde{M}_n\subseteq M_n$ of positive measure such that for all $x \in \tilde{M}_n$ we have
    \begin{equation*}
        R_n^* \mathds{1}_{B_n}(x) \geq \esssup_{x \in X} R_n^* \mathds{1}_{B_n}(x) - \frac{1}{n}.
    \end{equation*}
    Since $\norm{R_n^*} = \esssup_{x \in X} R_n^* \mathds{1}_{\Theta}(x)$, the previous inequalities give
    \begin{equation*}
         R_n^* \mathds{1}_{B_n}(x) \geq \norm{R_n^*} - \frac{2}{n} = R_n^* \mathds{1}_{\Theta}(x) - \frac{2}{n},
    \end{equation*}
    for all $x \in \tilde{M}_n$. This implies that
    \begin{equation*}
        S^* \mathds{1}_{B_n}(x) \geq S^* \mathds{1}_{\Theta}(x) - \frac{2}{n}
    \end{equation*}
    for all $x \in \tilde{M}_n$.

    Again, for each $n \in \mathbb{N}$, let $\tilde{\Sigma}_n = \{A \cap \tilde{M}_n: A \in \Sigma\}$ and $J_{L_1(\tilde{M}_n, \tilde{\Sigma}_n, \mu)}$ be the inclusion operator of the projection band $L_1(\tilde{M}_n, \tilde{\Sigma}_n, \mu)$. Define the integral operator $L_n: L_1(\tilde{M}_n, \tilde{\Sigma}_n, \mu) \to L_1(\Theta, \gamma, \nu)$ by $L_n = \chi_{\Theta} K J_{L_1(\tilde{M}_n, \tilde{\Sigma}_n, \mu)}$. By Lemma \ref{lmm: kernel-operator-finite-case} we can find a measurable set $G_n \subseteq \Theta$ with $\nu(G_n) < \frac{2\norm{L_n}}{n} \leq \frac{2\norm{K}}{n}$ and a measurable set $\hat{M}_n \subseteq \tilde{M}_n$ of positive measure such that
    \begin{equation*}
        \int_{G_n} k(x,y) d\nu(y) \geq \int_{\Theta} k_n(x,y) d\nu(y) - \frac{1}{n}
    \end{equation*}
    for all $x \in \hat{M}_n$. Let $H_n = G_n \cup B_n$ and note that $\lim_{n \to \infty} \nu(H_n) = 0$. 

    For all $x \in \hat{M}_n$ we have
    \begin{align*}
        (T^* \chi^*_{Y \backslash \Theta} + T^* \chi^*_{H_n}) \mathds{1}_Y(x) &= S^* \mathds{1}_{Y \backslash \Theta}(x) + K^* \mathds{1}_{Y \backslash \Theta}(x) + S^* \mathds{1}_{B_n \cup G_n}(x) +  K^* \mathds{1}_{B_n \cup G_n}(x)
        \\
        &= S^* \mathds{1}_{Y \backslash \Theta}(x) + K^* \mathds{1}_{Y \backslash \Theta}(x) + S^* \mathds{1}_{B_n \cup G_n}(x) + \int_{B_n \cup G_n} k(x,y) d\nu(y)
        \\
        &\geq S^* \mathds{1}_{Y \backslash \Theta}(x) + K^* \mathds{1}_{Y \backslash \Theta}(x) + S^* \mathds{1}_{B_n}(x) + \int_{\Theta} k_n(x,y) d\nu(y) - \frac{1}{n}
        \\
        &\geq S^* \mathds{1}_{Y \backslash \Theta}(x) + K^* \mathds{1}_{Y \backslash \Theta}(x) + S^* \mathds{1}_{\Theta}(x) + \int_{\Theta} k_n(x,y) d\nu(y) - \frac{3}{n}
        \\
        &\geq (S^* \mathds{1}_{Y \backslash \Theta}(x) + K^* \mathds{1}_{Y \backslash \Theta}(x) + S^* \mathds{1}_\Theta(x) + K^*_n \mathds{1}_\Theta(x)) - \frac{3}{n}
        \\
        &= V_n \mathds{1}_Y(x) - \frac{3}{n} \geq \norm{V_n} - \frac{1}{n} - \frac{3}{n} = \norm{V_n} - \frac{4}{n}.
    \end{align*}
    Since
    \begin{align*}
        \norm{\chi_{Y \backslash \Theta} T + \chi_{H_n}T} = \norm{T^* \chi^*_{Y \backslash \Theta} + T^*\chi^*_{H_n}} = \esssup_{x \in X} \{(T^* \chi^*_{Y \backslash \Theta} + T^* \chi^*_{H_n}) \mathds{1}_Y(x)\}
    \end{align*}
    and $\hat{M}_n$ has positive measure, we get
    \begin{equation*}
        \norm{\chi_{Y \backslash \Theta} T + \chi_{H_n}T} \geq \norm{V_n} - \frac{4}{n},
    \end{equation*}
    and since $\nu(H_n) \to 0$, it follows that
    \begin{equation*}
        \alpha_\Theta(T) \geq \limsup_{n \to \infty }\norm{\chi_{Y \backslash \Theta} T + \chi_{H_n}T} \geq \lim_{n \to \infty} (\norm{V_n} - \frac{4}{n}) = \eta,
    \end{equation*}
    which proves the other inequality and finishes the proof.
\end{proof}

Using the representation of $\alpha_\Theta(T)$, we are now ready to prove the following.

\begin{theorem}\label{th: main-with-Theta}
    Let $(X, \Sigma, \mu)$ and $(Y, \Gamma, \nu)$ be strictly localizable measure spaces and $T: L_1(X, \Sigma, \mu) \to L_1(Y, \Gamma, \nu)$ be an operator. Then for any measurable subset $\Theta \subseteq Y$ of finite measure, there exists a weakly compact operator $H$ with range in $L_1(\Theta, \gamma, \nu)$ such that
    \begin{equation*}
        \norm{T - H} = \norm{T}_{w, \Theta} = \alpha_\Theta(T).
    \end{equation*}
\end{theorem}
\begin{proof}
    Express $T = S + K$ where $S$ is a singular operator and $K$ is an integral operator with kernel $k$. Consider the measurable function $g: X \to \mathbb{R}$ given by
    \begin{equation*}
        g(x) = |S^*|\mathds{1}_Y(x) + \int_{Y \backslash \Theta} |k(x,y)|d \nu(y)
    \end{equation*}
     and the function $f: X \times [0, \infty] \to \mathbb{R}$ given by
    \begin{equation*}
        f(x,s)= g(x) + \int_{\Theta}(|k(x,y)| - s)^+ d\nu(y).
    \end{equation*}
    Observe that $f$ is a measurable function on $X \times [0,\infty]$, while for fixed $x \in X$, the dominated convergence theorem gives that $f$ is continuous in $s$. Define $A = \{x \in X: |T^*|\mathds{1}_Y(x) > \alpha_{\Theta}(T)\}$, so that by \eqref{eq: OperatorOfOneFunction} we have
    \begin{equation*}
        f(x, 0) = |T^*| \mathds{1}_Y(x) > \alpha_{\Theta}(T)
    \end{equation*}
    for all $x \in A$.

    Observe that if $k_s$ is the truncated kernel of $K$, then $(|k(x,y)| - s)^+ \leq |k_s(x,y)|$. It follows that
    \begin{equation*}
        f(x,s) \leq g(x) + \int_{\Theta} |k_s(x,y)|d\nu(y) \leq \esssup_{x \in X} \{ g(x) + \int_{\Theta} |k_s(x,y)|d\nu(y)\}.
    \end{equation*}
    Recall that by Proposition \ref{prop: alpha-representation} we have
    \begin{align*}
        \alpha_{\Theta}(T)&=\inf_{n\in\mathbb{N}}\esssup_{x\in X}\{|S^*|\mathds{1}_{Y}(x)+\int_{Y\backslash \Theta}|k(x,y)|d\nu(y)+\int_{\Theta}|k_{n}(x,y)|d\nu(y)\}\\
        &= \inf_{n\in\mathbb{N}}\esssup_{x\in X}\{g(x)+\int_{\Theta}|k_{n}(x,y)|d\nu(y)\},
    \end{align*}
    and thus $\lim_{s \to \infty} f(x,s) \leq \alpha_{\Theta}(T)$ for all $x \in A$.
    
    Recall that, for fixed $x \in X$, the function $f(s,x)$ is continuous and decreasing in $s$. Therefore, for every $x \in A$ there exists some $s \in [0,\infty]$ such that $f(x,s) = \alpha_\Theta(T)$. Define the function $s: X \to [0, \infty]$ by
    \begin{equation*}
        s(x) = \begin{cases}
            \displaystyle \inf\{s \in [0,\infty]: f(x,s) = \alpha_\Theta(T)\} & \text{if } x \in A, \\
            0 & \text{if } x \in X \backslash A.
            \end{cases}
    \end{equation*}
    Since $A$ is a measurable set and the function $f$ is measurable, it follows that $s$ is a measurable function. We can also find a measurable unimodular function $\pi: X \times Y \to \mathbb{R}$ such that $k(x,y) = \pi(x,y) |k(x,y)|$ for all $(x,y) \in X \times Y$. Finally define the function $h: X \times Y \to \mathbb{R}$ by
    \begin{equation*}
        h(x,y) = \begin{cases}
            \displaystyle \pi(x,y)[|k(x,y)| \wedge s(x)] & \text{if } (x,y) \in X \times \Theta, \\
            0 & \text{if } (x,y) \in X \times (Y \backslash \Theta).
            \end{cases}
    \end{equation*}
    In particular, note that $h(x,y) = 0$ whenever $x \in X \backslash A$. Furthermore, observe that $h$ is a measurable function, corresponding to the kernel of an integral operator $H: L_1(X, \Sigma, \mu) \to L_1(\Theta, \gamma, \nu) \subseteq L_1(Y, \Gamma, \nu)$, in other words
    \begin{equation*}
        (Hf)(y) = \int_{X} h(x,y) f(x) d\mu(x).
    \end{equation*}
    
    Indeed, that $h$ is measurable is clear from the definition, while it corresponds to a kernel of an integral operator since $k$ was the kernel of an integral operator. 
    
    We claim that $H$ is weakly compact. Assume the contrary, so that there exists $C > 0$ such that for all $n \in \mathbb{N}$ we have
    \begin{equation*}
        \esssup_{x \in A} \int_\Theta |h_n(x,y)|d\nu(y) \geq C > 0,
    \end{equation*}
    where naturally $h_n$ are the truncated kernels of $H$. For each $n \in \mathbb{N}$ choose $M_n \subseteq A$ of positive measure such that
    \begin{equation*}
        \int_\Theta |h_n(x,y)|d\nu(y) \geq C/2
    \end{equation*}
    for all $x \in M_n$. In particular, we have $n < s(x)$ for all $x \in M_n$ and therefore
    \begin{equation*}
        |k_n(x,y)| = [|k(x,y)| - s(x)]^+ + |h_n(x,y)|
    \end{equation*}
    whenever $x \in M_n$ and $y \in \Theta$. We obviously have
    \begin{align*}
        \esssup_{x \in X} \{g(x) + \int_{\Theta} |k_n(x,y)|d\nu(y)\} \geq \esssup_{x \in M_n} \{g(x) + \int_{\Theta} |k_n(x,y)|d\nu(y)\},
    \end{align*}
    while for $x \in M_n$ we get
    \begin{align*}
        g(x) + \int_{\Theta} |k_n(x,y)|d\nu(y) &= g(x) + \int_{\Theta} [|k(x,y)| - s(x)]^+ d\nu(y) + \int_\Theta |h_n(x,y)| d\nu(y)
        \\
        &= f(x, s(x)) + \int_\Theta |h_n(x,y)| d\nu(y) = \alpha_\Theta(T) + \int_\Theta |h_n(x,y)| d\nu(y) \\
        &\geq \alpha_\Theta(T) + C/2,
    \end{align*}
    and thus
    \begin{equation*}
        \esssup_{x \in X} \{g(x) + \int_{\Theta} |k_n(x,y)|d\nu(y)\} \geq \alpha_\Theta(T) + C/2.
    \end{equation*}
    Taking the infimum on the left side gives
    \begin{equation*}
        \alpha_\Theta(T) \geq \alpha_\Theta(T) + C/2,
    \end{equation*}
    a contradiction. Therefore, $H$ is a weakly compact operator, and we will show it is the best approximant. For this, observe that for $x \not \in A$ then $s(x) = 0$ and $|T|^* \mathds{1}_Y(x) \leq \alpha_{\Theta}(T)$, so that
    \begin{equation*}
        \esssup_{x \in X \backslash A} |T - H|^* \mathds{1}_Y(x) =  \esssup_{x \in X \backslash A} |T|^* \mathds{1}_Y(x)  \leq \alpha_{\Theta}(T),
    \end{equation*}
    while
    \begin{align*}
          \esssup_{x \in A} |T - H|^*\mathds{1}_Y(x) &= \esssup_{x \in A} \{ g(x) + \int_\Theta |k(x,y) - h(x,y)|d\nu(y)\} 
        \\
        &= \esssup_{x \in A} \{ g(x) + \int_\Theta [|k(x,y)| - s(x)]^+d\nu(y)\} = \esssup_{x \in A} \{f(x, s(x))\} \\
        &= \alpha_\Theta(T).
    \end{align*}
    It follows that
    \begin{equation*}
        \norm{T - H} = \norm{|T - H|^*} = \esssup_{x \in X} |T - H|^* \mathds{1}_{Y}(x) \leq \alpha_\Theta(T) \leq \norm{T}_{w, \Theta},
    \end{equation*}
    where the last inequality follows from Proposition \ref{prop: alpha-to-essential-theta}. This finishes the proof.
\end{proof}

\begin{corollary}\label{cor: alpha-decreasing-as-theta-increases}
     Let $(X, \Sigma, \mu)$ and $(Y, \Gamma, \nu)$ be strictly localizable measure spaces and $T: L_1(X, \Sigma, \mu) \to L_1(Y, \Gamma, \nu)$ be an operator. Then for any measurable subsets $\Pi \subseteq \Theta \subseteq Y$ of finite measure, we have
     \begin{equation*}
         \alpha_{\Theta}(T) \leq \alpha_\Pi(T).
     \end{equation*}
\end{corollary}
\begin{proof}
    From the definition, it is clear that $\norm{T}_{w,\Theta} \leq \norm{T}_{w, \Pi}$ whenever $\Pi \subseteq \Theta \subseteq Y$. The result now follows from Theorem \ref{th: main-with-Theta}.
\end{proof}

We are now also able to prove Theorem \ref{th: MainTheorem} when the underlying measure spaces have finite measure.

\begin{corollary}[General finite measure case]\label{cor: general-finitec-case}
    Let $(X, \Sigma, \mu)$ and $(Y, \Gamma, \nu)$ be measure spaces of finite measure. Then there exists a weakly compact operator $H$ such that
    \begin{equation*}
        \norm{T - H} = \norm{T} = \Delta(T).
    \end{equation*}
\end{corollary}
\begin{proof}
    This is a particular case of Theorem \ref{th: main-with-Theta} with $\Theta = Y$.
\end{proof}

We now prove that the range of weakly compact operators is always contained in a $\sigma$-finite subset.

\begin{proposition}\label{prop: weakly-compact-is-contained-in-sigma-finite}
    Let $(X, \Sigma, \mu)$ and $(Y, \Gamma, \nu)$ be strictly localizable measure spaces and $W: L_1(X, \Sigma, \mu) \to L_1(Y, \Gamma, \nu)$ be a weakly compact operator. Then there exists a $\sigma$-finite measurable set $\Theta \subseteq Y$, such that
    \begin{equation*}
        W[L_1(X, \Sigma, \mu)] \subseteq L_1(\Theta, \gamma, \nu),
    \end{equation*}
    where recall $\gamma = \{A \cap \Theta: A \in \Sigma\}$.
\end{proposition}
\begin{proof}
    As $(Y, \Gamma, \nu)$ is strictly localizable, we can decompose 
    \begin{equation*}
        L_1(Y, \Gamma, \nu) = \left(\bigoplus_{j \in J} L_1(Y_j, \Gamma_j, \nu)\right)_{\ell_1}
    \end{equation*}
    for some index set $J$ and disjoint sets of finite measure $Y_j \subseteq Y$.
    
    Since $W$ is weakly compact, the image of the closed unit ball $B_{L_1(X, \Sigma, \mu)}$ under $W$ is relatively weakly compact. By \cite[Lemma 7.2 (ii)]{kacena2013quantitative}, we can find a countable subset $J_0 \subseteq J$ such that $W[B_{L_1(X, \Sigma, \mu)}] \subseteq \left(\bigoplus_{j \in J_0} L_1(Y_j, \Gamma_j, \nu)\right)_{\ell_1}$. Setting $\Theta = \bigcup_{j \in J_0} Y_j$ gives the desired result.
\end{proof}

The previous result motivates the need for the following corollary, which shows that we can also find a best approximant when restricted to a $\sigma$-finite subset $\Theta \subseteq Y$.

\begin{corollary}\label{cor: sigma-finite-case}
     Let $(X, \Sigma, \mu)$ and $(Y, \Gamma, \nu)$ be strictly localizable measure spaces and $T: L_1(X, \Sigma, \mu) \to L_1(Y, \Gamma, \nu)$ be an operator. Then for any $\sigma$-finite measurable set $\Theta \subseteq Y$, there exists a weakly compact operator $H$ with range in $L_1(\Theta, \gamma, \nu)$, where $\gamma = \{A \cap \Theta: A \in \Gamma\}$, such that
     \begin{equation*}
         \norm{T - H} = \norm{T}_{w,\Theta}.
     \end{equation*}
\end{corollary}
\begin{proof}
    According to Proposition \ref{prop: sigma-finite-to-finite}, we can find an isometric lattice isomorphism $\Psi: L_1(Y, \Gamma, \nu) \to (Y, \Gamma, \tilde{\nu})$ where $\tilde{\nu}$ is a measure satisfying $\tilde{\nu}(\Theta) = 1$. Consider the operator $\tilde{T} = \Psi T: L_1(X, \Sigma, \mu) \to L_1(Y, \Gamma, \tilde{\nu})$ and note that by Theorem \ref{th: main-with-Theta} we can find a best approximant $\tilde{H}$ such that $\norm{\tilde{T} - \tilde{H}} = \norm{\tilde{T}}_{w, \Theta}$. It is now easy to check that $H = \Psi^{-1}\tilde{H}$ makes the result hold.
\end{proof}

We are finally ready to prove the first part of Theorem \ref{th: MainTheorem}.

\begin{proposition}\label{prop: first-part-MainTheorem}
    Let $(X, \Sigma, \mu)$ and $(Y, \Gamma, \nu)$ be strictly localizable measure spaces, and $T: L_1(X, \Sigma, \mu) \to L_1(Y, \Gamma, \nu)$ be an operator. Then there exists a weakly compact operator $H:  L_1(X, \Sigma, \mu) \to L_1(Y, \Gamma, \nu)$ such that
    \begin{equation*}
        \norm{T}_w = \norm{T - H}.
    \end{equation*}
\end{proposition}
\begin{proof}
    Choose a sequence of weakly compact operators $(W_n)_{n \in \mathbb{N}}$ such that $\norm{T - W_n} \to \norm{T}_w$. By Proposition \ref{prop: weakly-compact-is-contained-in-sigma-finite}, for each $n \in \mathbb{N}$, we can find a $\sigma$-finite measurable set $\Theta_n \subseteq Y$ such that $W_{n}[L_1(X, \Sigma, X)] \subseteq L_1(\Theta_n, \gamma_n, \nu)$ where $\gamma_n := \{A \cap \Theta_n: A \in \Gamma\}$. 

    Define $\Theta = \bigcup_{n \in \mathbb{N}} \Theta_n$ and $\gamma = \{A \cap \Theta: A \in \Gamma\}$. Observe that $\Theta$ is still $\sigma$-finite and  that $W_n[L_1(X, \Sigma, \mu)] \subseteq L_1(\Theta, \gamma, \nu)$ for all $n \in \mathbb{N}$. It follows that
    \begin{align*}
        \norm{T - W_n} &\geq \inf\{\norm{T - W} \hspace{5pt}| W: L_1(X, \Sigma, \mu) \to L_1(\Theta, \gamma, \nu) \text{ weakly compact} \} \\
        &= \norm{T}_{w, \Theta},
    \end{align*}
    and thus $\norm{T}_w \geq \norm{T}_{w, \Theta}$. Since the reverse inequality always holds, we conclude that $\norm{T}_w = \norm{T}_{w, \Theta}$. By Corollary \ref{cor: sigma-finite-case} we can find a weakly compact operator $H$ with range in $L_1(\Theta, \gamma, \nu)$ satisfying
    \begin{equation*}
        \norm{T - H} = \norm{T}_{w, \Theta} = \norm{T}_w. \qedhere
    \end{equation*}
\end{proof}

A duality argument now allows us to provide a proof of Corollary \ref{cor: bestapproximationC(K)}.

\begin{proof}[Proof of Corollary \ref{cor: bestapproximationC(K)}]
   Since $L$ is extremally disconnected, the space $C(L)$ is one complemented in its bidual $C(L)^{**}$. Denote by $P$ the corresponding norm one projection from $C(L)^{**}$ onto $C(L)$. Under these conditions \cite[Proposition 1.6 (i)]{tylli2001duality} gives $\norm{T}_w = \norm{T^*}_w$.

   Observe that $C(K)^*$ and $C(L)^*$ are $AL$-spaces, and therefore they can be isometrically identified with $L_1(X, \Sigma, \mu)$ and $L_1(Y, \Gamma, \nu)$, for some strictly localizable measure spaces $(X, \Sigma, \mu)$ and $(Y, \Gamma, \nu)$. Therefore, Proposition \ref{prop: first-part-MainTheorem} gives a weakly compact operator $H: C(L)^* \to C(K)^*$ such that $\norm{T^* - H} = \norm{T^*}_w$. Since $\norm{T^* - H} = \norm{T^{**} - H^*}$ it follows that
   \begin{equation*}
       \norm{T - PH^*|_{C(K)}} = \norm{PT^{**}|_{C(K)} - PH^*|_{C(K)}} \leq \norm{T}_w.
   \end{equation*}
   Taking $W = PH^*|_{C(K)}$ gives the result.
\end{proof}

The rest of this section is dedicated to proving the second part of Theorem \ref{th: MainTheorem}, namely, we want to show
\begin{equation*}
    \norm{T}_w = \inf_{\Theta\subseteq Y, \nu(\Theta)<\infty} \alpha_{\Theta}(T).
\end{equation*}

Before we can do so, we will need the following preliminary results.

\begin{lemma}\label{lmm: g-located-finite}
    Let $(Y, \Gamma, \nu)$ be a strictly localizable measurable space, $g \in L_1(Y, \Gamma, \nu)$ and $\Pi, \Theta \subseteq Y$ measurable sets of finite measure.

    Then, for every $\varepsilon > 0$, there exists a measurable subset of finite measure $\Theta_\varepsilon \subseteq Y$ such that $\Pi \cup \Theta \subseteq \Theta_\varepsilon$ and
    \begin{equation*}
        \norm{\chi_{\Theta_\varepsilon \backslash \Theta} g - \chi_{Y \backslash \Theta} g} < \varepsilon.
    \end{equation*}
\end{lemma}
\begin{proof}
    Let $((Y_j, \Gamma_j, \nu))_{j \in J}$ be a decomposition of $Y$, so that $g = \sum_{j \in J} \chi_{Y_j} g$. Only countably many summands in this sum are not equal to zero, so that we can write
    \begin{equation*}
        g = \sum_{n = 1}^\infty \chi_{Y_{j_n}}g.
    \end{equation*}
    For any $\varepsilon > 0$ we can choose $N \in \mathbb{N}$ large enough such that
    \begin{equation*}
        \norm{g - \sum_{n = 1}^N \chi_{Y_{j_n}}g} < \varepsilon.
    \end{equation*}
    Define $\Theta_\varepsilon = \Pi \cup \Theta \cup \left( \bigcup_{n = 1}^N Y_{i_n} \right)$. Then $\Theta_\varepsilon$ has finite measure and
    \begin{equation*}
        \norm{\chi_{\Theta_\varepsilon \backslash \Theta} g - \chi_{Y \backslash \Theta}g} = \norm{\chi_{Y \backslash \Theta_\varepsilon} g} \leq \norm{\chi_{Y \backslash \bigcup_{n = 1}^NY_{i_n}}g} < \varepsilon.
    \end{equation*}
\end{proof}

\begin{lemma}\label{lmm: witness-alpha-theta}
    Let $(X, \Sigma, \mu)$ and $(Y, \Gamma, \nu)$ be strictly localizable measure spaces, and $T: L_1(X, \Sigma, \mu) \to L_1(Y, \Gamma, \nu)$ be an operator. Denote
    \begin{equation*}
        \beta = \inf_{\Theta\subseteq Y, \nu(\Theta)<\infty} \alpha_{\Theta}(T).
    \end{equation*}
    Then there exist an increasing sequence of sets of finite measure $\Theta_n \subseteq Y$, a sequence $B_n \subseteq \Theta_n$ of measurable sets with $\nu(B_n) \to 0$ and a sequence of norm one elements $(f_n)_{n \in \mathbb{N}} \subseteq L_1(X, \Sigma, \mu)$ such that
    \begin{equation*}
        \lim_{n \to \infty} \norm{ \chi_{\Theta_{n+1} \backslash \Theta_n} T - \chi_{B_n}T} = \lim_{n \to \infty} \norm{ \chi_{\Theta_{n+1} \backslash \Theta_n} Tf_n - \chi_{B_n}T f_n } = \lim_{n \to \infty} \alpha_{\Theta_n}(T) = \beta.
    \end{equation*}
\end{lemma}
\begin{proof}
    Recall that, by definition, $\alpha_{\Theta}(T) = \limsup_{\nu(B) \to 0, B \subseteq \Theta} \norm{\chi_{Y \backslash \Theta} T + \chi_B T}$. Choose a sequence $(\Pi_n)_{n \in \mathbb{N}}$ of measurable sets of finite measure in $Y$ such that $\lim_{n \to \infty} \alpha_{\Pi_n}(T) = \beta$ and fix a sequence of positive scalars $(\varepsilon_n)_{n \in \mathbb{N}}$ converging to zero.

    Let $\Theta_1 = \Pi_1$. We can find a norm one function $f_1$ and $B_1 \subseteq \Theta_1$ with $\nu(B_1) < 1$ such that
    \begin{equation}\label{eq: technicaleq1}
        \norm{\chi_{Y \backslash \Theta_1} T + \chi_{B_1}T} - \varepsilon_1 < \alpha_{\Theta_1}(T) < \norm{\chi_{Y \backslash \Theta_1} Tf_1 + \chi_{B_1}Tf_1} + \varepsilon_1.
    \end{equation}
    By Lemma \ref{lmm: g-located-finite} we can find a measurable subset of finite measure $\Theta_2 \subseteq Y$ such that $\Theta_1 \cup \Pi_2 \subseteq \Theta_2$ and $\norm{\chi_{\Theta_2 \backslash \Theta_1} Tf_1 - \chi_{Y \backslash \Theta}Tf_1}$ is small enough to make \eqref{eq: technicaleq1} hold when we replace $\chi_{Y \backslash \Theta_1}Tf_1$ by $\chi_{\Theta_2 \backslash \Theta_1} Tf_1$, in other words we get
    \begin{align*}
        \norm{\chi_{\Theta_2 \backslash \Theta_1} T + \chi_{B_1}T} - \varepsilon_1&\leq \norm{\chi_{Y \backslash \Theta_1} T + \chi_{B_1}T} - \varepsilon_1 < \alpha_{\Theta_1}(T) \\
        &< \norm{\chi_{\Theta_2 \backslash \Theta_1} Tf_1 + \chi_{B_1}Tf_1} + \varepsilon_1.
    \end{align*}
    For $\Theta_2$, we can find a norm one function $f_2$ and $B_2 \subseteq \Theta_2$ with $\nu(B_2) < 1/2$ and such that
    \begin{equation*}
        \norm{\chi_{Y \backslash \Theta_2} T + \chi_{B_2}T} - \varepsilon_2 < \alpha_{\Theta_2}(T) < \norm{\chi_{Y \backslash \Theta_2} Tf_2 + \chi_{B_2}Tf_2} + \varepsilon_2.
    \end{equation*}
    In the same manner, we can choose a measurable set of finite measure $\Theta_3 \subseteq Y$ such that $\Theta_2 \cup \Pi_3 \subseteq \Theta_3$ and $\norm{\chi_{\Theta_3 \backslash \Theta_2} Tf_2 - \chi_{Y \backslash \Theta_2}Tf_1}$ is small enough so that again
    \begin{align*}
        \norm{\chi_{\Theta_3 \backslash \Theta_2} T + \chi_{B_2}T} - \varepsilon_2 &\leq \norm{\chi_{Y \backslash \Theta_2} T + \chi_{B_2}T} - \varepsilon_2 < \alpha_{\Theta_2}(T) \\
        &< \norm{\chi_{\Theta_3 \backslash \Theta_2} Tf_2 + \chi_{B_2}Tf_2} + \varepsilon_2.
    \end{align*}

    By induction, we obtain a sequence of measurable sets of finite measure $\Theta_n \subseteq \Theta_{n+1}$ with $\Pi_n \subseteq \Theta_n$, a sequence $B_n$ with $\nu(B_n) < 1/n$ and norm one functions $f_n$ such that
    \begin{align*}
        \norm{\chi_{\Theta_{n+1} \backslash \Theta_n} T + \chi_{B_n}T} - \varepsilon_n &\leq \norm{\chi_{Y \backslash \Theta_n} T + \chi_{B_n}T} - \varepsilon_n < \alpha_{\Theta_n}(T) \\
        &< \norm{\chi_{\Theta_{n+1} \backslash \Theta_n} Tf_n + \chi_{B_n}Tf_n} + \varepsilon_n.
    \end{align*}
    It follows that
    \begin{align*}
        \beta &\leq \alpha_{\Theta_n}(T) \leq \norm{\chi_{\Theta_{n+1} \backslash \Theta_n} Tf_n + \chi_{B_n}Tf_n} + \varepsilon_n \leq \norm{\chi_{\Theta_{n+1} \backslash \Theta_n} T + \chi_{B_n}T} + \varepsilon_n \\
        & \leq \norm{\chi_{Y \backslash \Theta_n} T + \chi_{B_n}T} + \varepsilon_n \leq \alpha_{\Theta_n}(T) + 2\varepsilon_n \leq \alpha_{\Pi_n}(T) + 2\varepsilon_n,
    \end{align*}
    where the last inequality follows from $\Pi_n \subseteq \Theta_n$ and Corollary \ref{cor: alpha-decreasing-as-theta-increases}. Since $\alpha_{\Pi_n}(T) \to \beta$ and $\varepsilon_n \to 0$, taking limits gives
    \begin{equation*}
        \lim_{n \to \infty} \norm{ \chi_{\Theta_{n+1} \backslash \Theta_n} T - \chi_{B_n}T} = \lim_{n \to \infty} \norm{ \chi_{\Theta_{n+1} \backslash \Theta_n} Tf_n - \chi_{B_n}T f_n } = \lim_{n \to \infty} \alpha_{\Theta_n}(T) = \beta,
    \end{equation*}
    which completes the proof.
\end{proof}

We are finally ready for the second part of our main result.

\begin{proposition}\label{prop: second-part-MainTheorem}
     Let $(X, \Sigma, \mu)$ and $(Y, \Gamma, \nu)$ be strictly localizable measure spaces, and $T: L_1(X, \Sigma, \mu) \to L_1(Y, \Gamma, \nu)$ be an operator. Then
     \begin{equation*}
        \norm{T}_w = \inf_{\Theta\subseteq Y, \nu(\Theta)<\infty} \alpha_{\Theta}(T).
     \end{equation*}
\end{proposition}
\begin{proof}
    By Lemma \ref{lmm: witness-alpha-theta}, we can find an increasing sequence $(\Theta_n)_{n = 1}^\infty$ of measurable sets of finite measure, a sequence $(B_n)_{n=1}^\infty$ of measurable sets with $B_n \subseteq \Theta_n$ and $\nu(B_n) \to 0$, and a sequence of norm one functions $(f_n)_{n=1}^\infty \subseteq L_1(X, \Sigma, \mu)$ such that
    \begin{equation*}
         \lim_{n \to \infty} \norm{ \chi_{\Theta_{n+1} \backslash \Theta_n} T - \chi_{B_n}T} = \lim_{n \to \infty} \norm{ \chi_{\Theta_{n+1} \backslash \Theta_n} Tf_n - \chi_{B_n}T f_n } = \lim_{n \to \infty} \alpha_{\Theta_n}(T) = \beta,
    \end{equation*}
    where $\beta = \inf_{\Theta \subseteq Y, \nu(\Theta) < \infty} \alpha_{\Theta(T)}$. Define the sets $A_n = \Theta_{n+1} \backslash \Theta_{n}$.

    By Proposition \ref{prop: separable-in-separable-sublattice} we can find a $\sigma$-finite measure space $(\Omega, \tilde{\Sigma}, \mu)$ with $\tilde{\Sigma}$ countably generated such that $\spn \{ f_n: n \in \mathbb{N}\} \subseteq L_1(\Omega, \tilde{\Sigma}, \mu) =: E$. Let $J_E: E \to L_1(X, \Sigma, \mu)$ be the inclusion map. Applying the same proposition, we can find a $\sigma$-finite measure space $(\Pi, \tilde{\Gamma}, \nu)$ with $\tilde{\Gamma}$ countably generated such that $F = L_1(\Pi, \tilde{\Gamma}, \nu)$ contains $T[E]$ and the characteristic functions $\mathds{1}_{A_n}$ and $\mathds{1}_{B_n}$ for all $n \in \mathbb{N}$. Denote by $P_F: L_1(Y, \Gamma, \nu) \to F$ the norm one projection.

    Arguing as in Proposition \ref{prop: sigma-finite-to-finite} we can find measures $\tilde{\mu}$ and $\tilde{\nu}$ such that $\tilde{\mu}(\Omega) = 1$ and $\tilde{\nu}(\Pi) = 1$, as well as isometric lattice isomorphisms $\Phi: \tilde{E} \to E$ and $\Psi: F \to \tilde{F}$, where $\tilde{E} = L_1(\Omega, \tilde{\Sigma}, \tilde{\mu})$ and $\tilde{F} = L_1(\Pi, \tilde{\Gamma}, \tilde{\nu})$. Observe that since $(A_n)_{n =1}^\infty$ are disjoint and $\tilde{\nu}$ is a finite measure, then $\tilde{\nu}(A_n) \to 0$, while since $\nu(B_n) \to 0$, by the construction of $\tilde{\nu}$ (which is measure decreasing), we also have $\tilde{\nu}(B_n) \to 0$. Let  $C_n = A_n \cup B_n$, so that $\tilde{\nu}(C_n) \to 0$ and define the operator $\tilde{T} = \Psi P_F T J_E \Phi: \tilde{E} \to \tilde{F}$. Corollary \ref{cor: general-finitec-case} gives $\norm{\tilde{T}}_w = \Delta(\tilde{T})$. Observe that by Proposition \ref{prop: sigma-finite-to-finite} we know that $\Psi$ commutes with $\chi_{C_n}$, while $\chi_{C_n} P_F = \chi_{C_n}$ since $C_n \in \tilde{\Gamma}$. From this, we get
      \begin{align*}
        \chi_{C_n}\tilde{T}=\chi_{C_n}\Psi P_F T J_{E}\Phi=\Psi (\chi_{C_n} P_F T J_{E}\Phi)
        =\Psi (\chi_{C_n} T J_{E}\Phi),
    \end{align*}
    and therefore with the fact that $\Psi$ and $\Phi$ are isometries
    \begin{align*}
        \beta &= \lim_{n \to \infty} \norm{\chi_{C_n}T} = \lim_{n \to \infty} \norm{\chi_{C_n}\tilde{T}} \leq \Delta(\Tilde{T}) = \norm{\tilde{T}}_w \leq \norm{\Psi}_w \norm{P_F}_w \norm{T}_w \norm{J_E}_w \norm{\Phi}_w \\
        &\leq  \norm{T}_w.
    \end{align*}
    By Theorem \ref{th: main-with-Theta}, for all $\Theta \subseteq Y$ of finite measure, we have 
    \begin{equation*}
        \alpha_{\Theta}(T) = \norm{T}_{w, \Theta} \geq \norm{T}_w.
    \end{equation*}
    Taking the infimum gives
    \begin{equation*}
        \beta = \inf_{\Theta \subseteq Y, \nu(\Theta) < \infty} \alpha_{\Theta}(T) \geq \norm{T}_w,
    \end{equation*} 
    which finishes the proof.
\end{proof}

\bigskip
\section{The Weak Calkin Algebra in $AL$-spaces}\label{sec: weak-calkinAL}

We now turn to the proof of Theorem \ref{th: AL-factorization}. We first require several preliminary results.

\begin{lemma}\label{lmm: almost-disjoint-L1-1}
    Let $(X, \Sigma, \mu)$ and $(Y, \Gamma, \nu)$ be strictly localizable measurable spaces, and $T: L_1(X, \Sigma, \mu) \to L_1(Y, \Gamma, \nu)$ be a non-weakly compact operator. Then there exist a sequence $(f_n)_{n=1}^\infty \subseteq L_1(X, \Sigma, \mu)$ of norm one functions and pairwise disjoint sets $(A_n)_{n=1}^\infty \subseteq \Gamma$ of finite measure such that
    \begin{equation*}
        \lim_{n \to \infty} \norm{\chi_{A_n} Tf_n} \geq \norm{T}_w
    \end{equation*}
\end{lemma}
\begin{proof}
     We construct $(f_n)_{n=1}^\infty$ inductively. Since $\norm{T} \geq \norm{T}_w$, there exists a norm one function $f_1 \in L_1(X, \Sigma, \mu)$ such that $\norm{Tf_1} > \norm{T}_w - 1/4$. In particular, we can find $B_{1} \in \Gamma$ of finite measure so that $\norm{\chi_{B_{1}}Tf_1} > \norm{T}_w - 1/4$. Set $C_1 = \emptyset$.

    Inductively, suppose we have constructed norm one functions $f_1, \dots, f_n$ and sets of finite measure $B_{1}, \dots, B_n \in \Gamma$, $C_{1}, \dots, C_n \in \Gamma$ such that:
    \begin{enumerate}[label=(\alph*), ref=(\alph*)]
            \item \label{it:1} $\norm{\chi_{B_{i}}Tf_i} > \norm{T}_w - 1/2^{i+1}$ for every $1 \leq i \leq n$.
            \item \label{it:2} $\norm{\chi_{C_{j}} Tf_i} < 1/2^{j + 1}$ for every $1 \leq i < j \leq n$.
            \item \label{it:3} $(B_{i} \backslash C_{j}) \cap B_{j} = \emptyset$ for every $1 \leq i < j \leq n$.
    \end{enumerate}
    We build $f_{n+1}$, $B_{n+1}$ and $C_{n+1}$. Let $\Theta = \bigcup_{i=1}^n B_i \in \Gamma$, which has finite measure, and find $\delta > 0$ small enough so that
    \begin{equation*}
        \norm{\chi_C Tf_i} < 1/2^{n+2}
    \end{equation*}
    for all $i=1, \dots, n$ and all sets $C \in \Gamma$ with $\nu(C) < \delta$.
    
    By Theorem \ref{th: MainTheorem}, $\norm{T}_w \leq \alpha_\Theta(T)$ and by definition of $\alpha_{\Theta}(T)$, we can find a measurable set $C \in \Gamma$ with $\nu(C) < \delta$ such that
    \begin{equation*}
        \norm{\chi_{Y \backslash \Theta}T + \chi_C T} > \alpha_{\Theta}(T) - 1/2^{n+2} \geq \norm{T}_w - 1/2^{(n+1)+1}.
    \end{equation*}
    Therefore, we can find a norm one function $f_{n+1}$ such that
    \begin{equation*}
         \norm{\chi_{Y \backslash \Theta}T f_{n+1} + \chi_C T f_{n+1}} > \norm{T}_w - 1/2^{(n+1)+1}.
    \end{equation*}
    By Lemma \ref{lmm: g-located-finite} we can find a measurable set of finite measure $D \subseteq Y$ such that $\norm{\chi_{D \backslash \Theta}T f_{n+1} - \chi_{Y \backslash \Theta}T f_{n+1}}$ is small enough to ensure that
    \begin{equation*}
         \norm{\chi_{D \backslash \Theta}T f_{n+1} + \chi_C T f_{n+1}} > \norm{T}_w - 1/2^{(n+1)+1}
    \end{equation*}
    holds.
    Define $B_{n+1} = C \cup (D \backslash \Theta)$ and $C_{n+1} = C$, it is clear that, with this definition, \ref{it:1}, \ref{it:2} and \ref{it:3} are satisfied. This finishes the inductive construction.

    Finally, let $A_n = B_n \backslash \bigcup_{m > n} C_m$ for each $n \in \mathbb{N}$. By construction, $A_n \cap A_m = \emptyset$ whenever $n \not = m$ and
    \begin{align*}
        \norm{\chi_{A_n} Tf_n} &\geq \norm{\chi_{B_n} Tf_n} - \sum_{m > n} \norm{\chi_{C_m} Tf_n} \\
        &\geq \norm{T}_w - 1/2^{n+1} - \sum_{m > n} 1/2^{m+1} = \norm{T}_w - 1/2^{n}.
    \end{align*}
    This shows that $\limsup_{n \to \infty} \norm{\chi_{A_n} Tf_n } \geq \|T\|_w \quad \text{as } n \to \infty$. Passing to a subsequence, we can assume that in fact $\lim_{n \to \infty} \norm{\chi_{A_n} Tf_n } \geq \|T\|_w$
\end{proof}

We will also need the following well-known result, see, e.g., \cite[Lemma 2.3]{wolff1984essential}.

\begin{proposition}\label{prop: to-zero-convex-sums-signed}
    Let $E$ be a Banach space and $(x_n)_{n = 1}^\infty \subseteq E$ be a weakly convergent sequence. Then there exist an increasing sequence of natural numbers $(i_n)_{n=1}^\infty$ and a sequence of scalars $(\alpha_i)_{i=1}^\infty$ satisfying
    \begin{enumerate}[label = (\roman*), ref = (\roman*)]
        \item $\sum_{i = i_n + 1}^{i_{n+1}} |\alpha_i| = 1$,
        \item the sequence $(x_n')_{n=1}^\infty$ defined by $x_n' = \sum_{i = i_{n} + 1}^{i_{n+1}} \alpha_i x_i$ converges in norm to $0$.
    \end{enumerate}
\end{proposition}

We next recall the classical splitting lemma of Pe\l czy\'nski and Kadec \cite{kadec1962bases}; see \cite[Lemma 5.2.8]{albiac2006topics} or \cite[Lemma 2.2]{wolff1984essential} for an explicit formulation.

\begin{lemma}[Subsequence Splitting Lemma]\label{lmm: splitting} 
    Let $(X, \Sigma, \mu)$ be a measure space of finite measure and $(f_n)_{n=1}^\infty$ be a bounded sequence in $L_1(X, \Sigma, \mu)$.
    Then, there exist a subsequence $(f_{n_k})_{k=1}^\infty$ and two other sequences $(g_k)_{k=1}^\infty$, $(h_k)_{k=1}^\infty$ such that the following three assertions hold:
    \begin{itemize}
        \item[(i)] $f_{n_k} = g_k + h_k$ for all $k \in \mathbb{N}$,
        \item[(ii)] $\inf \left( |g_k|, |h_k| \right) = 0 = \inf \left( |h_\ell|, |h_m| \right) \quad \text{for all } k,\ \text{and all } \ell \neq m$,
        \item[(iii)] $(g_k)_{k = 1}^\infty$ is weakly convergent.
    \end{itemize}
\end{lemma} 

\begin{definition}
    A sequence $(f_n)_{n=1}^\infty$ in $L_1(X, \Sigma, \mu)$ is called \emph{disjoint} if $|f_n| \wedge |f_m| = 0$ whenever $m \not = n$. A sequence $(g_n)_{n=1}^\infty$ in $L_1(X, \Sigma, \mu)$ is called \emph{almost disjoint} if there exists a disjoint sequence $(f_n)_{n=1}^\infty$ satisfying $\norm{f_n - g_n} \to 0$ as $n \to \infty$.
\end{definition}

\begin{lemma}\label{lmm: lemma-Weis1}
    Let $(X, \Sigma, \mu)$ and $(Y, \Gamma, \nu)$ be strictly localizable measure spaces and $T: L_1(X, \Sigma, \mu) \to L_1(Y, \Gamma, \nu)$ be a non-weakly compact operator. Then there exists a disjoint sequence $(u_n)_{n=1}^\infty \subseteq L_1(X, \Sigma, \mu)$ with $\norm{u_n} = 1$ such that $(Tu_n)_{n=1}^\infty$ is almost disjoint and $\norm{Tu_n} \to \norm{T}_w$.
\end{lemma}
\begin{proof}
    Choose a sequence of norm one functions $(f_n)_{n=1}^\infty \subseteq L_1(X, \Sigma, \mu)$ and disjoint measurable sets $(A_n)_{n=1}^\infty$ according to Lemma \ref{lmm: almost-disjoint-L1-1}.

    Arguing as in the proof of Proposition \ref{prop: second-part-MainTheorem}, we can find $\sigma$-finite measure spaces $(\Omega, \tilde{\Sigma}, \mu)$ and $(\Pi, \tilde{\Gamma}, \nu)$ with countably generated $\sigma$-algebras $\tilde{\Sigma} \subseteq \Sigma $ and $\tilde{\Gamma} \subseteq \Gamma $, such that
    \begin{equation*}
        \spn\{f_n : n \in \mathbb{N}\} \subseteq L_1(\Omega, \tilde{\Sigma}, \mu) =: E
    \end{equation*}
    and
    \begin{equation*}
        T[E] + \spn\{\mathds{1}_{A_n} : n \in \mathbb{N}\} \subseteq L_1(\Theta, \tilde{\Gamma}, \nu) =: F.
    \end{equation*}
    Note that, the previous conditions in particular imply that $(A_n)_{n=1}^\infty \subseteq \tilde{\Gamma}$.
    
    Moreover, again by the same arguments as in the proof of Proposition \ref{prop: second-part-MainTheorem}, we can find measures $\tilde{\mu}$ and $\tilde{\nu}$ such that $\tilde{\mu}(\Omega) = 1$ and $\tilde{\nu}(\Pi) = 1$, as well as isometric lattice isomorphisms $\Phi: \tilde{E} \to E$ and $\Psi: F \to \tilde{F}$, where $\tilde{E} = L_1(\Omega, \tilde{\Sigma}, \tilde{\mu})$ and $\tilde{F} = L_1(\Pi, \tilde{\Gamma}, \tilde{\nu})$. Furthermore, $\Psi$ and $\Phi$ commute with multiplication operators, in the sense of Proposition \ref{prop: sigma-finite-to-finite}.
    
    Observe that since $(A_n)_{n =1}^\infty$ are disjoint and $\tilde{\nu}$ is a finite measure, then $\tilde{\nu}(A_n) \to 0$. Let $J_E: E \to L_1(X, \Sigma, \mu)$ denote the inclusion operator, and let $P_F: L_1(Y, \Gamma, \nu) \to F$ be the norm-one projection onto $F$. Consider the operator $\tilde{T} = \Psi P_F T J_E \Phi: \tilde{E} \to \tilde{F}$ and define the sequence of norm one functions $(\tilde{f}_n)_{n=1}^\infty = (\Phi^{-1}(f_n))_{n=1}^\infty$. 
    
    Clearly, we have $\norm{\tilde{T}}_w \leq \norm{T}_w$, while by Corollary \ref{cor: general-finitec-case} we obtain
    \begin{equation*}
       \norm{\tilde{T}}_w = \Delta(\tilde{T}) \geq \lim_{n \to \infty} \norm{\chi_{A_n} \tilde{T} \tilde{f}_n} = \lim_{n \to \infty} \norm{\chi_{A_n} Tf_n} \geq \norm{T}_w,
    \end{equation*}
    and thus $\Delta(\tilde{T}) = \norm{\tilde{T}}_w = \norm{T}_w = \lim_{n \to \infty} \norm{\chi_{A_n} \tilde{T} \tilde{f}_n}$. 
    
    The rest of the proof now follows from an identical argument to that of \cite[Lemma 2.4]{wolff1984essential}; we include it here for the reader's convenience.

    Passing to a subsequence, we can assume $\tilde{f}_n = \tilde{g}_n + \tilde{h}_n$ according to Lemma \ref{lmm: splitting}. Since $\tilde{T}$ is weakly continuous, $(\tilde{T}\tilde{g}_n)_{n=1}^\infty$ is equi-integrable and thus \break $\lim_{n \to \infty} \norm{\chi_{A_n} \tilde{T} \tilde{g}_{n}} = 0$, so
    \begin{equation*}
        \lim_{n \to \infty} \norm{\chi_{A_n} \tilde{T} \tilde{h}_{n}} = \lim_{n \to \infty} \norm{\chi_{A_n} \tilde{T} \tilde{f}_{n}} = \Delta(\tilde{T}).
    \end{equation*}
    We claim that $\limsup_{n \to \infty} \norm{\tilde{h}_n} = 1$. Indeed, suppose that $\limsup_{n \to \infty} \norm{\tilde{h}_n} = c < 1$ taking $\tilde{h}'_n = \tilde{h}_n/\norm{\tilde{h}_n}$ we would obtain $\Delta(\tilde{T}) \geq \Delta(\tilde{T})/c$, a contradiction. 

    Thus, without loss of generality, we may assume $\norm{\tilde{h}_n} = 1$. By a similar reasoning, using the definition of $\Delta(\tilde{T})$, the set $\{\chi_{\Theta \backslash A_n} \tilde{T} \tilde{h}_n: n \in \mathbb{N}\}$ has to be relatively weakly compact. Passing to a further subsequence, we can assume that $(\chi_{\Theta \backslash A_n} \tilde{T} \tilde{h}_n)_{n=1}^\infty$ is weakly convergent. By Proposition \ref{prop: to-zero-convex-sums-signed} we obtain an increasing sequence $(i_n)_{n=1}^\infty$ and a sequence of scalars $(\alpha_i)_{i=1}^\infty$ such that for
    \begin{equation*}
        \tilde{v}_n := \sum_{i=i_{n}+1}^{i_{n+1}} \alpha_i \chi_{\Theta \backslash A_i} \tilde{T} \tilde{h}_i,
    \end{equation*}
   we get $\lim_{n \to \infty} \norm{\tilde{v}_n} = 0$. Define $\tilde{u}_n = \sum_{i=i_{n}+1}^{i_{n+1}} \alpha_i \tilde{h}_i$. By construction, $(\tilde{u}_n)_{n=1}^\infty$ is a sequence of norm one mutually disjoint functions. Consider 
   \begin{equation*}
       \tilde{w}_n = \tilde{T}\tilde{u}_n - \tilde{v}_n = \sum_{i= i_n + 1}^{i_{n+1}} \alpha_i \chi_{A_i} \tilde{T} \tilde{h}_i,
   \end{equation*}
   so that $\tilde{w}_n$ is a disjoint sequence. Recall that $\lim_{n \to \infty} \norm{ \chi_{A_n} \tilde{T} \tilde{h}_n} = \Delta(\tilde{T})$, which clearly implies $\lim_{n \to \infty} \norm{\tilde{w}_n} = \Delta(\tilde{T})$. Hence $(\tilde{T}\tilde{u}_n)_{n=1}^\infty = (\tilde{w}_n + \tilde{v}_n)_{n=1}^\infty$ is almost disjoint and satisfies $\lim_{n \to \infty} \norm{\tilde{T} \tilde{u}_n} = \Delta(\Tilde{T}) = \norm{T}_w$. Observe that, since both $\Phi$ and $\Psi$ are lattice isomorphisms, they preserve disjointness of supports. In particular, it is easy to see that the sequence
    \begin{equation*}
        (u_n)_{n=1}^\infty = \bigl(\Phi^{-1}(\tilde{u}_n)\bigr)_{n=1}^\infty \subseteq E,
    \end{equation*}
    viewed as a sequence in the ambient space $L_1(X,\Sigma,\mu)$, satisfy the conditions of the lemma. This completes the proof.
\end{proof}

Using the previous lemma, we can give a proof of the factorization theorem in $AL$-spaces.

\begin{proof}[Proof of Theorem \ref{th: AL-factorization}]
    Let $(X, \Sigma, \mu)$ and $(Y, \Gamma, \nu)$ be strictly localizable measure spaces such that $E = L_1(X, \Sigma, \mu)$ and $F = L_1(Y, \Gamma, \nu)$. The case when $T: E \to F$ is weakly compact is obvious; thus, by normalizing, we may assume $\norm{T}_w = 1$. By Lemma \ref{lmm: lemma-Weis1} we can find a normalized disjoint sequence $(u_n)_{n=1}^\infty \subseteq E$ such that $(Tu_n)_{n=1}^\infty$ is almost disjoint and $\lim_{n \to \infty} \norm{Tu_n} = 1$. Therefore, for any $\varepsilon > 0$ we can, after possibly passing to a subsequence, find a sequence of normalized disjoint functions $(g_n)_{n=1}^\infty$ such that
    \begin{equation*}
        2\sum_{n=1}^\infty \norm{Tu_n - g_n } < \varepsilon.
    \end{equation*}
    Observe that $Z = \overline{\spn \{g_n: n \in \mathbb{N}}\}$ is isometric to $\ell_1$ and $1$-complemented. Let $P: F \to Z$ be a norm-one projection and $B: Z \to \ell_1$ be the isometric isomorphism given by $g_n \mapsto e_n$, where $e_n$ is the canonical basis of $\ell_1$. By the Principle of small perturbations, see for example \cite[Theorem 1.3.9]{albiac2006topics}, we can find an isomorphism $A: F \to F$ such that $A(g_n) = Tu_n$ and $\norm{A^{-1}} \leq (1 - \varepsilon)^{-1}$. Define $V: \ell_1 \to E$ by $e_n \mapsto u_n$, and since $u_n$ are disjoint and normalized, observe that $V$ is an isometric embedding, in particular $\norm{V} = 1$. By construction $(BPA^{-1})TV = UTV = I_{\ell_1}$, where $U := BPA^{-1}$. Observe that $\norm{U} \leq \norm{B} \norm{P} \norm{A^{-1}} \leq (1 - \varepsilon)^{-1}$, in other words $\norm{U}^{-1} \norm{V}^{-1} \geq  1 - \varepsilon$. This shows that
    \begin{equation*}
        \norm{T}_w \leq \sup \{\norm{U}^{-1}\norm{V}^{-1}: UTV = I_{\ell_1} \},
    \end{equation*}  
    the reverse inequality is immediate, and thus the result follows.
\end{proof}

We now provide an example which shows that no such quantitative factorizations are possible unless additional constraints are imposed on the domain.

\begin{example}\label{ex: Not-main-L}
   By \cite[Corollary 3]{astala1990seminorms}, $L_1[0,1]$ does not have the weak approximation property, and thus by \cite[Corollary 3]{astala1990seminorms}, there exists a Banach space $F$ such that $\norm{\cdot}_w$ and $\omega(\cdot)$ are not equivalent, where recall that $\omega(\cdot)$ is the De Blasi measure of weak noncompactness \eqref{eq: deBlasi}. In particular, we can find a sequence of operators $T_n: F \to L_1[0,1]$ with $\norm{T_n}_w = 1$ and 
   \begin{equation*}
       \sup \{ \norm{U}^{-1}\norm{V}^{-1}: UT_nV = I_{\ell_1}\} \leq \omega(T_n) < 1/n,
   \end{equation*}
   where the first inequality follows since for any $U, V$ such that $UT_nV = I_{\ell_1}$ we have
   \begin{equation*}
    1 = \omega(I_{\ell_1}) = \omega(UT_nV) \leq \norm{U}\norm{V} \omega(T_n),
   \end{equation*}
   that is $\omega(T_n) \geq \norm{U}^{-1} \norm{V}^{-1}$. This shows that no such quantitative factorizations are possible using the weak essential norm.
\end{example}
\begin{rem}
    Incidentally, the same reasoning yields an identical conclusion for any $AL$-space without the Schur property, showing that no general factorization is possible in that case either. On the other hand, if an $AL$-space $E$ has the Schur property, then $E \cong \ell_1(\Gamma)$ for some set $\Gamma$. In this case, such a factorization is easily achieved, independently of the domain of the operator.
\end{rem}

Lastly, we argue that, from this factorization, the weak Calkin algebra admits a unique algebra norm.

\begin{proof}[Proof of Corollary \ref{cor: unique-algebra-normAL-space}]
    Let $E = L_1(X, \Sigma, \mu)$ for some strictly localizable measure space. It is well-known that (see \cite[332B]{fremlin2001v3} for the real case and \cite[Corollary, p.~136]{lacey2008isometric} for the complex case)
    \begin{equation*}
        E \cong \left(\ell_1(\Gamma)\oplus \left(\bigoplus_{m_\alpha \in A}L_1[0, 1]^{m_{\alpha}} \right)_{\ell_1}\right)_{\ell_1},
    \end{equation*}
    for some set $\Gamma$ and some indexed family of cardinal numbers $A$. In particular, it follows that $E \cong E \oplus E$, provided $E$ is infinite-dimensional; the finite-dimensional case is trivial. A classical theorem of B. E.~Johnson~\cite{johnson1967continuity} then ensures that every homomorphism from $\mathscr{B}(E)$, equipped with the operator norm, into a Banach algebra is continuous. Consequently, by \cite[Theorem 6.1.5]{palmer1994banach}, the weak essential norm is maximal on the weak Calkin algebra $\mathscr{B}(E)/\mathscr{W}(E)$. Theorem~\ref{th: AL-factorization}, together with~\cite[Proposition 2.12]{arnott2023uniqueness}, shows that the weak essential norm is incompressible (in fact, uniformly so), and the reasoning in~\cite[Corollary 2.10]{arnott2023uniqueness} then yields that $\mathscr{B}(E)/\mathscr{W}(E)$ admits a unique algebra norm.
\end{proof}

\subsection{Quantitative weak noncompactness of operators on $AL$-spaces.}

In this section, we will show that, in a broad sense, measures of weak noncompactness for operators on $AL$-spaces coincide. To this end, we first need to clarify what we mean by a measure of weak noncompactness. Throughout this section, we fix an $AL$-space $E = L_1(X, \Sigma, \mu)$. 

\begin{definition}\label{def: measureweakcompactness}
    A function $\alpha: \mathscr{B}(E) \to [0, \infty)$ is said to be a \emph{measure of weak noncompactness} if it satisfies:
    \begin{enumerate}[label=(\alph*), ref=(\alph*)]
        \item \label{it: en1} $\alpha(T + S) \leq \alpha(T) + \alpha(S)$ and $\alpha(\lambda T) = |\lambda| \alpha(T)$ for all $T, S \in \mathscr{B}(E)$ and $\lambda \in \mathbb{K}$.
        \item \label{it: en2} $\alpha(TS) \leq \min\{\norm{T}\alpha(S), \alpha(T)\norm{S} \}$ for all $T, S \in \mathscr{B}(E)$.
        \item \label{it: en3} $\alpha(T) = 0$ if and only if $T$ is weakly compact.
        \item \label{it: en4} There exist $J_{\ell_1}: \ell_1 \to E$ an isometric inclusion and $\pi_{\ell_1}: E \to \ell_1$ a norm one operator, satisfying $\pi_{\ell_1} J_{\ell_1} = I_{\ell_1}$ and $\alpha(J_{\ell_1} \pi_{\ell_1}) \geq 1$.
    \end{enumerate}
    We will say that $\alpha$ is normalized if, additionally, we have $\alpha(I_E) = 1$.
\end{definition}

\begin{rem}
    In \ref{it: en4}, we emphasise that it suffices to exhibit any such $J_{\ell_1}$ and $\pi_{\ell_1}$. In terms of set measures of weak noncompactness, this condition says that the unit ball of $\ell_1$ has measure at least $1$.
\end{rem}

We will need the following elementary observation.

\begin{lemma}\label{lmm: majorized-weak}
    Let $\alpha$ be a normalized measure of weak noncompactness and $T \in \mathscr{B}(E)$. Then
    \begin{equation*}
        \alpha(T) \leq \norm{T}_w.
    \end{equation*}
\end{lemma}
\begin{proof}
    For any weakly compact operator $W$, we have
    \begin{equation*}
        \alpha(T) = \alpha(T - W + W) \leq \alpha(T - W) + \alpha(W) = \alpha(T - W) \leq \alpha(T) + \alpha(-W) = \alpha(T),
    \end{equation*}
    and thus $\alpha(T) = \alpha(T - W)$. It follows that
    \begin{equation*}
        \alpha(T) = \alpha((T - W)I_E) \leq \norm{T - W} \alpha(I_E) = \norm{T - W}.
    \end{equation*}
    Taking the infimum over all weakly compact operators gives the result.
\end{proof}

\begin{proposition}\label{prop: normalised-measures-weak-non}
    Let $\alpha$ be a normalized measure of weak noncompactness and $T \in \mathscr{B}(E)$. Then $\alpha(T) = \norm{T}_w$.
\end{proposition}
\begin{proof}
    Let $T \in \mathscr{B}(E)$, normalising we may assume $\norm{T}_w = 1$. By Lemma \ref{lmm: majorized-weak} it is then enough to show $\alpha(T) \geq 1$. Fix $\varepsilon > 0$, by Theorem \ref{th: AL-factorization} we can find $U,V \in \mathscr{B}(E)$ such that $\norm{U} \norm{V} < 1 + \varepsilon$ and $UTV = J_{\ell_1} \pi_{\ell_1}$. It follows that
    \begin{equation*}
        1 \leq \alpha(J_{\ell_1} \pi_{\ell_1}) = \alpha(UTV) \leq \norm{U} \norm{V} \alpha(T) < (1 + \varepsilon)\alpha(T),
    \end{equation*}
    and since $\varepsilon$ was arbitrary, we get $\alpha(T) \geq 1$, which finishes the proof.
\end{proof}

This result has Corollary \ref{cor: BonsalMinimalityAL} as an easy consequence. In the next proof, for a norm $\vertiii{\cdot}$ on the weak Calkin algebra, we refer indiscriminately to $\vertiii{T}$ of any $T \in \mathscr{B}(E)$. This is to be understood as the seminorm this norm induces on $\mathscr{B}(E)$ or, alternatively, as the norm of the equivalence class represented by $T$.

\begin{proof}[Proof of Corollary \ref{cor: BonsalMinimalityAL}]
    Let $\vertiii{\cdot}$ be any algebra norm on the weak Calkin algebra $\mathscr{B}(E) / \mathscr{W}(E)$, with the property that $\vertiii{T} \leq \norm{T}_w$ for any operator $T \in \mathscr{B}(E)$. In particular, since $\vertiii{\cdot}$ is an algebra norm and $J_{\ell_1} \pi_{\ell_1}$ is a non-weakly compact idempotent, then
    \begin{equation*}
        \vertiii{J_{\ell_1} \pi_{\ell_1}} = \vertiii{(J_{\ell_1} \pi_{\ell_1})^2} \leq \vertiii{J_{\ell_1} \pi_{\ell_1}}^2,
    \end{equation*}
    that is $1 \leq \vertiii{J_{\ell_1} \pi_{\ell_1}}$. Clearly $\vertiii{I_E} = 1$ and thus $\vertiii{\cdot}$ is a normalized measure of weak noncompactness. The result now follows from Proposition \ref{prop: normalised-measures-weak-non}
\end{proof}

\begin{rem}
    Corollary \ref{cor: BonsalMinimalityAL} could also have been obtained by mimicking the arguments in the proof of \cite[Theorem 1.5]{acuaviva2024factorizations}. However, since the general form of the result for normalized measures of weak noncompactness will be needed later, and because other notions (not necessarily algebra norms) have been extensively studied in the literature, we have chosen to prove it starting from the more general formulation.
\end{rem}

Before moving on, we would like to point out that our perspective focuses on measures of weak noncompactness for operators. However, there exists a rich and extensive literature on measures of weak noncompactness for subsets $A$ of a Banach space $E$; see, for example, \cite{angosto2009measures} and \cite{kacena2013quantitative}. \\

Given a measure of weak noncompactness on sets, say $\alpha(\cdot)$, one defines the corresponding measure on an operator $T: E \to E$ by setting $\alpha(T) := \alpha(TB_E)$. Observe, however, that in general it is not possible to obtain the weak essential norm $\norm{\cdot}_w$ as the value of such a set measure of weak noncompactness applied to operators in this way. In this sense, many of the measures of weak noncompactness appearing in the literature, when applied to operators, fall within the scope of Definition~\ref{def: measureweakcompactness}.

\begin{example}\label{ex: deBlasi}
    The De Blasi measure of weak noncompactness $\omega(\cdot)$, when applied to operators, is a normalised measure of weak noncompactness in the sense of Definition~\ref{def: measureweakcompactness}. The only non-trivial property to verify is \ref{it: en2}. Let $S,T:E\to E$ be bounded linear operators. Since $SB_E\subseteq \|S\| B_E$, we have
    \begin{equation*}
        \omega(TS)
        = \omega\bigl( T(SB_E) \bigr)
        \le \omega\bigl( T(\|S\| B_E) \bigr)
        = \|S\|\,\omega\bigl( TB_E \bigr)
        = \|S\|\,\omega(T).
    \end{equation*}
    Thus, it suffices to prove the general inequality
    \begin{equation*}
        \omega(TA) \le \|T\|\,\omega(A)
        \qquad\text{for every bounded set } A\subseteq E.
    \end{equation*}
    
    Fix $\varepsilon>\omega(A)$. By definition of $\omega(A)$ there exists a weakly compact set $D\subseteq E$ such that
    \begin{equation*}
        A \subseteq D + \varepsilon B_E.
    \end{equation*}
    Since bounded linear operators are weak-to-weak continuous, $T(D)$ is weakly compact. Applying $T$ to the inclusion above yields
    \begin{equation*}
        TA \subseteq T(D) + \varepsilon\,T(B_E)
             \subseteq T(D) + \varepsilon\,\|T\| B_E.
    \end{equation*}
    Thus
    \begin{equation*}
        \omega(TA)\le \varepsilon\|T\|.
    \end{equation*}
    As this holds for all $\varepsilon>\omega(A)$, we obtain
    \begin{equation*}
        \omega(TA) \le \|T\|\,\omega(A).
    \end{equation*}
    
    Applying this with $A = SB_E$, we conclude that
    \begin{equation*}
        \omega(TS)
        = \omega(TSB_E)
        \le \|T\|\,\omega(SB_E)
        = \|T\|\,\omega(S),
    \end{equation*}
    which verifies property \ref{it: en2} for the De Blasi measure of weak noncompactness.
\end{example}

\medskip

For a subset $A \subseteq E$ define
\begin{equation*}
   \operatorname{wk}_{E}(A) := \hat{\operatorname{d}}(\overline{A}^{w^{*}},E),
\end{equation*}
and
\begin{equation*}
    \operatorname{wck}_E(A) :=\sup\{\operatorname{d}(\operatorname{clust}_{E^{**}}(x_k), E) : (x_k) \text{ is sequence in } A\},
\end{equation*}
where
\begin{equation*}
    \operatorname{d}(A,B) = \inf \{\norm{a - b}: a \in A, b \in B \} \hspace{5pt} \text{ and } \hspace{5pt} \hat{\operatorname{d}}(A, B) = \sup\{d(a, B) : a \in A\},
\end{equation*}
so that $\operatorname{wk}_{E}(\cdot)$ and $\operatorname{wck}_E(\cdot)$ are measures of weak noncompactness for sets. Proposition \ref{prop: normalised-measures-weak-non} and Example \ref{ex: deBlasi}, together with \cite[Proposition 7.1, Lemma 7.4 and Theorem 7.5]{kacena2013quantitative}, yield the following.

\begin{corollary}\label{cor: relMeasuresOfNonCom}
Let $E = L_{1}(X,\Sigma,\mu)$ for a strictly localizable measure space and $T\in \mathscr{B}(E)$. Then we have
 \begin{align*}
     \norm{T}_w &= \inf_{\Theta\subseteq Y, \nu(\Theta)<\infty} \alpha_{\Theta}(T)=
        \omega(TB_E)= \operatorname{wk}_E(TB_E)= \operatorname{wck}_E(TB_E)\\
        &=\inf_{c>0, \mu(E)<\infty}\sup\{\int_{X}(|Tf| - c\mathds{1}_{E})^{+}d\mu(x):f \in B_E\}.
    \end{align*}
    In the case that $\mu$ is finite, this simplifies to
    \begin{align*}
        \norm{T}_w&=\Delta(T)=\omega(TB_E) = \operatorname{wk}_E(TB_E) = \operatorname{wck}_E(TB_E)\\
        &=\inf_{c > 0}\sup\{\int_{X}(|Tf|-c)^{+}d\mu(x): f \in B_E)\}.
    \end{align*}
\end{corollary}
\bigskip
\section{The Weak Calkin Algebra in $C(K)$-spaces}\label{sec: weak-calkinC(K)}

In this section, we will show how the results of the previous section allow us to provide a proof of Theorem \ref{th: C(K)-factorization}. Recall that for a compact Hausdorff space $K$, we have $C(K)^* = \mathcal{M}(K)$, the space of all Borel regular measures on $K$, which is an $AL$-space.  

We introduce the following measure of weak noncompactness in measure spaces, which originates from the seminal work of Grothendieck \cite{grothendieck1953applications}.
\begin{definition}
    Let $K$ be a compact Hausdorff space, and let $A \subseteq \mathcal{M}(K)$ be a bounded subset. The \emph{modulus of uniform regularity} of $A$ is defined by
    \begin{equation*}
        \varepsilon(A) = \sup \left\{ \inf_{n \in \mathbb{N}} |\mu_n|(U_n) : (U_n)_{n=1}^\infty \text{ disjoint open subsets of } K,\ (\mu_n)_{n=1}^\infty \subseteq A \right\}.
    \end{equation*}
\end{definition}

In particular, the classical characterisation of weak compactness of Grothendieck states that $A$ is relatively weakly compact if and only if $\varepsilon(A) = 0$. We have the following result.

\begin{lemma}\label{lmm: C(K)-1}
    Let $K$ be a compact Hausdorff space. Then for any operator $T: \mathcal{M}(K) \to \mathcal{M}(K)$ we have
    \begin{equation*}
        \varepsilon(TB_{\mathcal{M}(K)}) = \omega(T B_{\mathcal{M}(K)}) = \norm{T}_w.
    \end{equation*}
\end{lemma}
\begin{proof}
    The first equality is proved in the same way as the inequalities in \cite[Proposition~5.2]{kalenda2012quantification}. 
    Observe that the factor $1/\pi$ appearing in their result arises from taking the supremum over $|\mu_n(U_n)|$ instead of over $|\mu_n|(U_n)$. The second equality follows directly from Proposition~\ref{prop: normalised-measures-weak-non} and Example~\ref{ex: deBlasi}. 
\end{proof}

We now establish the first part of Theorem \ref{th: C(K)-factorization}, which, using the previous lemma, essentially amounts to a careful examination of the proof of Pe{\l}czy\'{n}ski's theorem.

\begin{proposition}\label{prop: C(K)-part1}
    Let $T: C(K) \to C(L)$ be an operator, where $K, L$ are compact Hausdorff spaces. Then
    \begin{equation*}
        \norm{T^*}_w = \sup \{m(TV)|\hspace{5pt}  V: c_0 \to C(L), \hspace{5pt} \norm{V} \leq 1 \}.
    \end{equation*}
\end{proposition}
\begin{proof}
    Assume first that $K = L$. For convenience, let
    \begin{equation*}
        \eta := \sup \{m(TV)|\hspace{5pt}  V: c_0 \to C(L), \hspace{5pt} \norm{V} \leq 1 \}.
    \end{equation*}
    We start by showing that $\eta \leq \norm{T^*}_w$. Fix $\varepsilon > 0$, by definition we can find $V: c_0 \to C(K)$ such that $TV: c_0 \to C(K)$ is an isomorphism with $m(TV) \geq \eta - \varepsilon$. It follows that $V^*T^*: C(L)^* \to \ell_1$ is a quotient map such that $(\eta - \varepsilon)B_{\ell_1} \subseteq V^*T^*[B_{\mathcal{M}(K)}]$. Consequently, for any weakly compact operator $W$, we have 
    \begin{align*}
        \norm{T^* - W} & \geq \norm{V^*}\norm{T^* - W}  \geq \norm{V^*T^*  - V^*W} \geq \omega(V^*T^*  - V^*W) \\
        & = \omega(V^*T^*) \geq \eta - \varepsilon,
    \end{align*}
    where we recall that $\omega(\cdot)$ is the De Blasi measure of weak noncompactness \eqref{eq: deBlasi}. Taking the infimum over all weakly compact operators gives $\norm{T^*}_w \geq \eta - \varepsilon$, since $\varepsilon$ was arbitrary, we obtain $\norm{T^*}_w \geq \eta$.

    We focus now on the reverse implication. By Lemma \ref{lmm: C(K)-1} we have \break $\varepsilon(T^*) := \varepsilon(T^*  B_{\mathcal{M}(K)}) = \norm{T^*}_w$ for all $T \in \mathscr{B}(C(K))$, and thus it is enough to show that $\eta \geq \varepsilon(T^*)$. Fix $\varepsilon > 0$. By definition of $\varepsilon(T^*)$ we can find a sequence of disjoint open sets $(U_n)_{n =1}^\infty$ and measures $(\nu_n)_{n=1}^\infty \subseteq B_{\mathcal{M}(K)}$ such that $|T^* \nu_n|(U_n) > \varepsilon(T^*) - \varepsilon$ for all $n \in \mathbb{N}$. To simplify the notation, we set $\mu_n := T^* \nu_n$.
    
    By internal regularity, for each $n \in \mathbb{N}$, we can find a compact subset $F_n \subseteq U_n$ such that $| \mu_n|(U_n \backslash F_n) < \varepsilon$, which implies $|\mu_n|(F_n) > \varepsilon(T^*) - 2\varepsilon$.
    
    Observe that
    \begin{equation*}
        |\mu_n|(F_n) = \sup \{ |\int_{F_n} g d\mu_n|: g \in B_{C(K)}\},
    \end{equation*}
    so that, for each $n \in \mathbb{N}$, we can find $g_n \in B_{C(K)}$ such that $\int_{F_n} g_n d\mu_n > \varepsilon(T^*) - 2\varepsilon$. By Urysohn's lemma, there exists $h_n \in C(K)$, $0 \leq h_n \leq 1$ such that $h_n = 1$ in $F_n$ and $h_n = 0$ in $K \backslash U_n$. Let $f_n := h_n g_n$ and note that
    \begin{equation*}
        |\mu_n(f_n)| = |\int_{K} f_n d\mu_n| = |\int_{U_n} h_n g_n d\mu_n| \geq \int_{F_n} g_n d\mu_n - \int_{U_n \backslash F_n} |h_n g_n|d\mu_n > \varepsilon(T^*) - 3\varepsilon.
    \end{equation*}
    Define the map $V: c_0 \to C(K)$ by $Ve_n = f_n$, where $(e_n)_{n=1}^\infty$ is the canonical basis of $c_0$. Observe that, by construction, $\norm{V} \leq 1$.
    We have that $TV: c_0 \to C(K)$ satisfies
    \begin{align*}
        \inf_{n \in \mathbb{N}} \norm{TVe_n} &\geq \inf_{n \in \mathbb{N}}|\nu_n(TVe_n)| = \inf_{n \in \mathbb{N}}  |(T^*\nu_n)(f_n)| = \inf_{n \in \mathbb{N}}  |\mu_n(f_n)| > \varepsilon(T^*) - 3\varepsilon.
    \end{align*}
    It follows from a result of Rosenthal (see \cite{rosenthal1970relatively}, Remark 1 after Theorem 3.4) that $TV$ is an isomorphism on a copy of $c_0$. Moreover, as noted in \cite[Theorem 3.3]{arnott2023uniqueness}, we can ensure that, on this copy, $TV$ is bounded below by $\varepsilon(T^*) - 3\varepsilon$. In other words, by restricting the domain of $V$ to this copy, we obtain $m(TV) > \varepsilon(T^*) - 3\varepsilon$.
    Since $\varepsilon > 0$ was arbitrary, we conclude that $\eta \geq \varepsilon(T^*) = \norm{T^*}_w$, which completes the proof in the case $K = L$.

    For arbitrary compact spaces $K \not = L$, note that $C(K) \oplus C(L) = C(K \sqcup L)$, where the sum is taken with the $\ell_\infty$-norm and $K \sqcup L$ denotes the topological disjoint union of $K$ and $L$. Naturally, we may view $T: C(K) \to C(L)$ as an operator $T: C(K \sqcup L) \to C(K \sqcup L)$, so the previous case and straightforward manipulations yield the result.
\end{proof}

\begin{rem}
    The quantity
    \begin{equation*}
        \sup\{ m(TV) \;|\; V : c_0 \to C(L),\; \|V\|\le 1 \},
    \end{equation*}
    which measures how well an operator $T$ fixes a copy of $c_0$, has also been studied for operators between arbitrary Banach spaces $T : X \to Y$. For example, in \cite{krulivsova2017quantification} this quantity appears under the name $\operatorname{fix}_{c_0}(T)$, and it is shown that
    \begin{equation*}
    \operatorname{fix}_{c_0}(T) = \operatorname{uc}(T)
    \end{equation*}
    \cite[Theorem~3.6]{krulivsova2017quantification}, where $\operatorname{uc}(T)$ quantifies how far $T$ is from being unconditionally converging. 
    Moreover, it is shown that for operators $T : C(K)\to Y$, this quantity coincides with several other measures of weak noncompactness, see \cite[Corollary~5.6]{krulivsova2017quantification}.
    In particular, Proposition~\ref{prop: C(K)-part1} may be viewed as an extension of the aforementioned corollary in the case where the codomain is also a $C(K)$-space.
\end{rem}

In the previous proposition, we showed that the weak essential norm of the dual operator $T^*$ provides a measure of how well an operator $T$ fixes a copy of $c_0$. To obtain a factorization, we need to show that this copy is complemented with good constants (possibly after passing to a subcopy of $c_0$ inside of it). 

We will establish this in two general settings. First, when $K$ is sequentially compact (for instance, if $K$ is first countable), we will argue that every copy of $c_0$ admits a further subcopy that is complemented with good constants. Second, when $L$ is extremally disconnected, fixing a copy of $c_0$ implies fixing a copy of $\ell_\infty$ with the same constants. \\

We note that some structural assumptions on $K$ and $L$ are necessary. It is possible for an operator $T: C(K) \to C(L)$ to fix a copy of $c_0$ but not fix a complemented copy of $c_0$ nor a copy of $\ell_\infty$. For instance, by a remarkable construction of Haydon~\cite{haydon1981non}, there exists a compact Hausdorff space $K$ such that $C(K)$ is a Grothendieck space (hence it contains no complemented copy of $c_0$, see~\cite[Proposition 3.1.13]{gonzalez2021grothendieck}) not containing a copy of $\ell_\infty$. In this case, no factorization of the identity of $c_0$ or $\ell_\infty$ is possible through the identity operator $I_{C(K)}: C(K) \to C(K)$, while $I_{C(K)}$ trivially fixes a copy of $c_0$.

We begin with the case $K$ sequentially compact. For this, we require a mild variation of \cite[Theorem 6]{dowling1999remarks} that is better suited to our context. If $E$ is a Banach space, a subset $M \subseteq B_{E^*}$ is said to be $1$-norming if
\begin{equation*}
    \norm{x} = \sup_{x^* \in M} |x^*(x)|.
\end{equation*}

\begin{theorem}[Dowling, Randrianantoanina and Turett]\label{th: Downling-slighly-general}
    Let $E$ be a Banach space. Suppose there exists $M \subseteq B_{E^*}$ $1$-norming and weak$^*$ sequentially compact. Then for any $Y \subseteq E$, $Y \sim c_0$ and any $\varepsilon > 0$, there exist $Z \subseteq Y$ and a projection $P: E \to Z$ such that $Z$ is $(1 + \varepsilon)$-isometric to $c_0$ and $\norm{P} \leq 1 + \varepsilon$.
\end{theorem}
\begin{proof}
    The proof is essentially that of \cite[Theorem 6]{dowling1999remarks}, with natural modifications to accommodate the $1$-norming set $M$. We include an outline for the reader's convenience. Let $\delta > 0$ be chosen sufficiently small, depending on the desired $\varepsilon$. By James's distortion theorem (see, for example, \cite[Theorem 1]{diestel1984sequences}), there exists a sequence $(u_n)_{n = 1}^\infty \subseteq B_Y$ such that
    \begin{equation*}
        (1 - \delta) \sup_{n \in \mathbb{N}}|a_n| \leq \norm{\sum_{n=1}^\infty a_n u_n} \leq  \sup_{n \in \mathbb{N}}|a_n|,
    \end{equation*}
    for any $(a_n)_{n=1}^\infty \in c_0$. Choose $(u^*_n)_{n=1}^\infty \subseteq M$ such that $|u_n^*(u_n)| \geq 1 - \delta$. Multiplying $u_n$ by an unimodular scalar number, we may assume $u_n^*(u_n) \geq 1 - \delta$.  Observe that this implies
    \begin{equation}\label{eq: technicalc0}
        |u^*_n (\sum_{m \not = n} a_m u_m ) | \leq \delta
    \end{equation}
    for all $(a_m)_{m \in \mathbb{N}} \in B_{c_0}$ and all $n \in \mathbb{N}$. Since $M$ is weak$^*$ sequentially compact, by passing to a subsequence, we may assume $(u_n^*)_{n=1}^\infty$ converges weak$^*$ to some element.  Define $w^*_n =  (u_{2n} - u_{2n-1})/2$ and $w_n = u_{2n} - u_{2n-1}$. It is clear that $(w^*_n)_{n=1}^\infty$ is a weak$^*$ null sequence, while $w_n^*(w_n) \geq 1 - \delta$ for all $n \in \mathbb{N}$ and \eqref{eq: technicalc0} still holds with $u_n^*$ and $u_m$ replaced by $w_n^*$ and $w_m$ respectively.

    Define the operator $J: c_0 \to E$, $(a_n)_{n=1}^\infty \mapsto \sum_{n=1}^\infty a_n w_n$ and $Q: E \to c_0$, $x \mapsto (w^*_n(x))_{n=1}^\infty$. It is straightforward to verify that $\norm{I_{c_0} - QJ} \leq 2\delta$, and thus we can find an operator $R: c_0 \to c_0$ with $\norm{R} \leq (1 - 2\delta)^{-1}$ such that $RQJ = I_{c_0}$. It is easy to verify that $Z = \overline{\spn \{w_n: n \in \mathbb{N} \}} \subseteq Y$ is the desired copy of $c_0$ and $P = J(QR)$ the desired projection.
\end{proof}

\begin{rem}\label{rem: Ant-avil}
    The previous reformulation is stronger than the original version by Dowling, Randrianantoanina, and Turett, particularly in the context of $C(K)$-spaces. Indeed, the set $\{\delta_k : k \in K\} \subseteq C(K)^*$ is always $1$-norming, so the existence of a $1$-norming, weak$^*$ sequentially compact subset of the dual unit ball is automatically satisfied whenever $K$ is sequentially compact. In contrast, requiring the dual unit ball to be weak$^*$ sequentially compact can be a stronger condition. For instance, assuming the Continuum Hypothesis, Talagrand \cite{talagrand1980mesures} constructed a Corson compact space $K$ that is first countable, and hence sequentially compact, but for which the unit ball of $C(K)^*$ (endowed with the weak$^*$ topology) contains a topological copy of $\beta \mathbb{N}$ (see \cite[Theorem 5.6]{plebanek2024survey}). In particular, the dual unit ball fails to be sequentially compact while $K$ is sequentially compact. We thank Antonio Avil\'es for bringing this example to our attention.
\end{rem}

The case when $L$ is extremally disconnected will follow from a quantitative version of a result of Rosenthal, originally stated as \cite[Theorem 1.3]{rosenthal1970relatively}.

\begin{theorem}[Rosenthal]\label{th: Roshental}
    Let $\Gamma$ be an infinite set, and let $E$ and $F$ be Banach spaces with $E$ $1$-complemented in $E^{**}$. Let $T: E \to F$ and $V: c_0 \to E$, $\norm{V} \leq 1$, such that $m(TV) > \delta$. Then there exists $\tilde{V}: \ell_\infty \to E$ such that $T\tilde{V}: \ell_\infty \to F$ satisfies $m(T\tilde{V}) > \delta$.
\end{theorem}
\begin{proof}
    It suffices to observe that the proof of \cite[Proposition 1.2]{rosenthal1970relatively} is essentially quantitative. This proposition is the key ingredient in the proof of \cite[Theorem 1.3]{rosenthal1970relatively}, which now becomes quantitative upon incorporating the condition that the projection $P: E^{**} \to E$ has norm one.
\end{proof}

As an automatic consequence, we get the stated factorization result for operators between $C(K)$-spaces.

\begin{proof}[Proof of Theorem \ref{th: C(K)-factorization}]
    The first part was proven in Proposition \ref{prop: C(K)-part1}. The second part follows automatically from the first part, which tells us that we fix a copy of $c_0$ with controlled bounds, and Theorems \ref{th: Downling-slighly-general} and \ref{th: Roshental}. 

    Indeed, Theorem \ref{th: Downling-slighly-general} guarantees that when $K$ is sequentially compact, we can pass to a subcopy of the fixed copy of $c_0$ that is complemented with good constants. On the other hand, when $L$ is extremally disconnected, $C(L)$ is $1$-complemented in its bidual; thus, by Theorem \ref{th: Roshental}, we quantitatively fix a copy of $\ell_\infty$, which is $1$-injective and therefore complemented with good constants.

    The details are standard and hence omitted here; see, for example, the argument in \cite[Lemma 3.7]{arnott2023uniqueness}.
\end{proof}

Corollaries \ref{cor: algebra-normC(K)-space-extremally}, \ref{cor: algebra-normC(K)-space-sequentially}, \ref{cor: BonsalMinimality-extremally}, and \ref{cor: BonsalMinimality-sequentially} follow directly from Theorem \ref{th: C(K)-factorization}. In both cases, minimality is derived from the argument in \cite[Theorem 1.1]{acuaviva2024factorizations} -- see also \cite[Remark 8.3]{acuaviva2024factorizations}--, and observe that, in contrast with the proof of Corollary \ref{cor: unique-algebra-normAL-space}, we need these arguments to bypass the need for maximality. Bonsall's minimality property is established by the same reasoning used in the proof of \cite[Theorem 1.5]{acuaviva2024factorizations}. Alternatively, similar results regarding measures of weak noncompactness, akin to those in Proposition \ref{prop: normalised-measures-weak-non}, can be obtained, which also imply Bonsall's minimality property. Since these arguments closely follow the aforementioned ones, we omit the details. Theorem \ref{th: Residum-and-w} follows immediately from these results. 

\bigskip
\newline
\bigskip

\noindent\textbf{Acknowledgements.} This paper forms part of the first-named author's PhD research at Lancaster University, conducted under the supervision of Professor N. J. Laustsen. The authors are grateful to Professor Laustsen for his helpful comments and suggestions on the manuscript presentation. We are also thankful to W. B. Johnson for sharing an early draft of the manuscript \cite{johnsonandphillips} and to Antonio Avil\'es for making us aware of the example in Remark \ref{rem: Ant-avil}. Finally, we are grateful to the anonymous referee for a detailed and helpful report, and in particular for drawing our attention to the notion of rich families, which helped streamline some of our results.

The first-named author gratefully acknowledges funding from the EPSRC (grant number EP/W524438/1), which has supported his studies.\\

\noindent\textbf{Data availability.} No data was used for the research described in the article.


\end{document}